   \pdfoutput=1 

\documentclass{article}
\usepackage{excludeonly}

\usepackage{savesym}
\usepackage{scrextend}  
\usepackage{amsthm,amsmath}
\usepackage{amssymb,latexsym,graphicx}
\usepackage{amscd}
 \usepackage[all,cmtip]{xy}
\usepackage{tikz-cd}

\usepackage{accents}
\usepackage{cite}
\usepackage{mathtools}
\usepackage{stmaryrd} 

\usepackage{calligra}
\DeclareMathAlphabet{\mathcalligra}{T1}{calligra}{m}{n}
\DeclareMathAlphabet{\mathpzc}{OT1}{pzc}{m}{it}
\usepackage{hycolor}
\usepackage{xcolor}
\usepackage[
                       colorlinks=true,
                       linkcolor=black, 
                       citecolor=black, 
                       urlcolor=blue,
                     ]{hyperref}
\usepackage{soul} 

\usepackage{makeidx}

\usepackage[intoc]{nomentbl}
\makenomenclature

\usepackage{stackengine} 
\usepackage{booktabs} 

\newtheorem{theoremABC}{Theorem}

\newtheorem{theorem}{Theorem}[section]

\newtheorem{corollary}[theorem]{Corollary}

\newtheorem{lemma}[theorem]{Lemma}
\newtheorem{proposition}[theorem]{Proposition}

\theoremstyle{definition}
\newtheorem{definition}[theorem]{Definition}

\newtheorem{remark}[theorem]{Remark}

\newtheorem{example}[theorem]{Example}

\newtheorem{question}[theorem]{Question}

\theoremstyle{remark}

\newcommand{\A}{{\mathbb{A}}}
\newcommand{\B}{{\mathbb{B}}}
\newcommand{\C}{{\mathbb{C}}}

\newcommand{\N}{{\mathbb{N}}}

\newcommand{\R}{{\mathbb{R}}}

\newcommand{\Z}{{\mathbb{Z}}}
\newcommand{\Dd}{{\mathcal{D}}}
\newcommand{\Ee}{{\mathcal{E}}}
\newcommand{\Ff}{{\mathcal{F}}}
\newcommand{\Hh}{{\mathcal{H}}}
\newcommand{\Ii}{{\mathcal{I}}}

\newcommand{\Ll}{{\mathcal{L}}}   

\newcommand{\Ss}{{\mathcal{S}}}
\newcommand{\Tt}{{\mathcal{T}}}
\newcommand{\Uu}{{\mathcal{U}}}
\newcommand{\Vv}{{\mathcal{V}}}
\newcommand{\Ww}{{\mathcal{W}}}

\newcommand{\coker}{{\rm coker\, }}  
\newcommand{\im}{{\rm im\, }}             
\newcommand{\Id}{{\rm Id}}
\newcommand{\INDEX}{\mathop{\mathrm{index}}}     

\newcommand{\cgraph}[1]{\Gamma_{\kern-.5ex{}#1}}     
\renewcommand{\Re}{{\rm Re}}       
\renewcommand{\Im}{{\rm Im}}       

\newcommand{\Map}{{\rm Map}}          
\newcommand{\mat}[1]{\mathbf{#1}}  
\newcommand{\Fix}{{\rm Fix}}            
\newcommand{\spec}{\mathrm{spec}\,}        
\newcommand{\Arnold}{{Arnol$'$d}}           
\newcommand{\norm}{{\rm norm}}

\newcommand{\eps}{{\varepsilon}}

\newcommand{\inner}[2]{\langle #1, #2\rangle}   
\newcommand{\INNER}[2]{\left\langle #1, #2\right\rangle}

\newcommand{\SC}{{\mathrm{sc}}}                  
\newcommand{\mbf}[1]{\text{\boldmath $#1$}}  

\def\NABLA#1{{\mathop{\nabla\kern-.5ex\lower1ex\hbox{$#1$}}}}
\def\Nabla#1{\nabla\kern-.5ex{}_{#1}}
\def\Tabla#1{\Tilde\nabla\kern-.5ex{}_{#1}}
\def\Babla#1{\widebar\nabla\kern-.5ex{}_{#1}}
\def\abs#1{\mathopen|#1\mathclose|}   
\def\Abs#1{\left|#1\right|}            
\def\norm#1{\mathopen\|#1\mathclose\|}
\def\Norm#1{\left\|#1\right\|}

\renewcommand{\Tilde}{\widetilde}

\newcommand{\quarter}{{\frac{1}{4}}}

\newcommand{\INTO}{\hookrightarrow}              

\newlength\eqshift
\renewcommand\theequation{\thesection.\arabic{equation}}
\let\savetheequation\theequation

\makeatletter
\renewcommand*\env@matrix[1][\arraystretch]{%
  \edef\arraystretch{#1}%
  \hskip -\arraycolsep
  \let\@ifnextchar\new@ifnextchar
  \array{*\c@MaxMatrixCols c}}
\makeatother

\makeatletter
\let\save@mathaccent\mathaccent
\newcommand*\if@single[3]{%
  \setbox0\hbox{${\mathaccent"0362{#1}}^H$}%
  \setbox2\hbox{${\mathaccent"0362{\kern0pt#1}}^H$}%
  \ifdim\ht0=\ht2 #3\else #2\fi
  }
\newcommand*\rel@kern[1]{\kern#1\dimexpr\macc@kerna}
\newcommand*\widebar[1]{\@ifnextchar^{{\wide@bar{#1}{0}}}{\wide@bar{#1}{1}}}
\newcommand*\wide@bar[2]{\if@single{#1}{\wide@bar@{#1}{#2}{1}}{\wide@bar@{#1}{#2}{2}}}
\newcommand*\wide@bar@[3]{%
  \begingroup
  \def\mathaccent##1##2{%
    \let\mathaccent\save@mathaccent
    \if#32 \let\macc@nucleus\first@char \fi
    \setbox\z@\hbox{$\macc@style{\macc@nucleus}_{}$}%
    \setbox\tw@\hbox{$\macc@style{\macc@nucleus}{}_{}$}%
    \dimen@\wd\tw@
    \advance\dimen@-\wd\z@
    \divide\dimen@ 3
    \@tempdima\wd\tw@
    \advance\@tempdima-\scriptspace
    \divide\@tempdima 10
    \advance\dimen@-\@tempdima
    \ifdim\dimen@>\z@ \dimen@0pt\fi
    \rel@kern{0.6}\kern-\dimen@
    \if#31
      \overline{\rel@kern{-0.6}\kern\dimen@\macc@nucleus\rel@kern{0.4}\kern\dimen@}%
      \advance\dimen@0.4\dimexpr\macc@kerna
      \let\final@kern#2%
      \ifdim\dimen@<\z@ \let\final@kern1\fi
      \if\final@kern1 \kern-\dimen@\fi
    \else
      \overline{\rel@kern{-0.6}\kern\dimen@#1}%
    \fi
  }%
  \macc@depth\@ne
  \let\math@bgroup\@empty \let\math@egroup\macc@set@skewchar
  \mathsurround\z@ \frozen@everymath{\mathgroup\macc@group\relax}%
  \macc@set@skewchar\relax
  \let\mathaccentV\macc@nested@a
  \if#31
    \macc@nested@a\relax111{#1}%
  \else
    \def\gobble@till@marker##1\endmarker{}%
    \futurelet\first@char\gobble@till@marker#1\endmarker
    \ifcat\noexpand\first@char A\else
      \def\first@char{}%
    \fi
    \macc@nested@a\relax111{\first@char}%
  \fi
  \endgroup
}
\makeatother

\long\def\symbolfootnote[#1]#2{\begingroup%
\def\thefootnote{\fnsymbol{footnote}}\footnote[#1]{#2}\endgroup}

\tikzset{
  symbol/.style={
    draw=none,
    every to/.append style={
      edge node={node [sloped, allow upside down, auto=false]{$#1$}}}
  }
}

\newcommand{\gf}[1]{#1}                       

\newcommand{\gfset}{\mathcal{G}}         
\newcommand{\gtset}{\mathfrak{G}   }    

\newcommand{\gt}[1]{\mathfrak #1}      

\newcommand{\HhMorse}{\mathcal{H}^{\text{Morse}}}          
\newcommand{\HhcoMorse}{\mathcal{H}^{\text{co-Morse}}} 

\begin{document}
\sloppy
\author{\quad Urs Frauenfelder \quad \qquad\qquad
             Joa Weber\footnote{
  Email: urs.frauenfelder@math.uni-augsburg.de
  \hfill
  joa@unicamp.br
  }
    \\
    Universit\"at Augsburg \qquad\qquad
    UNICAMP
}

\title{Growth of eigenvalues of Floer Hessians}

\date{\today}

\maketitle 
%


%
%

%





\begin{abstract}
In this article we prove that the space of Floer Hessians
has infinitely many connected components.
\end{abstract}

\tableofcontents

\section{Introduction}

\subsection{Main results}

We consider a Hilbert space pair $H:=(H_0,H_1)$,
that is $H_0$ and $H_1$ are both infinite dimensional
Hilbert spaces such that $H_1\subset H_0$
and inclusion $\iota\colon H_1\INTO H_0$ is compact and dense.
Then there exists an unbounded monotone function
$$
   \gf{h}=\gf{h}(H_0,H_1) \colon \mathbb{N} \to (0,\infty),
$$
called \textbf{pair growth function},
such that the pair $(H_0,H_1)$ is isometric to the pair
$(\ell^2,\ell^2_{\gf{h}})$,
see Appendix~\ref{sec:HS-pairs}, and from now on we identify the pairs
$$
   (H_0,H_1)=(\ell^2, \ell^2_\gf{h})
   ,\qquad
   \gf{h}=\gf{h}(H_0,H_1) .
$$
Here $\ell^2_\gf{h}$ is defined as follows. In general,
for $\gf{f}\colon\mathbb{N}\to(0,\infty)$ unbounded monotone
$\ell^2_\gf{f}$ is the space of all sequences $x=(x_\nu)_{\nu \in \mathbb{N}}$ with
$\sum_{\nu=1}^\infty \gf{f}(\nu)x_\nu^2 <\infty$.
The space $\ell^2_\gf{f}$ becomes a Hilbert space if we endow it
with the inner product
\begin{equation*}
   \langle x,y \rangle_\gf{f}=\sum_{\nu=1}^\infty \gf{f}(\nu)x_\nu y_\nu
   ,\qquad
   \Norm{x}_\gf{f}:=\sqrt{\langle x,x \rangle_\gf{f}} .
\end{equation*}
Therefore we can define for every real number $r$ a Hilbert space 
$
   H_r:=\ell^2_{\gf{f}^r}
$.
For every $s<r$ the inclusion $H_r\subset H_s$
is compact and dense; see Section~\ref{sec:frac}.

\begin{definition}[Weak Hessian]\label{def:weak-Hessian}
The topological\footnote{
  $\Ll(H_1,H_0)\cap \Ll(H_2,H_1)$ is a Banach space under the norm
  $\max\{\norm{\cdot}_{\Ll(H_1,H_0)},\norm{\cdot}_{\Ll(H_2,H_1)}\}$.
  }
\textbf{space of weak Hessians} 
$$
   \Hh_\gf{h}\subset\Ll(H_1,H_0)\cap \Ll(H_2,H_1)
$$
consists of all bounded linear maps 
$A\colon H_1\to H_0$ which are Fredholm\footnote{
  Closed image and finite dimensional kernel and cokernel,
  $\INDEX :=\dim\ker -\dim\coker$.
  }
of index zero and \textbf{\boldmath$H_0$-symmetric}, that is
\begin{equation}\label{eq:0-symmetric}
   \INNER{Ax}{y}_0=\INNER{x}{Ay}_0
\end{equation}
for all $x,y\in H_1$. Furthermore, it is required that $A$
restricts to a bounded linear map $H_2\to H_1$,
notation $A_1$, which is Fredholm of index zero as well.
\end{definition}

Since the inclusion $H_1=\ell^2_\gf{h}\INTO H_0=\ell^2$
is compact, the resolvent $R\colon\ell^2\to\ell^2$ of a weak Hessian $A$ is a compact
operator. Hence there exists an orthonormal basis $\Vv$ of $\ell^2$
composed of eigenvectors $v_\nu\in \ell^2_\gf{h}$ of
$A\colon \ell^2_\gf{h}\to\ell^2$ and the spectrum of $A$
consists of discrete real eigenvalues $a_\nu$ of finite multiplicity each.
Arrange the squares of eigenvalues of $A$ with multiplicity in an increasing
sequence $\{a_\nu^2\}_{\nu\in\N}$ and order the orthonormal eigenvector
basis  $\Vv=(v_\nu)_{\nu\in\N}\subset \ell^2_\gf{h}$ of $\ell^2$ accordingly.
\\
Assume that $A$ is invertible. Then the squares of all eigenvalues of
$A$ are positive and
we define a positive monotone increasing function, 
called \textbf{total growth of a weak Hessian} $\mbf{A}\in\Hh_\gf{h}$,
by
\begin{equation}\label{eq:f_A}
   \gf{f}_A\colon \N\to (0,\infty)
   ,\quad
   \gf{f}_A(\nu)
   :=a_\nu^2 .
\end{equation}

\begin{lemma}\label{le:f_a-f}
Any invertible weak Hessian $A\in\Hh_\gf{h}$ has a total growth function
$\gf{f}_A$ of the same \textbf{growth type} as $\gf{h}$,
notation~${\gf{f}_A\sim \gf{h}}$.\footnote{
  That is, there exists a constant $c>0$
  such that $\frac{1}{c}\cdot\gf{h}(\nu)
  \le \gf{f}_A(\nu)\le c\cdot\gf{h}(\nu)$,
  $\forall \nu\in\N$.
  }
\end{lemma}

\begin{proof}
By Theorem~\ref{thm:f-g-equiv} 
we need to find an isomorphism of Hilbert space pairs
$$
   T\colon 
   (\ell^2,\ell^2_{\gf{f}_{A}})
   \to
   (\ell^2,\ell^2_{\gf{h}}) .
$$
By assumption $A\colon \ell^2_{\gf{h}}\to\ell^2$ is an isomorphism.
Hence the pull-back inner product $A^*\INNER{\cdot}{\cdot}_{\ell^2}$
is on $\ell^2_{\gf{h}}$ an inner product equivalent to the inner
product $\INNER{\cdot}{\cdot}_{\gf{h}}$.
Therefore the identity gives rise to an isomorphism
$$
   I\colon
   (\ell^2,(\ell^2_{\gf{h}},A^*\INNER{\cdot}{\cdot}_{\ell^2}))
   \to
   (\ell^2,\ell^2_{\gf{h}})
$$
of Hilbert pairs.

Let $(e_\nu)_{\nu\in\N}$ be the standard basis of $\ell^2$ and $(v_\nu)_{\nu\in\N}$
the orthonormal eigenvector basis introduced above
($Av_\nu=a_\nu v_\nu$).
Consider the isometry given by
$$
   \Psi\colon\ell^2\to\ell^2,\quad e_\nu\mapsto  v_\nu .
$$
For $\mu,\nu\in\N$ we compute
$$
   \INNER{A\Psi e_\mu}{A\Psi e_\nu}_{\ell^2}
   =\INNER{A v_\mu}{A v_\nu}_{\ell^2}
   =a_\mu a_\nu \delta_{\mu \nu}
   =a_\mu^2 \delta_{\mu \nu}
   \stackrel{\text{(\ref{eq:f_A})} \atop \text{(\ref{eq:ell^2_f})}}{=}
   \INNER{e_\nu}{e_\mu}_{\gf{f}_A}.
$$
This
shows that $\Psi^* A^*\INNER{}{}_{\ell^2}=\INNER{\cdot}{\cdot}_{\gf{f}_A}$.
Therefore the restriction
$$
   \Psi|_{\ell^2_{\gf{f}_A}}\colon 
   \ell^2_{\gf{f}_A}
   \to
   (\ell^2_{\gf{h}},A^*\INNER{\cdot}{\cdot}_{\ell^2})
$$
is an isometry.
To summarize we have an isometry of pairs
$$
   \Psi\colon
   (\ell^2,\ell^2_{\gf{f}_{A}})
   \to
   \left(\ell^2, (\ell^2_{\gf{h}},A^*\INNER{\cdot}{\cdot}_{\ell^2})\right) .
$$
Hence $T:=I\Psi\colon (\ell^2,\ell^2_{\gf{f}_{A}})\to(\ell^2,\ell^2_{\gf{h}})$
is an isomorphism of Hilbert pairs.
\end{proof}

According to the lemma the total growth type of $A$
is completely determined by the ambient Hilbert spaces.
Hence it cannot be used to distinguish different connected components of
the space $\Hh_\gf{h}$.
\\
In order to get some information on the topology of the space
$\Hh_\gf{h}$, in particular its connected components,
we are taking into account as well the signs of the eigenvalues.
We distinguish the following three \textbf{types}.
\begin{enumerate}
\item
  \textbf{Morse.} Finitely many negative, infinitely many positive
  eigenvalues.
\item
  \textbf{Co-Morse.} Finitely many positive, infinitely many negative
  eigenvalues.

\item
  \textbf{Floer.} Infinitely many negative and positive 
  eigenvalues.
\end{enumerate}

An application of spectral flow~\cite{robbin:1995a}
shows the following.

\begin{lemma}
Suppose that $A^0$ and $A^1$ are connected by a continuous path in
$\Hh_\gf{h}$.
Then $A^0$ and $A^1$ are either both Morse, both Co-Morse, or both Floer.
\end{lemma}

According to the above three types the set of weak Hessians $\Hh_\gf{h}$
decomposes into three open and closed subsets denoted by
$$
   \Hh_\gf{h}
   =\HhMorse_\gf{h}
   \cup\HhcoMorse_\gf{h}
   \cup \mathcal{H}^{\text{Floer}}_\gf{h} .
$$
In this article we are interested in the Floer Hessians.
Therefore we abbreviate
$$
   \Ff_\gf{h}:=\mathcal{H}^{\mathrm{Floer}}_\gf{h}
   ,\qquad
   \Ff_\gf{h}^*:=\{ \text{$A\in \Ff_\gf{h}$ invertible}\} .
$$
In the Floer case, if $A\in\Ff_\gf{h}$ is invertible, i.e. zero is not an
eigenvalue, we look separately at the \textbf{positive} and
\textbf{negative growth functions}
$\gf{f}_A^+$ and $\gf{f}_A^-$---given by ordering the
positive eigenvalues of the weak Hessian $A$, with multiplicities,
respectively the \emph{absolute values} of the 
negative eigenvalues.

\begin{definition}
A \textbf{growth function} is a monotone unbounded function
$\gf{f} \colon \mathbb{N} \to (0,\infty)$. Having the same growth
type, symbol $\gf{f}\sim\gf{g}$, see previous footnote,
defines an equivalence relation on
the set $\gfset$ of all growth functions. The elements
$\gt{f}=[\gf{f}]$ of the quotient $\gtset:=\gfset/\sim$ are called
\textbf{growth types}.
\end{definition}

We are interested in weak Hessians of given negative and positive growth types.
That is, given $\gt{a},\gt{b}\in\gtset$, we are interested in the sets
\begin{equation}\label{eq:F^*}
   (\Ff_\gf{h}^{\gt{a}\gt{b}})^*
   :=\{A\in\Ff_\gf{h}^*
  \mid \text{$[\gf{f}_A^-]=\gt{a}$, $[\gf{f}_A^+]=\gt{b}$}\} .
\end{equation}
Given a growth function $\gf{f} \colon \mathbb{N} \to
(0,\infty)$, we define the shifted growth function
$\gf{f}_1(\nu):=\gf{f}(\nu+1)$.
We say that the growth type is \textbf{shift invariant}
if $\gf{f}\sim \gf{f}_1$.
We denote the \textbf{set of shift-invariant growth types} by
$\gtset_1\subset\gtset$.

\smallskip
For $\lambda\in\R$ consider the translation
\begin{equation}\label{eq:AA_lambda}
   T_\lambda\colon \Ff_\gf{h}\to \Ff_\gf{h}
   ,\quad
   A\mapsto A-\lambda\iota=:\A_\lambda
\end{equation}
which is a homeomorphism with inverse $T_{-\lambda}$.
Observe that $\A_\lambda\colon H_1\to H_0$ is in fact a weak Floer
Hessian operator. Indeed the inclusion $\iota\colon H_1\INTO H_0$ is
compact and $H_0$-symmetric and so is $-\lambda\iota$.
But the sum of a Fredholm operator $A$
and a compact operator $-\lambda\iota$
is Fredholm of the same index.
Moreover, the sum of $H_0$-symmetric operators is $H_0$-symmetric.
The inclusion restricts to a map $\iota|\colon H_2\INTO H_1$
and the same arguments show that $A-\lambda\iota\colon H_2\to H_1$
is Fredholm of index zero and $H_0$-symmetric.
This shows that $A-\lambda\iota \in \Ff_\gf{h}$.

\begin{proposition}\label{prop:cascasca-intro-NEW}
Let $\gt{a},\gt{b}\in\gtset_1$ be shift-invariant growth types.
If $A$ is a weak Hessian in $(\Ff_\gf{h}^{\gt{a}\gt{b}})^*$,
then $\forall\lambda\in\R\setminus\spec A$ so is the shifted operator
$
   \A_\lambda
   \in (\Ff_\gf{h}^{\gt{a}\gt{b}})^*
$.
\end{proposition}

\begin{proof}
Section~\ref{sec:Hessian-growth-types}.
\end{proof}

If $\gt{a},\gt{b}\in\gtset_1$ are \emph{shift-invariant} growth types we define
$$
   \Ff_\gf{h}^{\gt{a}\gt{b}}
   :=\{
   A\in \Ff_\gf{h} \mid
   \forall \lambda\in\R\setminus\spec A\colon
   \text{$[\gf{f}^+_{\A_\lambda}]=\gt{a}$ $\wedge$
   $[\gf{f}^-_{\A_\lambda}]=\gt{b}$}
   \} .
$$
Note that the resolvent set $\R\setminus\spec A\not=\emptyset$
is non-empty by Remark~\ref{rem:EV}.
Observe that, by Proposition~\ref{prop:cascasca-intro-NEW},
for $A\in \Ff_\gf{h}$ to lie in $\Ff_\gf{h}^{\gt{a}\gt{b}}$ it
suffices that for one $\lambda\in\R\setminus\spec A$ the positive and
negative growth types are
$[\gf{f}^+_{\A_\lambda}]=\gt{a}$ and $[\gf{f}^-_{\A_\lambda}]=\gt{b}$.
Therefore, alternatively, we can define $\Ff_\gf{h}^{\gt{a}\gt{b}}$
for $\gt{a},\gt{b}\in\gtset_1$ as well in the form
\begin{equation}\label{eq:Ffab-0}
   \Ff_\gf{h}^{\gt{a}\gt{b}}
   :=\{
   A\in \Ff_\gf{h} \mid
   \exists \lambda_0\in\R\setminus\spec A\colon
   \text{$[\gf{f}^+_{\A_{\lambda_0}}]=\gt{a}$ $\wedge$
   $[\gf{f}^-_{\A_{\lambda_0}}]=\gt{b}$}
   \} .
\end{equation}
Note that $(\Ff_\gf{h}^{\gt{a}\gt{b}})^*\subset \Ff_\gf{h}^{\gt{a}\gt{b}}$
by choosing $\lambda_0=0$ in which case $\A_{\lambda_0}=A$.

\smallskip
Observe that $\Ff_\gf{h}^{\gt{a}\gt{b}}$ is translation invariant in the following sense
\begin{equation}\label{eq:translation-invariance}
   A\in\Ff_\gf{h}^{\gt{a}\gt{b}}
   \qquad\Leftrightarrow\qquad
   \A_\lambda \in\Ff_\gf{h}^{\gt{a}\gt{b}}\quad\forall\lambda\in\R.
\end{equation}
Indeed '$\Leftarrow$' is obvious.
For '$\Rightarrow$' assume that $A\in \Ff_\gf{h}^{\gt{a}\gt{b}}$.
This means that there exists $\lambda_0\in\R\setminus\spec A$ such that
$\A_{\lambda_0}\in(\Ff_\gf{h}^{\gt{a}\gt{b}})^*$.
Given $\lambda\in\R$, note that
$T_{\lambda_0-\lambda} T_\lambda A=T_{\lambda_0} A=\A_{\lambda_0}\in
(\Ff_\gf{h}^{\gt{a}\gt{b}})^*$,
thus $T_\lambda A\in \Ff_\gf{h}^{\gt{a}\gt{b}}$ by~(\ref{eq:Ffab-0})
with $A$ and $\lambda_0$ replaced by $T_\lambda A$ and
$\lambda_0-\lambda$, respectively.

\begin{theoremABC}\label{thm:A}
Suppose that $\gt{a},\gt{b}\in\gtset_1$ are shift-invariant growth types.
Then $\Ff_\gf{h}^{\gt{a}\gt{b}}$ is open and closed in $\Ff_\gf{h}$.
\end{theoremABC}

The proof of Theorem~\ref{thm:A} is based on
Proposition~\ref{prop:cascasca-intro-NEW} and the following result.

\begin{theoremABC}\label{thm:B}
Let $A\in\Ff_\gf{h}^*$.
Then there exists an open neighborhood $\Uu_A$ of $A$ in $\Ff_\gf{h}^*$
such that the spectral projection selector maps
$$
    \Pi_\pm\colon\Uu_A\to \Ll(H_\frac{1}{2})\cap\Ll(H_\frac{3}{2})
   ,\quad
   B\mapsto \Pi_\pm(B)
$$
defined by~(\ref{eq:Pi+}) are continuous.
\end{theoremABC}

In the next theorem we discuss the question when $\Ff_\gf{h}^{\gt{a}\gt{b}}$ is non-empty.

\begin{theoremABC}\label{thm:C}
Let $\gf{h}=\gf{h}(H_0,H_1)$ be a pair growth function.
\begin{itemize}\setlength\itemsep{0ex} 
\item[\rm a)]
  If $\gt{a},\gt{b}\in\gtset_1$ are shift-invariant growth
  types and $\Ff_\gf{h}^{\gt{a}\gt{b}}\not=\emptyset$,
  then the pair growth type $[\gf{h}]=[\gf{h}]_1$ is shift-invariant.
\end{itemize}
If the pair growth type $[\gf{h}]=[\gf{h}]_1$ is shift-invariant, then
the following are true:
\begin{itemize}\setlength\itemsep{0ex} 
\item[\rm b)]
  If a growth type $\gt{a}\in\gtset$ satisfies
  $(\Ff_\gf{h}^{\gt{a}\gt{a}})^*\not=\emptyset$, then $\gt{a}=\gt{a}_1$ is shift-invariant.
\item[\rm c)]
  There are infinitely many shift-invariant pairs $\gt{a},\gt{b}\in\gtset_1$
  with $\Ff_\gf{h}^{\gt{a}\gt{b}}\not=\emptyset$.
\end{itemize}
\end{theoremABC}

\begin{corollary}
If the pair growth type $[\gf{h}]=[\gf{h}]_1$ is shift invariant,
then the set $\Ff_\gf{h}$ of Floer type weak Hessians
has infinitely many connected components.
\end{corollary}

\begin{proof}
By Theorem~\ref{thm:A} and
Theorem~\ref{thm:C}~c):
$(\gt{a},\gt{b})\not=(\gt{a}^\prime,\gt{b}^\prime)$ $\Rightarrow$
$\Ff_\gf{h}^{\gt{a}\gt{b}}\cap \Ff_\gf{h}^{\gt{a}^\prime\gt{b}^\prime}=\emptyset$.
\end{proof}

It is interesting to contrast the corollary
with the result of Atiyah and Singer \cite{Atiyah:1969a}
for the case of real self-adjoint Fredholm operators from separable infinite
dimensional Hilbert space into itself.
There, in the corollary on p.\,307,
they show that the complement of the essentially positive and
essentially negative Fredholm operators
is a classifying space for the real $\mathrm{K}$-theory functor 
$\mathrm{K}\R^{-7}$, see also~\cite{Phillips:1996a}.

\begin{question}
Given shift-invariant growth types $\gt{a},\gt{b}\in\gtset_1$,
is the space $\Ff_\gf{h}^{\gt{a}\gt{b}}$ connected?
\end{question}

This question is related to the question if there is a scale version
of the Atiyah-J\"anich theorem, \cite[Satz\,2]{Janich:1965a}
and \cite[Thm.\,A1 p.\,154]{Atiyah:1967a}.
\\
The Atiyah-J\"anich theorem also plays an important role
in the result of Atiyah-Singer mentioned above.
Both results are based on the result of Kuiper~\cite{Kuiper:1965a} on the
contractibility of the infinite dimensional orthogonal group. 
\\
A scale version of the Kuiper theorem can be found in the book by
Kronheimer and Mrowka~\cite[Prop.\,33.1.5]{Kronheimer:2007a}.

\begin{question}
Is the assumption that our Fredholm operators restrict to operators
from $H_2$ to $H_1$ necessary?
\end{question}

\smallskip\noindent
{\bf Acknowledgements.}
We would like to thank Gilles Pisier for explaining
to us an obstruction for a matrix to be a Schur multiplier
(Corollary~\ref{cor:obstruction}).
UF~acknowledges support by DFG grant
FR~2637/4-1.

\newpage
\subsection{Motivation and general perspective}

Nowadays we have many powerful examples of Floer theories,
although we do not have a good notion what a Floer theory actually is.
This issue especially gets important in the extension to Hamiltonian
delay equations.

\smallskip
Hamiltonian delay equations have a long history.
One of their first occurrences was a paper by Carl Neumann
in 1868~\cite{Neumann:1868a}, \href{http://www.ifi.unicamp.br/~assis/Weber-in-English-Vol-4.pdf}{English translation}~\cite{assis:2021a},
in which Carl Neumann deduces Wilhelm Weber's force law in
electrodynamics from a retarded Coulomb potential;
see~\cite{Frauenfelder:2019a} for a modern treatment.
Quite recently Hamiltonian delay equations found surprising
new applications due to an interesting paper by Barutello, Ortega, and
Verzini~\cite{Barutello:2021b}.
In this paper the authors discovered a new non-local regularization
for two-body collisions.

\smallskip
For autonomous Hamiltonian systems there are various techniques to
regularize two-body collisions by blowing up the energy hypersurface;
see~\cite{Frauenfelder:2018b} for a discussion of such techniques from
the point of view of symplectic and contact topology.
For non-autonomous Hamiltonian systems there is no preserved energy
and therefore Barutello, Ortega, and Verzini don't blow up an energy
surface, but the loop space by reparametrizing every loop
individually.
The new equation they obtain is a Hamiltonian delay equation
which for non-collisional orbits can be transformed back to
the original Hamiltonian ODE. This has the advantage that it also
covers collision orbits.
The resulting action functionals have fascinating
properties~\cite{Frauenfelder:2021e,Cieliebak:2023a}.
How Floer theory can be used to prove existence for some Hamiltonian
delay equations was already discussed in~\cite{Albers:2020a}.
With the new techniques discovered by Barutello, Ortega, and Verzini
the question of Floer homologies for Hamiltonian delay equations
becomes now a topic of major concern.

\smallskip
In view of the growing interest in Hamiltonian delay equations
the question becomes pressing what a Floer theory actually is.
When Floer introduced his celebrated homology to prove the \Arnold\,
conjecture~\cite{floer:1988c,floer:1989c}
he considered the gradient flow equation with respect to a very weak metric.
This was in sharp contrast to previous work by Conley and
Zehnder~\cite{Conley:1983a,conley:1984a} where they considered an
interpolation metric that led to a gradient flow equation which can
be interpreted as an ODE on a Hilbert manifold.
Different from that the Hessian in Floer's theory is not a
self-adjoint operator from a fixed Hilbert space into itself,
but a symmetric operator on a scale of Hilbert spaces.

\smallskip
For Hamiltonian delay equations the structure of the space of such
Hessians is now a topic of major interest.
As we show in the present article this space differs drastically from
the case of self-adjoint operators on a fixed Hilbert space.
While in our case the space has infinitely many different connected
components, it is connected in the later case~\cite{Atiyah:1967a,Phillips:1996a}.

\smallskip
The proof of our main result uses interpolation theory and Schur
multipliers. That Floer theory is related to interpolation theory was
used before in several
contexts~\cite{robbin:1995a,albers:2013a,Simcevic:2014a,Hofer:2021a}.
The connection between Floer theory and Schur multipliers seems to be
completely new.

\subsection{Outline}

This article is organized as follows.

\smallskip
In Section~\ref{sec:frac} ``Growth functions and growth types''
we recall in \S\,\ref{sec:growth-functions}
the pointwise $\gf{f}\cdot\gf{g}$ and the Kang, or union, product $\gf{f}*\gf{g}$
on the set $\gfset$ of growth functions.
The latter is the growth function $\N\to(0,\infty)$
whose values are the union of values $\gf{f}(\N)\cup\gf{g}(\N)$
enumerated increasingly with multiplicities.
The two products are compatible in the sense that $(f*g)^2=f^2*g^2$.
\\
In \S\,\ref{sec:growth-types}
we endow the set of growth types $\gtset:=[\gfset]:=\gfset/\sim$
with a partial order ``$\le$''.
Both, pointwise and Kang product, descend to $[\gfset]$
as commutative associative products, same notation.
In \S\,\ref{sec:shift} we look at shift invariance of growth types and
their Kang and pointwise products.
In \S\,\ref{sec:scale} we introduce scale invariance of growth types
an example of which are known Floer growth types.
\S\,\ref{sec:Kang-product-representations}
establishes the existence of infinitely many
pairs of shift invariant growth types $[\gf{f}_i]$ representing a given
shift invariant growth type $[\gf{f}]=[\gf{f}_i*\gf{f}_j]$.

\smallskip
Section~\ref{sec:weak-Hess}
``Weak Hessians'' provides properties
of general weak Hessians $A\in\Hh_\gf{h}$
such as ``regularity and spectrum'' in
\S\,\ref{sec:spectrum} which also deals with properties
inherited from $A\colon H_1\to H_0$
by the restriction $A\colon H_2\to H_1$
such as invertibility, eigenspaces, and discrete, finite multiplicity,
pure eigenvalue spectrum.
A crucial side fact is that $\R\setminus\spec\,A$ is non-empty.
\\
\S\,\ref{sec:Hessian-growth-types}
``Proof of Proposition~\ref{prop:cascasca-intro-NEW}''
shows that the set $(\Ff_\gf{h}^{\gt{a}\gt{b}})^*$ of
invertible Floer Hessians $A$ with prescribed negative and positive
growth types $\gt{a}$ and $\gt{b}$ is preserved under translation
$T_\lambda A:=A-\lambda\iota$ by the elements $\lambda$ of $\R\setminus\spec\,A$.

\smallskip
In Section~\ref{sec:thmA}
``Proof of Theorem~\ref{thm:A}''
we show using Theorem~\ref{thm:B} (continuity of $B\mapsto \Pi_+^B$)
that the restriction of the projection $\Pi^B_+$
to the image of  the projection $\Pi^A_+$ gives rise to a scale
isomorphism 
$$
   \Pi_+^B|_A\colon
   (\Pi_+^A  H_{\frac{1}{2}}, \Pi_+^A  H_{\frac{3}{2}})
   \to (\Pi_+^B H_{\frac{1}{2}}, \Pi_+^B H_{\frac{3}{2}})
$$
whenever $B$ is close enough to $A$.
Then  Theorem~\ref{thm:A} follows from the result that a scale
isomorphism preserves the growth type of the pair, as proved in
the appendix in Theorem~\ref{thm:f-g-equiv}.

\smallskip
Section~\ref{sec:Schur}
``Schur multipliers''
deals with the Schur product of two matrices which is defined by
taking the entry-wise product.
A Schur multiplier is a matrix which, via Schur multiplication,
keeps the space of matrices corresponding to bounded linear maps
$\ell^2\to\ell^2$ invariant.
In this section we explain a criterium of Schur to show that a given
matrix is a Schur multiplier and mention a result by Grothendieck
which shows that the criterium of Schur is not just sufficient but
also necessary. Via the result of Grothendieck we discuss a criterium
explained to us by Gilles Pisier to deduce that a given matrix is not
a Schur multiplier.

The positive results of Schur are one of the major ingredients
to prove Theorem~\ref{thm:B}.
While the obstruction is the reason that our proof of
Theorem~\ref{thm:B} does not extend from $H_\frac12$ to $H_0$.
This obstruction also the main reason why for Theorem~\ref{thm:A}
we need that our weak Hessians restrict to Fredholm operators of index
zero from $H_2$ to $H_1$.

\smallskip
Section~\ref{sec:thmB}
``Proof of Theorem~\ref{thm:B}''
is the analytical heart of the present article.
To prove Theorem~\ref{thm:B}, inspired by Kato's perturbation theory
for linear operators, we write the projection $\Pi_+^A$ for $A$ invertible
as an operator-valued integral in the complex plane
where the contour $\gamma_+$ encircles all positive eigenvalues of $A^{-1}$,
but no negative eigenvalue of $A^{-1}$, as illustrated by
Figure~\ref{fig:fig-concatenation}.

The spectrum of $A^{-1}$ has the origin as its unique accumulation
point. Therefore we decompose our contour $\gamma_+=\beta\#\alpha_+$
into two paths, where the dangerous part $\beta$ of the path is going
through the origin and the harmless part $\alpha_+$ stays at
bounded distance of the spectrum of $A^{-1}$.
Using restriction and inclusion,
we interpret the projection as a linear operator from $H_1$ to $H_0$,
see Figure~\ref{fig:Pi_+},
and explicitly compute its differential with
respect to variation $\Delta$ of $A$.

The main ingredient to prove Theorem~\ref{thm:B}
is to show that the images of this operator are actually bounded
linear maps from $H_\frac12$ to $H_\frac12$.
Along the harmless part $\alpha_+$ this follows by interpolation
theory. While along the dangerous part $\beta$ we need the theory of
Schur multipliers.

\smallskip
In Appendix~\ref{sec:HS-pairs}
``Hilbert space pairs'',
following~\cite{Frauenfelder:2009b},
we define the growth type $[\gf{h}]$ of a pair of Hilbert spaces $(H_0,H_1)$
and show that up to isomorphism a Hilbert space pair
is uniquely characterized by its growth type

\smallskip
In Appendix~\ref{sec:interpolation}
``Interpolation''
we explain a special case of an interpolation
theorem by Stein going back to a method discovered by Thorin.

\section{Growth functions and growth types}
\label{sec:frac}

A \textbf{growth function} is a monotone unbounded 
function
$$
   \gf{f} \colon \mathbb{N} \to (0,\infty) .
$$
We abbreviate by $\ell^2_\gf{f}=\ell^2_\gf{f}(\N)$ 
the space of all sequences $x=(x_\nu)_{\nu \in \mathbb{N}}$ satisfying
$$\sum_{\nu=1}^\infty \gf{f}(\nu)x_\nu^2 <\infty.$$
The space $\ell^2_\gf{f}$ becomes a Hilbert space if we endow it with the inner product
\begin{equation}\label{eq:ell^2_f}
     \langle x,y \rangle_\gf{f}=\sum_{\nu=1}^\infty \gf{f}(\nu)x_\nu y_\nu.
\end{equation}
Let $\Norm{x}_\gf{f}:=\sqrt{\langle x,x \rangle_\gf{f}}$ be the associated norm.
For $r \in \mathbb{R}$ 
we abbreviate
\begin{equation}\label{eq:ell^2_f^r}
     H_r:=\ell^2_{\gf{f}^r}
     ,\qquad
     \INNER{\cdot}{\cdot}_r:=\INNER{\cdot}{\cdot}_{\gf{f}^r}
     ,\qquad
     \Norm{\cdot}_r:=\Norm{\cdot}_{\gf{f}^r} .
\end{equation}
Because the function $\gf{f}$ is monotone increasing,
it follows that whenever $s\le r$ there is an inclusion 
\[
     H_r \hookrightarrow H_s
\]
of Hilbert spaces and the corresponding linear inclusion operator is bounded,
a bound being~$1$ since $\gf{f}\ge 1$.
For strict inequality $s<r$, by unboundedness of $\gf{f}$, the inclusion
operator is compact. Moreover, its image is
dense. For integers $r,s\in\N_0$ these facts are proved
in~\cite[Thm.\,8.1]{Frauenfelder:2021b}. In the real case 
the same arguments prove density and compactness
of the inclusion $\ell^2_{\gf{f}^r}\INTO\ell^2$ for $r>0$ and
$\ell^2\INTO\ell^2_{\gf{f}^r}$ for $r<0$.
To deal with the case $s<r$ use the isometric isomorphisms $\phi_r$ below.
In particular,
$$H=H_0 \supset H_1 \supset H_2 \supset \ldots$$
defines an $\SC$-structure, also called a Banach scale, in the sense of Hofer,
Wysocki and Zehnder, see~\cite{Hofer:2007a,Hofer:2021a}.
As all levels are Hilbert spaces, one speaks of a Hilbert scale.
For $r \in \mathbb{R}$ we abbreviate by $H^r$ the Hilbert scale with levels
$$
     H^r_m:=H_{r+m}, \quad m\in \mathbb{N}_0.
$$
For $r \in \mathbb{N}$ these are the shifted scale spaces
introduced in \cite[\S 2.1]{Hofer:2007a}.
In the set-up of fractal Hilbert scales, introduced in 
\cite{Frauenfelder:2021b},
the shift works more generally for every $r \in \mathbb{R}$.
Observe that the $\SC$-Hilbert spaces $H^r$ and $H$ are isomorphic via
the levelwise isometries
$$\phi_r \colon H \to H^r, \quad (x_\nu)_\nu \mapsto
\big(\gf{f}(\nu)^{-r/2}x_\nu\big)_\nu.$$

\subsection{Growth functions}
\label{sec:growth-functions}

Let $\gfset$ be the \textbf{set of growth functions}.
The \textbf{shift} $\gf{f}_1$ of $\gf{f}\in\gfset$ is defined by
$\gf{f}_1(\nu):=\gf{f}(\nu+1)$.
On $\gfset$ there are two products.
The \textbf{pointwise product}
$$
   \cdot\colon\gfset\times\gfset\to\gfset
   ,\quad
   (\gf{f},\gf{g})
   \mapsto \gf{f}\cdot\gf{g}
   ,\qquad
   (\gf{f}\cdot\gf{g})(\nu):=\gf{f}(\nu)\cdot\gf{g}(\nu)
   ,\qquad
   \gf{f}^2:=\gf{f}\cdot \gf{f},
$$
and the \textbf{Kang product}~\cite[(3.1)]{kang:2019a}
$$
   *\colon\gfset\times\gfset\to\gfset
   ,\quad
   (\gf{f},\gf{g})
   \mapsto \gf{f}*\gf{g} .
$$
The product $\gf{f}*\gf{g}$ is the monotone
unbounded function $\gf{f}*\gf{g}\colon\N\to(0,\infty)$
whose value at $\nu$ is the $\nu^{\rm th}$ smallest number, counted
with multiplicities, among the members of the two lists
$(\gf{f})=\left(\gf{f}(1), \gf{f}(2),\dots \right)$ and
$(\gf{g})=\left(\gf{g}(1), \gf{g}(2),\dots\right)$.\footnote{
  If $\gf{f}$ and $\gf{g}$ are strictly monotone and their images
  $\gf{f}(\N)\cap\gf{g}(\N)=\emptyset$ are disjoint, then
  the Kang product is given by the formula
  $(\gf{f}*\gf{g})(\nu)
  =\min\left\{(\gf{f}(\N)\cup\gf{g}(\N))\setminus
  \bigcup_{\mu=1}^{\mu-1}(\gf{f}*\gf{g})(\mu)\right\}$.
  }
\\
For example the lists
$$
   (\gf{f})=\left(1,1,2,3,5,\dots \right)
   ,\quad
   (\gf{g})=\left(1,2,2,3,4,5,\dots \right),
$$
have Kang product
$$
   (\gf{f}*\gf{g})=\left(1,1,1,2,2,2,3,3,4,5,5,\dots
   \right) .
$$
Recall that for a real $x$ \textbf{ceiling} $\lceil x \rceil$ and
\textbf{floor} $\lfloor x \rfloor$ denote the integer next to $x$,
larger or equal $x$, respectively less or equal $x$, in symbols
$$
   \lceil x \rceil
   :=\min\{n\in\Z\mid n\ge x\}
   ,\;\;
   \lfloor x \rfloor
   :=\max\{n\in\Z\mid n\le x\}
   ,\;\;
   \lfloor \tfrac{x}{2} \rfloor+1=\lceil \tfrac{x}{2}+\tfrac12 \rceil .
$$
The Kang square and its shift satisfy
\begin{equation*}
   (\gf{f}*\gf{f})(\nu)
   =f(\lceil \tfrac{\nu}{2} \rceil)
   ,\qquad
   (\gf{f}*\gf{f})_1(\nu)
   =f(\lfloor \tfrac{\nu}{2}\rfloor+1)
   =f_1(\lfloor \tfrac{\nu}{2}\rfloor) .
\end{equation*}
In particular, it holds
\begin{equation}\label{eq:Kang-square-2}
   (\gf{f}*\gf{f})(2\nu)
   =f(\nu)
   ,\qquad
   (\gf{f}*\gf{f})_1(2\nu)
   =f(\nu+1)=f_1(\nu) .
\end{equation}
Kang product and pointwise square are compatible in the sense that
\begin{equation}\label{eq:Kang-ptw_square}
   (\gf{f}*\gf{g})^2=\gf{f}^2*\gf{g}^2 .
\end{equation}
This holds since squaring preserves the ordering.

\subsection{Growth types}
\label{sec:growth-types}

On $\gfset$ there is the equivalence relation of having the \textbf{same growth type}
\begin{equation}\label{eq:same-growth-type}
   \gf{f}\sim\gf{f}^\prime
   \quad:\Leftrightarrow\quad
   \exists c>0\colon
   \frac{1}{c} \gf{f}(\nu)\le\gf{f}^\prime(\nu)\le c \gf{f}(\nu)
   ,\; \forall \nu\in\N.
\end{equation}
An equivalence class $\gt{f}:=[\gf{f}]$ is referred to
as a \textbf{growth type}. The set of equivalence classes
$\gtset:=[\gfset]:=\gfset/\sim$ is the \textbf{set of growth types}.

Simple examples of growth functions of the same growth type
are growth functions which just differ by an additive constant
as we show next.

\begin{lemma}\label{le:growth-constant}
Let $\gf{f}$ and $\gf{g}$ be growth functions.
Assume that there exists $\lambda\in\R$ such that
$\gf{f}(\nu)=\gf{g}(\nu)+\lambda$ for every $\nu\in\N$.
Then $\gf{f}$ and $\gf{g}$ have the same growth type,
in symbols $\gf{f}\sim\gf{g}$.
\end{lemma}

\begin{proof}
We can assume without loss of generality that $\lambda$ is positive.
Indeed if $\lambda=0$ the assertion is trivial and if $\lambda<0$ we
can switch the role of $\gf{f}$ and $\gf{g}$.
\\
Since $\gf{g}>0$
we have $0<\gf{g}(1)=\gf{f}(1)-\lambda$
and therefore $\lambda\in(0,\gf{f}(1))$.
\\
The function $x\mapsto \frac{x}{x-\lambda}$ for $x>\lambda$
is monotone decreasing. Indeed its derivative is given by
$-\lambda/(x-\lambda)^2<0$. As $\gf{f}$ is monotone increasing
we have $\frac{\gf{f}(\nu)}{\gf{f}(\nu)-\lambda}
\le\frac{\gf{f}(1)}{\gf{f}(1)-\lambda}$.
Hence
$\gf{f}(\nu)
=\frac{\gf{f}(\nu)}{\gf{f}(\nu)-\lambda}
(\gf{f}(\nu)-\lambda)
\le
\frac{\gf{f}(1)}{\gf{f}(1)-\lambda}(\gf{f}(\nu)-\lambda)
=\frac{\gf{f}(1)}{\gf{f}(1)-\lambda} \gf{g}(\nu)
$.  
Therefore the growth function $\gf{f}$ is bounded from above
by the growth function $\gf{g}$ up to a multiplicative constant.
\\
Since $\lambda>0$ we have $\gf{g}(\nu)\le\gf{f}(\nu)$
and therefore the two growth functions have the same growth type.
\end{proof}

\begin{definition}[Partial order]
On the set $[\gfset]$ of growth types we define a partial order as follows.
Assume $[\gf{f}]$ and $[\gf{g}]$ are growth types.
We write $[\gf{f}]\le[\gf{g}]$ if there exists a constant $d>0$ such that
$\gf{f}(\nu)\le d \gf{g}(\nu)$ for every $\nu\in\N$.
We write $[\gf{f}]<[\gf{g}]$ if $[\gf{f}]\le[\gf{g}]$ and $[\gf{f}]\not=[\gf{g}]$.
\end{definition}

Pointwise and Kang product both descend to $[\gfset ]$
as commutative associative products, still denoted by
$$
   \cdot,*\colon [\gfset]\times[\gfset]\to[\gfset]
   ,\quad
   [\gf{f}]\cdot[\gf{g}]:=[\gf{f}\cdot\gf{g}]
   ,\quad
   [\gf{f}]*[\gf{g}]:=[\gf{f}*\gf{g}] .
$$
For the pointwise product this is easy,\footnote{
   In~(\ref{eq:same-growth-type}),
   if $\gf{f}\sim\gf{f}^\prime$ with constant $c$
   and $\gf{g}\sim\gf{g}^\prime$ with constant $d$,
   then $\gf{f}\cdot \gf{g}\sim
   \gf{f}^\prime\cdot\gf{g}^\prime$ with $cd$.
  }
for the Kang product see~\cite[\S\,3]{kang:2019a}.

\begin{lemma}[Kang~{\cite{kang:2019a}}]\label{le:Kang}
The Kang product is compatible with the partial order on the set of
growth types $\gtset$ as follows.
\begin{itemize}\setlength\itemsep{0ex} 
\item[\rm (i)]
  $[\gf{g}_1]\le[\gf{g}_2]$ $\Rightarrow$ $[\gf{f}*\gf{g}_1]\le[\gf{f}*\gf{g}_2]$;
\item[\rm (ii)]
  $\gf{f}*\gf{g}\le \gf{f}$ pointwise at any $\nu\in\N$ and for all
  growth functions $\gf{f}$ and $\gf{g}$.
\end{itemize}
\end{lemma}

\subsubsection{Shift invariance}
\label{sec:shift}

Note that if two growth functions $\gf{f}$ and $\gf{g}$
have the same growth type, then their shifts have the same growth type
as well, in symbols
$$
   \gf{f}\sim \gf{g}
   \quad\Rightarrow\quad
   \gf{f}_1\sim \gf{g}_1 .
$$
Thus the shift descends to growth type, so the following is well defined
$$
   [\gf{f}]_1:=[\gf{f}_1].
$$

\begin{definition}
A growth type $[\gf{f}]$ is called \textbf{shift invariant}
if $[\gf{f}]_1=[\gf{f}]$.
\end{definition}

\begin{lemma}\label{le:si}
For growth functions $\gf{f},\gf{g}\colon\N\to(0,\infty)$ the following
is true.
\begin{itemize}\setlength\itemsep{0ex} 
\item[\rm (i)]
  $(\gf{f}\cdot\gf{g})_1=\gf{f}_1\cdot\gf{g}_1$.
\item[\rm (ii)]
  $(\gf{f}*\gf{g})_1=\begin{cases}
  \gf{f}_1*\gf{g}&\text{, if $\gf{f} (1)\le \gf{g}(1)$,}\\
  \gf{f}*\gf{g}_1&\text{, if $\gf{g}(1)\le \gf{f} (1)$.}\\
  \end{cases}$
\item[\rm (iii)]
  If the growth types $[\gf{f}]$ and
  $[\gf{g}]$ are shift invariant, then, by~(ii), their
  pointwise product and their Kang product are shift invariant, too,
  in symbols
  $$
     [\gf{f}]=[\gf{f}]_1
     \;\;\wedge\;\;
     [\gf{g}]=[\gf{g}]_1
     \qquad\Rightarrow\qquad
     [\gf{f}\cdot\gf{g}]_1=[\gf{f}\cdot\gf{g}]
     \;\;\wedge\;\;
     [\gf{f}*\gf{g}]_1=[\gf{f}*\gf{g}] .
  $$
\item[\rm (iv)]
  A Kang square $[\gf{f}*\gf{f}]$ is shift invariant iff
  the growth type $[\gf{f}]$ itself is, in symbols
  $$
     [\gf{f}*\gf{f}]_1=[\gf{f}*\gf{f}]
     \quad\Leftrightarrow\quad
     [\gf{f}]_1=[\gf{f}] .
  $$
\item[\rm (v)]
  A pointwise square $[\gf{f}\cdot\gf{f}]$ is shift invariant iff
  the growth type $[\gf{f}]$ itself is, in symbols
  $$
     [\gf{f}\cdot\gf{f}]_1=[\gf{f}\cdot\gf{f}]
     \quad\Leftrightarrow\quad
     [\gf{f}]_1=[\gf{f}] .
  $$
\end{itemize}
\end{lemma}

\begin{proof}
(i) and (ii) follow directly by definition. 
(iii) follows from (i) and (ii).
We prove (iv). '$\Rightarrow$' 
Suppose $\gf{f}*\gf{f}\sim (\gf{f}*\gf{f})_1$,
that is there is a constant $c>0$ such that
$\frac{1}{c}(\gf{f}*\gf{f})(\nu)
\le (\gf{f}*\gf{f})_1(\nu)
\le c (\gf{f}*\gf{f})(\nu)$
for every $\nu\in\N$.
In particular, 
$\frac{1}{c}(\gf{f}*\gf{f})(2\nu)
\le (\gf{f}*\gf{f})_1(2\nu)
\le c (\gf{f}*\gf{f})(2\nu)$
for every $\nu\in\N$.
Hence, by~(\ref{eq:Kang-square-2}),
$\frac{1}{c}\gf{f}(\nu)
\le \gf{f}_1(\nu)
\le c \gf{f}(\nu)$
for every $\nu\in\N$.
Thus $\gf{f}\sim \gf{f}_1$.
 '$\Leftarrow$' 
Choose $\gf{g}=\gf{f}$ in~(iii).
\\
(v)~This follows by taking roots.
\end{proof}

\subsubsection{Scale invariance}
\label{sec:scale}

\begin{definition}
A growth type $[\gf{f}]$ is called \textbf{scale invariant}
if there exists a constant $c$ such that
$\gf{f}(2\nu)\le c\gf{f}(\nu)$ for every $\nu\in\N$.\footnote{
  Well defined: In~(\ref{eq:same-growth-type}),
  if $\gf{f}\sim\gf{f}^\prime$ with constant $d$,
  then $\gf{f}^\prime(2\nu)\le cd^2\gf{f}^\prime(\nu)$.
  }
\end{definition}

\begin{lemma}\label{le:fg456}
Suppose that $[\gf{f}]$ is a scale invariant growth type.
Then $[\gf{f}]$ is (i) shift-invariant and (ii) Kang idempotent.
\end{lemma}

\begin{proof}
(i)~It holds $\gf{f}_1(\nu):= \gf{f} (\nu+1)\le \gf{f} (2\nu)\le c
\gf{f} (\nu)$ where inequality one uses monotonicity of $\gf{f}$.
(ii)~It holds~$
   (\gf{f} *\gf{f})(2\nu)
   \le \gf{f} (2\nu)
   \le c \gf{f} (\nu)
   =c (\gf{f}*\gf{f})(2\nu)
$.
\end{proof}

\begin{example}[Floer growth types]
Given a real parameter $r$, the growth function
$\gf{f}_r(\nu)=\nu^r$ for $\nu\in\N$ 
is scale invariant with $c=2^r$.
The growth type of $\gf{f}_r$ satisfies
$[f_r]=[\gf{f}_r*\gf{f}_s]$ whenever $s\ge r$.
The growth types of Floer theories studied so far
are of this type, see~\cite[p.\,378]{Frauenfelder:2021b}.
\end{example}

\begin{example}
The growth function $\gf{f}(\nu)=e^\nu$ is not scale invariant.
\end{example}

\begin{lemma}\label{le:gujg45}
Assume $[\gf{f}]$ is a scale invariant growth type.
Then there exists a strictly larger scale invariant growth type $[\gf{g}]$,
in symbols $[\gf{g}]>[\gf{f}]$.
\end{lemma}

\begin{proof}
Pick $\kappa>0$. Define $\gf{g}(\nu):=\nu^\kappa \gf{f}(\nu)$.
Then $[\gf{f}]\le [\gf{g}]$ with constant $d=1$.
On the other hand $[\gf{f}]\not= [\gf{g}]$, since
$\nu\mapsto \nu^\kappa$ is not bounded by any constant $c$.
For scale invariance note that
$
   \gf{g}(2\nu)=(2\nu)^\kappa \gf{f}(2\nu)
   \le 2^\kappa \nu^\kappa c \gf{f}(\nu)
   =2^\kappa c \gf{g}(\nu)
$.
\end{proof}

\subsubsection{Kang product representations}
\label{sec:Kang-product-representations}

\begin{theorem}\label{thm:infinitely-many}
Assume $[\gf{f}]\in\gtset_1$ is a shift invariant growth type.
Then there exist infinitely many pairs $([\gf{f}_1], [\gf{f}_2])\in
\gtset_1\times \gtset_1$ such that $[\gf{f}]=[\gf{f}_1*\gf{f}_2]$.
\end{theorem}

\begin{proof}
\smallskip
\noindent
\textbf{Case~1 ($\mbf{[\gf{f}]}$ scale invariant).} 
By Lemma~\ref{le:gujg45} there exists an infinite
sequence of scale invariant growth types
$[\gf{f}]=[\gf{g}_1]<[\gf{g}_2]<[\gf{g}_3]<\dots$.
\\
We claim that $[\gf{f}]=[\gf{f}*\gf{g}_k]$ whenever $k\in\N$:
Since $[\gf{f}]$ is scale invariant it is Kang idempotent by Lemma~\ref{le:fg456}.
So the assertion is true for $k=1$.
Hence, using parts (i) and (ii) of Lemma~\ref{le:Kang}, we have
$[\gf{f}]=[\gf{f}*\gf{f}]\le [\gf{f}*\gf{g}_k]\le [\gf{f}]$.
This proves that $[\gf{f}]=[\gf{f}*\gf{g}_k]$ for every $k\in\N$.
Since scale invariant growth types are shift invariant, by
Lemma~\ref{le:fg456}~(i), Case 1 is proved.

\smallskip
\noindent
\textbf{Case~2 ($\mbf{[\gf{f}]}$ not scale invariant).} 
Pick an integer $k\ge 2$.
We decompose the image of $\gf{f}$ into two parts as follows.
For $\nu\in\N$ we define $f_k(\nu):=f(k\nu)$
and
$$
   \gf{g}_k(\nu)
   :=\begin{cases}
     \gf{f}(\nu)&\text{, $\nu\in\{1,\dots,k-1\}$}\\
     \gf{f}(\nu+1)&\text{, $\nu\in\{k,\dots,2k-2\}$}\\
     \gf{f}(\nu+2)&\text{, $\nu\in\{2k-1,\dots,3k-3\}$}\\
     \gf{f}(\nu+3)&\text{, $\nu\in\{3k-2,\dots,4k-4\}$}\\
     \quad\vdots&\quad\vdots
   \end{cases}
$$
Observe that
$
   \gf{g}_k(\nu)
   =\gf{f}(\nu+\ell)
$
whenever $\ell k-(\ell-1)\le\nu\le(\ell+1)k-(\ell+1)$ for some $\ell\in\N_0$.
Then $\gf{f}_k*\gf{g}_k=\gf{f}$.
Moreover, shift invariance of $\gf{f}$ implies that
$\gf{f}_k$ and $\gf{g}_k$ are shift invariant.\footnote{
  e.g. $\gf{f}_k(\nu+1)=\gf{f}(k\nu+k)\sim \gf{f}(k\nu)=\gf{f}_k(\nu)$
  }

\smallskip\noindent
\textsc{Claim.}
Let $k_2>k_1\ge 2$. Then the growth types of $\gf{f}_{k_1}$
and $\gf{f}_{k_2}$ are different.

\smallskip\noindent
\textsc{Proof of Claim.}
We argue by contradiction and assume that there exists a constant $c>0$
such that
\begin{equation}\label{eq:hg44re}
   {\color{gray}\gf{f}(k_2\nu)=\,}
   \gf{f}_{k_2}(\nu)\le c\gf{f}_{k_1}(\nu)
   {\color{gray}\, =c\gf{f}(k_1\nu)}
\end{equation}
for every $\nu\in\N$.
To derive a contradiction we show that this implies that there exists a
constant $C>0$ such that
\begin{equation}\label{eq:hg4554re}
   \gf{f}(2\nu)\le C\gf{f}(\nu)
\end{equation}
for every $\nu\in\N$ which contradicts the assumption that $\gf{f}$ is not scale
invariant.

\medskip\noindent
\textsc{Subclaim.}
Suppose $\nu_1\in \N$ satisfies $\nu_1\ge k_1\sqrt{k_2}/(\sqrt{k_2}-\sqrt{k_1})$.
Then there exists $\nu_2\in \N$ such that a)
$\gf{f}(\nu_1)\ge \frac{1}{c} \gf{f}(\nu_2)$
and b) $\nu_2\ge \sqrt{k_2/k_1} \nu_1$.

\smallskip\noindent
\textsc{Proof of Subclaim.}
There exist unique $m\in\N_0$ and
$r\in\{0,\dots,k_1-1\}$ such that $\nu_1=k_1m+r$.
We set $\nu_2:=mk_2$. 
\\
a)~We estimate
$
   \gf{f}(\nu_1)
   \ge\gf{f}(mk_1)
   \ge\tfrac{1}{c}\,\gf{f}(mk_2)
$
where the first step is by monotonicity of $\gf{f}$ and the second
step is by~(\ref{eq:hg44re}) for $\nu=m$.
b)~We estimate
$
   \frac{k_1\sqrt{k_2}}{\sqrt{k_2}-\sqrt{k_1}}\le\nu_1=k_1m+r<(m+1)k_1
$.
Resolving for $m$ we get
\begin{equation}\label{eq:dvbhr44}
   m\ge \frac{\sqrt{k_2}}{\sqrt{k_2}-\sqrt{k_1}}-1
   =\frac{\sqrt{k_1}}{\sqrt{k_2}-\sqrt{k_1}} =:m_0.
\end{equation}
The function $\rho\colon (-1,\infty)\to\R$, $x\mapsto x/(1+x)$ is strictly
monotone increasing, indeed its derivative is $\rho^\prime(x)=1/(1+x)^2>0$.
Therefore
$$
   \frac{m}{m+1}=\rho(m)
   \ge \rho(m_0)
   =\sqrt{k_1/k_2} .
$$
We use this estimate in the last step of what follows to get
$$
   \nu_2:=mk_2=\frac{mk_2}{\nu_1}\nu_1
   \ge \frac{m}{m+1} \frac{k_2}{k_1}\nu_1
   \ge\sqrt{\frac{k_2}{k_1}}\nu_1 
$$
where we used $\nu_1<(m+1)k_1$ in the first inequality.
This proves the Subclaim.

\medskip
Since $\sqrt{k_2/k_1}>1$ there exists $N\in\N$ such that
$\sqrt{k_2/k_1}^N\ge 2$.
Now suppose that $\nu\ge k_1\sqrt{k_2}/(\sqrt{k_2}-\sqrt{k_1})$.
Then applying the Subclaim iteratively $N$ times, we obtain the scale
invariance criterium 
$\gf{f}(\nu)\ge \frac{1}{c^N}\gf{f}(2\nu)$, but only for large $\nu$.
But along finitely many $\nu$'s the function $\gf{f}$ is bounded.
So~(\ref{eq:hg4554re}) follows with
$$
   C:=\max\left\{c^N,
   \frac{\gf{f}(2 \lceil k_1\sqrt{k_2}/(\sqrt{k_2}-\sqrt{k_1}) \rceil)}{\gf{f}(1)}
   \right\}.
$$
This proves the Claim and Theorem~\ref{thm:infinitely-many}.
\end{proof}

\subsection{Proof of Theorem~\ref{thm:C}}
\label{sec:thm-C}

a) Since $\Ff_\gf{h}^{\gt{a}\gt{b}}$ is non-empty, there exists
$A\in\Ff_\gf{h}^{\gt{a}\gt{b}}$.
Without loss of generality we may assume that $A$ is 
invertible.\footnote{
  Otherwise, consider $\A_\lambda:=A-\lambda\iota$ for
  $\lambda\in \R\setminus\spec A$.
  }
By definition of $\gf{f}_A$ (order the eigenvalue squares) there
is the first identity in
\begin{equation}\label{eq:square-44}
   \gf{f}_A
   =(\gf{f}_A^+)^2*(\gf{f}_A^-)^2
   \stackrel{\text{(\ref{eq:Kang-ptw_square})}}{=}
      (\gf{f}_A^+*\gf{f}_A^-)^2 .
\end{equation}
Hence, by Lemma~\ref{le:si} (iii), the growth type of $\gf{f}_A$ is shift-invariant
and therefore the same is true for the pair growth function $\gf{h}$
by Lemma~\ref{le:f_a-f}.
This proves a).

\medskip\noindent
b)
By Lemma~\ref{le:f_a-f} it holds $\gf{f}_A\sim \gf{h}$.
Therefore, by the assumption $[\gf{h}]=[\gf{h}]_1$, the growth type of
$\gf{f}_A$ is shift invariant, in symbols $[\gf{f}_A]_1=[\gf{f}_A]$.
\\
By assumption $\gf{f}_A^-\sim \gf{f}_A^+$.
Thus, by~(\ref{eq:square-44}) and Lemma~\ref{le:si}~(v),
it follows that the growth type of $\gf{f}_A^+*\gf{f}_A^+$
is shift invariant, in symbols
$[\gf{f}_A^+*\gf{f}_A^+]_1=[\gf{f}_A^+*\gf{f}_A^+]$.
Hence by  Lemma~\ref{le:si}~(iv) it follows that the growth type of
$\gf{f}_A^+$ itself is shift invariant, in symbols
$[\gf{f}_A^+]_1=[\gf{f}_A^+]$.

\medskip\noindent
c)~Theorem~\ref{thm:infinitely-many}.

\section{Weak Hessians}
\label{sec:weak-Hess}

\subsection{Regularity and spectrum}
\label{sec:spectrum}

\begin{lemma}[Regularity]\label{le:regularity}
Let $A$ be a weak Hessian, in symbols $A\in\Hh_\gf{h}$.
If $\xi$ and $A\xi$ both lie in $H_1$,
then $\xi$ actually lies in $H_2$.
\end{lemma}

\begin{proof}
Since $H_2$ is separable the operator $A$, as a map $H_2\to H_1$,
can only have countably many eigenvalues.
Hence there exists $\lambda\in\R$ which is not an eigenvalue.
We claim that $A-\lambda\iota\colon H_2\to H_1$ is an isomorphism.

Since $\lambda$ is not an eigenvalue of $A$ there is no eigenvector,
in symbols $\ker (A-\lambda\iota)=\{0\}$.
Since $\iota\colon H_2\to H_1$ is compact
and $A\colon H_2\to H_1$ is Fredholm of index zero,
by stability of the Fredholm index under compact perturbation,
$A-\lambda\iota\colon H_2\to H_1$ is as well Fredholm of index zero.
Therefore, since the kernel of $A-\lambda\iota$ is trivial so must be the cokernel.
In particular $A-\lambda\iota\colon H_2\to H_1$ is bijective,
hence by the open mapping theorem the inverse is continuous,
hence $A-\lambda\iota\colon H_2\to H_1$ is an isomorphism.

Now the proof of the lemma is straightforward, indeed
$$
   \xi
   =(A-\lambda\iota)^{-1} (\underbrace{A\xi-\lambda\xi}_{\in H_1}) \in H_2.
$$
\end{proof}

\begin{corollary}[Kernel and invertibility]\label{cor:ker-inv}
Let $A\in\Hh_\gf{h}$.
Then $\ker (A\colon H_1\to H_0)=\ker (A\colon H_2\to H_1)$.
In particular $A\colon H_1\to H_0$ is invertible iff $A\colon H_2\to H_1$ is invertible.
\end{corollary}

\begin{proof}
'$\supset$' is obvious.
'$\subset$' Suppose that $\xi\in H_1$ satisfies $A\xi=0\in H_1$.
Thus $\xi\in H_2$ by Lemma~\ref{le:regularity}.
Concerning invertibility,
since $A$ is Fredholm of index zero,
as a map $H_1\to H_0$ and as a map $H_2\to H_1$,
invertibility is equivalent to triviality of the kernel;
see the second paragraph in the proof of Lemma~\ref{le:regularity}.
\end{proof}

\begin{remark}[Resolvent set non-empty]\label{rem:EV}
The second paragraph of the proof of Lemma~\ref{le:regularity} shows that
the spectrum\footnote{
  $\spec A:=\{\lambda\in\R\mid \text{$A-\lambda\iota\colon
  H_2\to H_1$ is not an isomorphism}\}$
  }
of $A$ only consists of eigenvalues.\footnote{
  Eigenvalues are those $\lambda\in\spec A$ for which $A-\lambda\iota\colon
  H_2\to H_1$ is not injective.
  }
Moreover, in view of the first paragraph of the proof of Lemma~\ref{le:regularity}
the resolvent set $\R\setminus\spec A\not=\emptyset$ is non-empty.
\end{remark}

\begin{corollary}[Spectrum]
Let $A\in\Hh_\gf{h}$.
Then $\spec (A\colon H_2\to H_1)$
is discrete, consists of eigenvalues of finite multiplicity,
and is equal to $\spec (A\colon H_1\to H_0)$.
\\
Moreover, the eigenspaces of $A\colon H_1\to H_0$ and
$A\colon H_2\to H_1$ are equal, in particular multiplicities are equal.
\end{corollary}

\begin{proof}
Choose $\lambda\in \R\setminus \spec A$.
This is possible since the resolvent set of $A\colon H_1\to H_0$
is non-empty by the previous remark.
Since $A$ is $H_0$-symmetric the resolvent operator
$R_\lambda:=\iota\circ (A-\lambda\iota)^{-1}\colon H_0\to H_1\to H_0$
is symmetric as well and it is compact since $\iota$ is.
Therefore, by the spectral theory for symmetric compact operators, the
spectrum of $R_\lambda$ consists of real eigenvalues $\mu$ of finite multiplicity
which only accumulate at the origin.
Observe the equivalence
$$
   R_\lambda\xi
   :=(A-\lambda\iota)^{-1}\xi
   =\mu\xi
   \quad\Leftrightarrow\quad
   A\xi=\tfrac{1+\lambda\mu}{\mu}\xi .
$$
Hence the spectrum of $A\colon H_1\to H_0$
consists of discrete eigenvalues $a_\ell$ of finite multiplicities
which accumulate at $\pm\infty$.

For any $\lambda\in\R$, applying Corollary~\ref{cor:ker-inv}
to $A-\lambda\iota$, we see that the kernels
of $A-\lambda\iota\colon H_1\to H_0$ and of
$A-\lambda\iota\colon H_2\to H_1$ coincide.
In particular, the eigenvalues,
together with their multiplicities, of $A\colon H_1\to H_0$
and of $A\colon H_2\to H_1$ are equal.
\end{proof}

\subsection{Proof of Proposition~\ref{prop:cascasca-intro-NEW}}
\label{sec:Hessian-growth-types}

Suppose $\gt{a},\gt{b}\in\gtset_1$ are shift-invariant
and $A\in (\Ff_\gf{h}^{\gt{a}\gt{b}})^*$.
For any non-eigenvalue $\lambda\in \R\setminus\spec A$
the shifted operator
$
   \A_\lambda:=A-\lambda\iota
$
is invertible.

\medskip\noindent
\textbf{Claim $\mbf{+}$.}
$\gf{f}^+_{\A_\lambda} \sim \gf{f}^+_A$
\hfill{\small\color{gray} equivalently $[\gf{f}^+_{\A_\lambda}]=[\gf{f}^+_A]=\gt{b}$}

\begin{proof}[Proof of Claim $+$]
Enumerate for $\ell\in\Z$ the eigenvalue $a_\ell\in\R$ of $A$ in the form
\begin{equation}\label{eq:A^0-spectrum}
   \dots\le a_{-2}\le a_{-1}<0< a_1\le a_2\le\dots,\qquad
   \Z^*:=\Z\setminus\{0\} .
\end{equation}
Since $\lambda$ is not an eigenvalue, either
$$
   a_{-1}<\lambda<a_1
   \qquad\text{or}\qquad
   a_{\ell_\lambda}<\lambda<a_{{\ell_\lambda}+1}
$$
for some $\ell_\lambda\in \Z^*\setminus \{-1\}=\Z\setminus\{-1,0\}$.
Hence, either
$$
   a_{-1}-\lambda<0<a_1-\lambda=\gf{f}^+_{\A_\lambda}(1)
   \quad\text{or}\quad
   a_{\ell_\lambda}-\lambda < 0<a_{{\ell_\lambda}+1}-\lambda
      =\gf{f}^+_{\A_\lambda}(1) .
$$
So the growth function of the \emph{positive} eigenvalues of
$\A_\lambda=A-\iota\lambda$ is either
\begin{equation}\label{eq:f_A^+-lambda}
   \gf{f}^+_{\A_\lambda}(\nu)
   =a_{\nu}-\lambda
   \qquad\text{or}\qquad
   \gf{f}^+_{\A_\lambda}(\nu)
   =a_{\ell_\lambda +\nu}-\lambda
\end{equation}
where $\nu\in\N$.
Because $[\gf{f}_A^+]_1=[\gf{f}_A^+]$ by Step~1, there exists
$c>0$ such that for each $k\ge 1$ we have
$a_{k+1}\le c a_k$, hence 
\begin{equation}\label{eq:a_k-c^i}
   a_{k+i}\le c^i a_k
\end{equation}
for every $i\in\N_0$.

\smallskip
\noindent
\textsc{Case~1 (${\ell_\lambda}\ge 1$).} So by our enumeration
$0<a_{\ell_\lambda}$ ($<\lambda$). Therefore we estimate
\begin{equation}\label{eq:growth-nu-1}
   \gf{f}^+_{\A_\lambda}(\nu) =a_{\ell_\lambda +\nu}-\lambda
   \le c^{\ell_\lambda} a_\nu -\lambda\le c^{\ell_\lambda} a_\nu
   =c^{\ell_\lambda} \gf{f}^+_A(\nu) .
\end{equation}
Since $\ell_\lambda\ge 1$ we get the first inequality in
\begin{equation}\label{eq:growth-nu-2}
   \gf{f}^+_{A}(\nu)
   =a_\nu
   \le a_{\ell_\lambda +\nu}
   \le \frac{a_{1+\ell_\lambda}}{a_{1+\ell_\lambda}-\lambda} 
   (a_{\ell_\lambda +\nu}-\lambda)
   =\frac{a_{1+\ell_\lambda}}{a_{1+\ell_\lambda}-\lambda} \gf{f}^+_{\A_\lambda}(\nu)
\end{equation}
and the second inequality is obtained by applying Lemma~\ref{le:growth-constant}
to the by $\ell_\lambda$ shifted positive growth function, namely to
$(\gf{f}_A^+)_{\ell_\lambda}$.
\\
The two inequalities~(\ref{eq:growth-nu-1}) and~(\ref{eq:growth-nu-2})
show that $\gf{f}_{\A_\lambda}^+$ and $\gf{f}_A^+$ are of
the same growth type.
This proves the $+$ assertion of Proposition~\ref{prop:cascasca-intro-NEW} in
case $\ell_\lambda\ge 1$.

\smallskip
\noindent
\textsc{Case~2 (${\ell_\lambda}\le -2$).}
In this case $\lambda<a_{\ell_\lambda+1}\le a_{-1}<0$ is negative and
\begin{equation*}
   \gf{f}^+_{\A_\lambda}(\nu)
   =a_{\ell_\lambda +\nu}-\lambda
   \le a_{\nu}-\lambda
   =\tfrac{a_\nu-\lambda}{a_\nu}\, a_\nu
   \le \tfrac{a_1-\lambda}{a_1}\, a_\nu
   =\tfrac{a_1-\lambda}{a_1}\, \gf{f}_A^+(\nu) 
\end{equation*}
for each $\nu\in\N$.
For the second inequality we use that
for $\mu=-\lambda>0$ the function $x\mapsto \frac{x+\mu}{x}$
is monotone decreasing. Indeed its derivative is
$\mu/x^2>0$.
\\
For $\nu\ge 1-\ell_\lambda$ we obtain by~(\ref{eq:a_k-c^i}) the first
inequality in what follows
$$
   \gf{f}_A^+(\nu)
   =a_\nu
   \le c^{-\ell_\lambda} a_{\ell_\lambda +\nu}
   < c^{-\ell_\lambda} (a_{\ell_\lambda +\nu} -\lambda)
   \stackrel{\text{(\ref{eq:f_A^+-lambda})}}{=}
   c^{-\ell_\lambda}\gf{f}_{\A_\lambda}^+(\nu)
$$
where the second inequality holds since $-\lambda>0$.
\\
Since this inequality is true for all $\nu\in\N$, except finitely many
exceptions, there exists a constant $C>0$ such that
$\gf{f}_A^+(\nu)\le C \gf{f}_{\A_\lambda}^+(\nu)$ $\forall\nu\in\N$.

The above two inequalities
show that $\gf{f}_{\A_\lambda}^+$ and $\gf{f}_A^+$ are of
the same growth type.
This proves the $+$ assertion of Proposition~\ref{prop:cascasca-intro-NEW}
in case $\ell_\lambda\le -2$.

\smallskip
\noindent
\textsc{Case~3 ($a_{-1}<\lambda<a_1$).}
In this case $\gf{f}_{\A_\lambda}^+(\nu)=a_\nu-\lambda
=\gf{f}_A^+(\nu)-\lambda$.
Hence the two growth functions just differ by a constant and therefore,
by Lemma~\ref{le:growth-constant},
they have the same growth type.

\smallskip
\noindent
This proves Claim $+$, so
Proposition~\ref{prop:cascasca-intro-NEW}
in the $(+)$ case, namely $[\gf{f}^+_{\A_\lambda}]=\gt{b}$.
\end{proof}

\noindent
The $(-)$ case of Proposition~\ref{prop:cascasca-intro-NEW}
follows from the $(+)$ case by replacing $A$ by~$-A$.
The proof of Proposition~\ref{prop:cascasca-intro-NEW} is complete.

\section{Proof of Theorem~\ref{thm:A}}
\label{sec:thmA}

We prove Theorem~\ref{thm:A} in three steps.

\smallskip
\noindent
\textbf{Step~1.} 
Given growth types $\gt{a},\gt{b}\in\gtset$, shift invariant or not,
and $A\in (\Ff_\gf{h}^{\gt{a}\gt{b}})^*$, in particular $A$ is invertible.
Then there exists an open neighborhood $\Vv_A$ of $A$ in $\Ff_\gf{h}^*$ which lies in
$(\Ff_\gf{h}^{\gt{a}\gt{b}})^*$, in symbols $\Vv_A\subset
(\Ff_\gf{h}^{\gt{a}\gt{b}})^*\subset \Ff_\gf{h}$.

\begin{proof}[Proof of Step~1]
In Case~1 shift invariance is not required.
By Theorem~\ref{thm:B}, there is an open neighborhood $\Uu_A$ of
$A$ in $\Ff_\gf{h}^*$ such that both maps
$$
    \Pi_\pm\colon\Uu_A\to \Ll(H_\frac{1}{2})\cap\Ll(H_\frac{3}{2})
   ,\quad
   B\mapsto \Pi_\pm(B)=:\Pi_\pm^B
$$
are continuous.
Hence, given $\eps>0$, there exists an open neighborhood
$\Uu^+_\eps(A)\subset\Uu_A$ of $A$ such that
$$
   \Norm{\Pi^B_+-\Pi^A_+}_{\Ll(H_\frac{1}{2})\cap\Ll(H_\frac{3}{2})}\le\eps
$$
for every $B\in\Uu_+^\eps(A)$.
We fix
$
   \eps_+:=\min\{\tfrac12,
   1/4\norm{\Pi^A_+}_{\Ll(H_\frac{1}{2})\cap\Ll(H_\frac{3}{2})}\}
$
and set
$$
   \Vv_A^+:=\Uu^+_{\eps_+}(A) .
$$

\smallskip\noindent
\textbf{Claim \boldmath$+$.}
For every $B\in \Vv_A^+$ it holds that $\gf{f}^+_B\sim
\gf{f}^+_A$.

\smallskip\noindent
\textsc{Proof of Claim~$+$.}
Pick $B\in \Vv_A^+$.
We restrict the projections $\Pi_+^B,\Pi_+^A\colon H_\frac12\to H_\frac12$
and consider them as maps between the following subspaces
$$
   \Pi_+^B|_A:=
   \Pi_+^B|_{\Pi_+^A H_\frac12}\colon
   \Pi_+^A H_\frac12 \to \Pi_+^B H_\frac12
$$
and
$$
   \Pi_+^A|_B:=
   \Pi_+^A|_{\Pi_+^B H_\frac12}\colon
   \Pi_+^B H_\frac12 \to \Pi_+^A H_\frac12.
$$
We prove that the following compositions are isomorphisms
\begin{equation}\label{eq:a)b)}
   \text{a)
   $\Pi_+^A|_B\circ\Pi_+^B|_A
   \in\Ll(\Pi_+^A  H_\frac12)$}
   ,\qquad
   \text{b) 
   $\Pi_+^B|_A\circ\Pi_+^A|_B
   \in\Ll(\Pi_+^B H_\frac12)$.}
\end{equation}
This is true for $B=A$ since
the image of the projection $\Pi_+^A $ is its
fixed point set, in symbols
$\Pi_+^A H_\frac12=\Fix\, \Pi_+^A $.
Hence $\Pi_+^A|_A=\Id$ on $\Pi_+^A  H_\frac12$.

To see that the above two compositions are invertible
we show that they differ in norm from the identity by less than $1$.
Indeed
\begin{equation*}
\begin{split}
   &\Norm{\Pi_+^A|_B\circ\Pi_+^B|_A-\Id_{\Pi_+^A  H_\frac12}}
   _{\Ll(H_\frac12)\cap \Ll(H_{\frac{3}{2}})}
\\
   &=\Norm{\Pi_+^A|_B\circ\Pi_+^B|_A-\Pi_+^A |_A\circ\Pi_+^A |_A}
   _{\Ll(H_\frac12)\cap \Ll(H_{\frac{3}{2}})}
\\
   &=\Norm{\Pi_+^A \left(\Pi_+^B|_A-\Pi_+^A |_A\right)}
   _{\Ll(H_\frac12)\cap \Ll(H_{\frac{3}{2}})}
\\
   &\le 
   \Norm{\Pi_+^A }
   _{\Ll(H_\frac12)\cap \Ll(H_{\frac{3}{2}})}
   \Norm{\Pi_+^B|_A-\Pi_+^A |_A}
   _{\Ll(H_\frac12)\cap \Ll(H_{\frac{3}{2}})}
\\
   &\le 
   \Norm{\Pi_+^A }
   _{\Ll(H_\frac12)\cap \Ll(H_{\frac{3}{2}})}
   \Norm{\Pi_+^B-\Pi_+^A }
   _{\Ll(H_\frac12)\cap \Ll(H_{\frac{3}{2}})}
\\
   &\le \Norm{\Pi_+^A }_{\Ll(H_\frac12)\cap \Ll(H_{\frac{3}{2}})}\cdot \eps_+ 
   \le\tfrac{1}{4} 
\end{split}
\end{equation*}
and
\begin{equation*}
\begin{split}
   &\Norm{\Pi_+^B|_A\circ\Pi_+^A|_B-\Id_{\Pi_+^B H_\frac12}}
   _{\Ll(H_\frac12)\cap \Ll(H_{\frac{3}{2}})}
\\
   &=\Norm{\Pi_+^B|_A\circ\Pi_+^A |_B-\Pi_+^B|_B\circ\Pi_+^B|_B}
   _{\Ll(H_\frac12)\cap \Ll(H_{\frac{3}{2}})}
\\
   &=\Norm{\Pi_+^B\left(\Pi_+^A |_B-\Pi_+^B|_B\right)}
   _{\Ll(H_\frac12)\cap \Ll(H_{\frac{3}{2}})}
\\
   &\le 
   \Norm{\Pi_+^B}
   _{\Ll(H_\frac12)\cap \Ll(H_{\frac{3}{2}})}
   \Norm{\Pi_+^A |_B-\Pi_+^B|_B}
   _{\Ll(H_\frac12)\cap \Ll(H_{\frac{3}{2}})}
\\
   &\le 
   \Norm{\Pi_+^B}
   _{\Ll(H_\frac12)\cap \Ll(H_{\frac{3}{2}})}
   \Norm{\Pi_+^A -\Pi_+^B}
   _{\Ll(H_\frac12)\cap \Ll(H_{\frac{3}{2}})}
\\
   &\le \Norm{\Pi_+^B{\color{gray}-\,\Pi_+^A +\Pi_+^A}}_{\Ll(H_\frac12)\cap\Ll(H_{\frac{3}{2}})}
   \cdot \eps_+
\\
   &\le \eps_+^2+\eps_+ \Norm{\Pi_+^A }
   _{\Ll(H_\frac12)\cap \Ll(H_{\frac{3}{2}})}
   \le\tfrac{1}{4}+\tfrac{1}{4}=\tfrac{1}{2} .
\end{split}
\end{equation*}
This proves that both compositions 
a) and b) in~(\ref{eq:a)b)})
are invertible.
Therefore $\Pi_+^B|_A\colon \Pi_+^A  H_\frac12\to \Pi_+^B H_\frac12$
is a) injective and b) surjective, hence an isomorphism.
\\
By the same argument
$\Pi_+^B|_A\colon \Pi_+^A  H_{\frac{3}{2}}\to \Pi_+^B
H_{\frac{3}{2}}$ is an isomorphism.

\smallskip
\noindent
\textsc{Conclusion.} 
For any $B$ in the open neighborhood $\Vv_A^+$ of $A$
in $\Ff_\gf{h}$ the restriction
$$
   \Pi_+^B|_A\colon
   (\Pi_+^A  H_{\frac{1}{2}}, \Pi_+^A  H_{\frac{3}{2}})
   \to (\Pi_+^B H_{\frac{1}{2}}, \Pi_+^B H_{\frac{3}{2}})
$$
is an isomorphism of Hilbert pairs.
\\
Hence, by Theorem~\ref{thm:f-g-equiv},
the two Hilbert pairs have the same growth type.
The growth type of the pair $(\Pi_+^A H_0, \Pi_+^A H_1)$
is $[(\gf{f}^+_A)^2]$ and by shift-invariance of scales,
see Example~\ref{ex:h-pair},
the same is true for the growth type of $(\Pi_+^A  H_{\frac{1}{2}},\Pi_+^A
H_{\frac{3}{2}})$.
By the same reasoning the growth type of $(\Pi_+^B  H_{\frac{1}{2}},\Pi_+^B
H_{\frac{3}{2}})$ is $[(\gf{f}^+_B)^2]$.

Hence $(\gf{f}^+_A)^2\sim(\gf{f}^+_B)^2$ and therefore $\gf{f}^+_A\sim \gf{f}^+_B$.
This proves Claim~$+$.

\smallskip\noindent
\textbf{Claim \boldmath$-$.}
There exists an open neighborhood $\Vv_A^-$ of $A$ in $\Ff_\gf{h}$
such that for every $B\in \Vv_A^-$ it holds that $\gf{f}^-_A\sim
\gf{f}^-_B$.

\smallskip\noindent
\textsc{Proof of Claim~$-$.}
Same argument as in Claim~$+$.

\medskip 
Define $\Vv_A:=\Vv_A^+\cap \Vv_A^-$.
Then $\Vv_A$ is an open neighborhood of $A$ in $\Ff_\gf{h}^*$
and for every $B\in\Vv_A$ we have 
$[\gf{f}^+_B]=[\gf{f}^+_A]=\gt{b}$
and $[\gf{f}^-_B]=[\gf{f}^-_A]=\gt{a}$.
Thus $B$ belongs to $(\Ff_\gf{h}^{\gt{a}\gt{b}})^*$.
This proves Step~1.
\end{proof} 

\noindent
\textbf{Step~2.}
Given shift invariant growth types $\gt{a},\gt{b}\in\gtset_1$,
then $\Ff_\gf{h}^{\gt{a}\gt{b}}$ is open in $\Ff_\gf{h}$.

\begin{proof}[Proof of Step~2]
Pick any non-eigenvalue $\lambda\in \R\setminus\spec A$;
see Remark~\ref{rem:EV}.
Then the (invertible) operator $\A_\lambda:=A-\lambda\iota$ 
lies in $(\Ff_\gf{h}^{\gt{a}\gt{b}})^*$,
by Proposition~\ref{prop:cascasca-intro-NEW}.
Hence Case~2 follows by applying Case~1 to $\A_\lambda$ as we explain next.
By Case~1 we get an open neighborhood
$\Vv_{\A_\lambda}\subset \Ff_\gf{h}^{\gt{a}\gt{b}}$
of $\A_\lambda$ in $\Ff_\gf{h}$.
The translated set $\Vv_A:=\Vv_{\A_\lambda}+\lambda\iota$ is still
open and contains $A$. 
By translation invariance~(\ref{eq:translation-invariance}) of
$\Ff_\gf{h}^{\gt{a}\gt{b}}$ the set $\Vv_A$ is an open neighborhood of
$A$ in $\Ff_\gf{h}^{\gt{a}\gt{b}}$. This proves Step~2.
\end{proof} 

Before we continue the proof to show closedness we introduce some notation.

\begin{definition}[Wild Floer Hessians]
A Floer Hessian $A\in\Ff_\gf{h}$ is called \textbf{wild} 
if there exist $\lambda_1,\lambda_2\in \R\setminus\spec A$ such that
\begin{equation}\label{eq:gtuy77h}
   \left([\gf{f}^+_{\A_{\lambda_1}}],[\gf{f}^-_{\A_{\lambda_1}}]\right)
   \not=
   \left([\gf{f}^+_{\A_{\lambda_2}}],[\gf{f}^-_{\A_{\lambda_2}}]\right)
   \in\gtset\times\gtset .
\end{equation}
Equivalently, by Proposition~\ref{prop:cascasca-intro-NEW},
a Floer Hessian $A$ is wild if for one, hence for every, $\lambda\in
\R\setminus\spec A$ at least one of the growth types
$[\gf{f}^+_{\A_{\lambda}}]$ or $[\gf{f}^-_{\A_{\lambda}}]$ is not shift invariant.
We call $A\in\Ff_\gf{h}$ \textbf{tame} if it is not wild.
\end{definition}

Therefore we have the disjoint union in wild and tame Floer Hessians
$$
   \Ff_\gf{h}=\Ww_\gf{h} \cup \Tt_\gf{h}
$$
where $\Ww_\gf{h}=\{A\in \Ff_\gf{h}\mid\text{$A$ wild}\}$ and
$\Tt_\gf{h}=\{A\in \Ff_\gf{h}\mid\text{$A$ tame}\}$.
Note that $A$ is tame iff there exist shift-invariant growth types
$\gt{a}$ and $\gt{b}$ such that $A\in \Ff_\gf{h}^{\gt{a}\gt{b}}$.

\begin{lemma}\label{le:wild-open}
The set of wild Floer Hessians $\Ww_\gf{h}$ is open in $\Ff_\gf{h}$.
\end{lemma}

\begin{proof}
Suppose that $A\in \Ww_\gf{h}$.
Then there exist $\lambda_1,\lambda_2\in \R\setminus\spec A$
such that~(\ref{eq:gtuy77h}) holds true. By definition~(\ref{eq:F^*})
we have that
$$
   \A_{\lambda_1}
   \in \left(\Ff_\gf{f}^{[\gf{f}^-_{\A_{\lambda_1}}]\,[\gf{f}^+_{\A_{\lambda_1}}]}\right)^*
$$
and so Step~1 applies:
there is an open neighborhood $\Vv_{\A_{\lambda_1}}$ of
$\A_{\lambda_1}$ in $\Ff_\gf{h}$ with 
$$
   \Vv_{\A_{\lambda_1}}\subset
   \left(\Ff_\gf{f}^{[\gf{f}^-_{\A_{\lambda_1}}]\,[\gf{f}^+_{\A_{\lambda_1}}]}\right)^*.
$$
Similarly there exists an open neighborhood $\Vv_{\A_{\lambda_2}}$ of
$\A_{\lambda_2}$ in $\Ff_\gf{h}$ with 
$$
   \Vv_{\A_{\lambda_2}}\subset
   \left(\Ff_\gf{f}^{[\gf{f}^-_{\A_{\lambda_2}}]\,[\gf{f}^+_{\A_{\lambda_2}}]}\right)^*.
$$
For $\lambda\in\R$ consider the translation
$$
   T_\lambda\colon \Ff_\gf{h}\to \Ff_\gf{h}
   ,\quad
   B\mapsto B-\lambda\iota=:\B_\lambda .
$$
Note that $A=T_{-\lambda_1}\A_{\lambda_1}$  and $A=T_{-\lambda_2}\A_{\lambda_2}$. So
$A\in T_{-\lambda_1}\Vv_{\A_{\lambda_1}}\cap T_{-\lambda_2}\Vv_{\A_{\lambda_2}}=:\Vv_A$
and $\Vv_A$ is open in $\Ff_\gf{h}$.

It remains to show that $\Vv_A$ is contained in $\Ww_\gf{h}$.
To see this pick $B\in\Vv_A$.
Then $T_{\lambda_1} B\in \Vv_{\A_{\lambda_1}}\subset 
\left(\Ff_\gf{f}^{[\gf{f}^-_{\A_{\lambda_1}}]\,[\gf{f}^+_{\A_{\lambda_1}}]}\right)^*$.
So $([\gf{f}^-_{T_{\lambda_1} B}], [\gf{f}^+_{T_{\lambda_1} B}])
=([\gf{f}^-_{\A_{\lambda_1}}],[\gf{f}^+_{\A_{\lambda_1}}])$.
Similarly for $\lambda_2$.
Therefore, together with~(\ref{eq:gtuy77h}), it follows that
$([\gf{f}^-_{T_{\lambda_1} B}], [\gf{f}^+_{T_{\lambda_1} B}])
\not=([\gf{f}^-_{T_{\lambda_2} B}], [\gf{f}^+_{T_{\lambda_2} B}])$
which means that $B$ is wild.
This proves Lemma~\ref{le:wild-open}.
\end{proof}

\begin{corollary}\label{cor:Tt_h-closed}
The set of tame Floer Hessians $\Tt_\gf{h}$ is closed in $\Ff_\gf{h}$.
\end{corollary}

We are now in position to finish the proof of the theorem.

\smallskip
\noindent
\textbf{Step~3.}
Given shift invariant growth types $\gt{a},\gt{b}\in\gtset_1$,
then $\Ff_\gf{h}^{\gt{a}\gt{b}}$ is closed in $\Ff_\gf{h}$.

\begin{proof}[Proof of Step~3]
Since $\Tt_\gf{h}$ is closed in $\Ff_\gf{h}$, it suffices to show that
$\Ff_\gf{h}^{\gt{a}\gt{b}}$ is closed in $\Tt_\gf{h}$.
We write
$$
   \Tt_\gf{h}=\bigcup_{(\gt{c},\gt{d})\in\atop\gtset_1\times\gtset_1}
   \Ff_\gf{h}^{\gt{c}\gt{d}}.
$$
Since for every pair $(\gt{c},\gt{d})\in\gtset_1\times\gtset_1$
the subset $\Ff_\gf{h}^{\gt{c}\gt{d}}$ is open in $\Tt_\gf{h}$ by Step~2, it follows
that the complement of $\Ff_\gf{h}^{\gt{a}\gt{b}}$, namely
$$
   \Tt_\gf{h}\setminus \Ff_\gf{h}^{\gt{a}\gt{b}}
   =\bigcup_{(\gt{c},\gt{d})\not=(\gt{a},\gt{b})\in\atop\gtset_1\times\gtset_1}
   \Ff_\gf{h}^{\gt{c}\gt{d}}.
$$
is open in $\Tt_\gf{h}$ and therefore $\Ff_\gf{h}^{\gt{a}\gt{b}}$ is
closed in $\Tt_\gf{h}$.
This proves Step~3.
\end{proof} 

This concludes the proof of Theorem~\ref{thm:A}.

\section{Schur multipliers}\label{sec:Schur}

Consider the Hilbert space $\ell^2=\ell^2_\N(\R)$
of square summable real sequences enumerated by the positive integers.
Let $\{e_\nu\}_{\nu\in\N}$ be the standard basis of $\ell^2$.
Let $\Ll(\ell^2)$ be the vector space of bounded linear operators on $\ell^2$.
Let $\ell^\infty_{\N\times\N}=\ell^\infty_{\N\times\N}(\R)$ the
vector space of (infinite) real matrices.

\begin{definition}\label{def:Schur}
A bounded linear operator $\mat{a}\in\Ll(\ell^2)$
we identify with its (infinite) matrix $(a_{\mu\nu})=(a_{\mu\nu})_{\mu,\nu}$ with entries
$a_{\mu\nu}:=\INNER{e_\mu}{\mat{a} e_\nu}$ for $\mu,\nu\in\N$.
Thus we write $\mat{a}=(a_{\mu\nu})$.
A matrix $\mat{b}=(b_{\mu\nu})$ is called a \textbf{Schur multiplier}
if it gives rise to a continuous (bounded) linear operator on $\Ll(\ell^2)$ defined by
$$
   S_{\mat{b}}\colon \Ll(\ell^2)\to\Ll(\ell^2) ,\quad
   (a_{\mu\nu})\mapsto (b_{\mu\nu} a_{\mu\nu})=:\mat{b}\odot\mat{a} .
$$
The operator $S_{\mat{b}}\in\Ll(\Ll(\ell^2))$ is called a \textbf{Schur operator}.

\smallskip
The set $\Ss$ of all Schur operators $S_{\mat{b}}$ is a vector subspace of
$\Ll(\Ll(\ell^2))$. It is endowed with the
commutative\footnote{
  Since multiplication in $\R$ is associative it holds
  $S_{\mat{b}}\circ S_{\mat{\tilde b}}=S_{\mat{b}\odot\mat{\tilde b}}$
  and since it is commutative it holds
  $\mat{b}\odot\mat{\tilde b}=\mat{\tilde b}\odot\mat{b}$.
  }
product
\begin{equation*}
\begin{split}
   \Ss\times\Ss&\to\Ss 
\\
   (S_{\mat{b}}, S_{\mat{\tilde b}})
   &\mapsto S_{\mat{b}}\circ S_{\mat{\tilde b}}
   =S_{\mat{b}\odot\mat{\tilde b}} .
\end{split}
\end{equation*}
\end{definition}

The following lemma is due to Schur.

\begin{lemma}[Schur~{\cite[Thm.\,VI]{Schur:1911a}}]\label{le:Schur}
Let $I\subset\R$ be an interval and $\kappa$ a constant.
Suppose $f_\nu, g_\mu\in L^2(I,\R)$ are two uniformly bounded function sequences
\begin{equation}\label{eq:Schur-hyp}
   \norm{f_\nu}_{L^2(I)}^2\le \kappa
   ,\qquad
  \norm{g_\mu}_{L^2(I)}^2\le \kappa ,
\end{equation}
for all $\mu,\nu\in\N$. Then the matrix entries
$$
   b_{\mu\nu}:=\int_I g_\mu(s) f_\nu(s)\, ds
$$
define a Schur multiplier $\mat{b}=(b_{\mu\nu})$ whose Schur operator satisfies 
$$
   \norm{S_\mat{b}}_{\Ll(\Ll(\ell^2))}
   \le \kappa.
$$
\end{lemma}

\begin{proof}
We follow the proof of Schur~\cite[p.\,13]{Schur:1911a}.
Let $\mat{a}\in\Ll(\ell^2)$. For unit vectors $x,y\in\ell^2$ we get
\begin{equation*}
\begin{split}
   \Abs{\INNER{y}{S_{\mat{b}}\mat{a} x}}
   &=\Abs{\sum_{\mu,\nu} y_\mu b_{\mu\nu}a_{\mu\nu}x_\nu}\\
   &\stackrel{{\color{gray} 1}}{\le}
   \int_I\Abs{\sum_{\mu,\nu} y_\mu g_\mu(s) a_{\mu\nu}x_\nu f_\nu(s)}\, ds\\
   &\stackrel{{\color{gray} 2}}{\le}
   \int_I \norm{\mat{a}}_{\Ll(\ell^2)}\sqrt{\sum_\mu y_\mu^2 g_\mu(s)^2}
   \sqrt{\sum_\nu x_\nu^2 f_\nu(s)^2}\, ds\\
   &\stackrel{{\color{gray} 3}}{\le}
   \frac12 \norm{\mat{a}}_{\Ll(\ell^2)}
   \int_I\left(\sum_\mu y_\mu^2 g_\mu(s)^2+\sum_\mu x_\mu^2
     f_\mu(s)^2\right) ds\\
   &=\frac12 \norm{\mat{a}}_{\Ll(\ell^2)}\left(
   \sum_\mu y_\mu^2 \int_I g_\mu(s)^2\, ds
   +\sum_\mu x_\mu^2 \int_I f_\mu(s)^2\, ds
    \right) \\
   &\stackrel{{\color{gray} 4}}{\le}
   \frac{\kappa}{2} \norm{\mat{a}}_{\Ll(\ell^2)}
   \left(\sum_\mu y_\mu^2 +\sum_\mu x_\mu^2 \right)\\
   &=\kappa \norm{\mat{a}}_{\Ll(\ell^2)} .
\end{split}
\end{equation*}
Here inequality {\color{gray} 1} is by definition of $b_{\mu\nu}$.
Inequality {\color{gray} 2} is by definition of the operator
norm of $\mat{a}$ applied to the $\ell^2$ elements $(y_\mu g_\mu(s))$
and $(x_\nu f_\nu(s))$. Inequality {\color{gray} 3} uses that $rs\le (r^2+s^2)/2$.
Inequality {\color{gray} 4}  uses hypothesis~(\ref{eq:Schur-hyp}).
The final identity holds since $x,y\in\ell^2$ are unit vectors.

The above estimate implies that
$\norm{S_{\mat{b}}\mat{a}}_{\Ll(\ell^2)} \le\kappa \norm{\mat{a}}_{\Ll(\ell^2)}$
and therefore we obtain $\norm{S_{\mat{b}}}_{\Ll(\Ll(\ell^2))} \le\kappa$.
\end{proof}

The following corollary essentially follows the argument of
Schur~{\cite[p.\,14]{Schur:1911a} when he showed that the Hilbert
matrix $\mat{b}=(\frac{1}{\mu+\nu})$ is a Schur multiplier.

\begin{corollary}\label{cor:Schur-m}
Suppose that $a_\mu$, $b_\nu$ are sequences of positive reals,
then
$$
   b_{\mu\nu}:=\frac{\sqrt{a_\mu b_\nu}}{a_\mu+b_\nu}
$$
defines a Schur multiplier $\mat{b}=(b_{\mu\nu})$ such that
$\norm{S_{\mat{b}}}_{\Ll(\Ll(\ell^2))}\le\frac12$.
\end{corollary}

\begin{proof}
Consider the interval $I=[0,\infty)$ and the functions
$f_\nu(s):=\sqrt{b_\nu} e^{-b_\nu s}$ and $g_\mu(s):=\sqrt{a_\mu} e^{-a_\mu s}$.
Then calculation shows that
$$
   \norm{f_\nu}^2_{L^2(I)}=\tfrac12, \quad
   \norm{g_\mu}^2_{L^2(I)}=\tfrac12, \quad
   b_{\mu\nu}:=\int_I f_\nu(s) g_\mu(s)\, ds=\tfrac{\sqrt{a_\mu b_\nu}}{a_\mu+b_\nu}.
$$
Now the Schur Lemma~\ref{le:Schur} applies with $\kappa=\frac12$.
\end{proof}

\begin{corollary}\label{cor:Schur-abstr}
Suppose there exists a separable
Hilbert space $H$, sequences of vectors $f_\nu, g_\mu \in H$,
and a constant $\kappa$ such that
\begin{equation*}
   \norm{f_\nu}_H^2\le \kappa
   ,\qquad
  \norm{g_\mu}_H^2\le \kappa ,
\end{equation*}
for all $\nu,\mu\in\N$.
Then the matrix entries 
$$
   b_{\mu\nu}:=\INNER{g_\mu}{f_\nu}_H
$$
define a Schur multiplier $\mat{b}=(b_{\mu\nu})$ whose Schur operator satisfies 
$$
   \norm{S_\mat{b}}_{\Ll(\Ll(\ell^2))}
   \le \kappa.
$$
\end{corollary}

\begin{proof}
An infinite dimensional separable Hilbert space $H$
is isometric to $L^2(I,\R)$ for any interval $I=[a,b]$ with $a<b$.
Now Lemma~\ref{le:Schur} applies.
If $H$ has finite dimension, embed it into an infinite
dimensional separable Hilbert space.
\end{proof}

\begin{corollary}\label{cor:Schur-abstr-2}
Suppose there exists a separable
Hilbert space $H$, sequences of vectors $f_\nu, g_\mu \in H$,
and constants $\kappa_1$ and $\kappa_2$ such that
\begin{equation*}
   \norm{f_\nu}_H\le \kappa_1
   ,\qquad
  \norm{g_\mu}_H\le \kappa_2 ,
\end{equation*}
for all $\nu,\mu\in\N$.
Then the matrix entries 
$$
   b_{\mu\nu}:=\INNER{g_\mu}{f_\nu}_H
$$
define a Schur multiplier $\mat{b}=(b_{\mu\nu})$ whose Schur operator satisfies 
$$
   \norm{S_\mat{b}}_{\Ll(\Ll(\ell^2))}
   \le \kappa_1\kappa_2.
$$
\end{corollary}

\begin{proof}
The rescaled sequences $F_\nu:=\sqrt{\kappa_2/\kappa_1} f_\nu$ and
$G_\mu:=\sqrt{\kappa_1/\kappa_2} g_\mu$ satisfy
$$
   \norm{F_\nu}_H^2\le \kappa_1\kappa_2
   ,\qquad
  \norm{g_\mu}_H^2\le \kappa_1\kappa_2 ,
$$
and the matrix entries are the same
$$
   B_{\mu\nu}:=\INNER{G_\mu}{F_\nu}_H
   =\INNER{g_\mu}{f_\nu}_H
   =:b_{\mu\nu}.
$$
Now apply Corollary~\ref{cor:Schur-abstr} to the sequences $F_\nu$ and $G_\mu$.
\end{proof}

\begin{example}\label{ex:Schur}
Let $a_\mu$ and $b_\nu$ be $\N$-sequences of positive reals. Let $H=\R$.

\smallskip
I.~For $f_\nu:=1$ and $g_\mu:=\arctan a_\mu$ we have
$\abs{f_\nu}_\R=1$ and $\abs{g_\mu}_\R\le \pi/2$. Thus, by
Corollary~\ref{cor:Schur-abstr-2}, the matrix entries
$b_{\mu\nu}^{\rm I}:=g_\mu f_\nu=\arctan a_\mu$
define a Schur multiplier $\mat{b}^{\rm I}=(b_{\mu\nu}^{\rm I})$
whose Schur operator satisfies $\norm{S_{\mat{b}^{\rm I}}}_{\Ll(\Ll(\ell^2))}\le\pi/2$.

\smallskip
II.~For $f_\nu:=\arctan b_\nu$ and $g_\mu:=1$ we have
$\abs{f_\nu}_\R\le\pi/2$ and $\abs{g_\mu}_\R=1$. Thus, by
Corollary~\ref{cor:Schur-abstr-2}, the matrix entries
$b_{\mu\nu}^{\rm I}:=g_\mu f_\nu=\arctan b_\nu$
define a Schur multiplier $\mat{b}^{\rm II}=(b_{\mu\nu}^{\rm II})$
whose Schur operator satisfies $\norm{S_{\mat{b}^{\rm II}}}_{\Ll(\Ll(\ell^2))}\le\pi/2$.

\smallskip
III.~Combining I and II with Corollary~\ref{cor:Schur-m}
and using that Schur multipliers form an algebra
we conclude that the matrix entries
$$
   \tilde b_{\mu\nu}
   :=\frac{\sqrt{a_\mu b_\nu}}{a_\mu+b_\nu}\left(\arctan a_\mu+\arctan b_\nu\right)
$$
define a Schur multiplier $\mat{\tilde b}=(\tilde b_{\mu\nu})$
and $\norm{S_{\mat{\tilde b}}}_{\Ll(\Ll(\ell^2))}\le\pi/2$.
\end{example}

For the proofs of Theorem~\ref{thm:A} and Theorem~\ref{thm:B}
Corollary~\ref{cor:Schur-abstr-2} is sufficient.
The reason that we cannot, with the method of Schur multipliers,
improve Theorem~\ref{thm:B} from $H_{\frac12}$ to $H_0$ lies in the
fact that there is a converse to Corollary~\ref{cor:Schur-abstr-2}
which then gives obstructions for matrices to be a Schur multiplier.

Indeed Grothendieck~\cite{Grothendieck:1956a}
proved, see also~\cite[Thm.\,5.1]{Pisier:2001a},
the converse

\begin{theorem}\label{thm:Grothendieck-Schur}
A matrix $\mat{b}=(b_{\mu\nu})$ is a Schur multiplier
whose Schur operator satisfies 
$$
   \norm{S_\mat{b}}_{\Ll(\Ll(\ell^2))}
   \le \kappa
$$
if and only if there exists a Hilbert space $H$ and
sequences $f_\nu, g_\mu \in H$ satisfying
\begin{equation*}
   \norm{f_\nu}_H\le \kappa_1
   ,\qquad
  \norm{g_\mu}_H\le \kappa_2
   ,\qquad
   \forall \nu,\mu\in\N ,
\end{equation*}
such that $b_{\mu\nu}=\INNER{g_\mu}{f_\nu}_H$
for all $\nu,\mu\in\N$ and $\kappa=\kappa_1\kappa_2$.
\end{theorem}

The following corollary of Theorem~\ref{thm:Grothendieck-Schur}
was explained to us by Gilles Pisier.

\begin{corollary}\label{cor:Schur-criterium}
Let $\mat{b}=(b_{\mu\nu})$ be a Schur multiplier
such that both limits exist
$$
   \lim_{\mu\to\infty} \lim_{\nu\to\infty} b_{\mu\nu}
   ,\qquad
   \lim_{\nu\to\infty}\lim_{\mu\to\infty} b_{\mu\nu} .
$$
Then these two limits are equal.
\end{corollary}

\begin{proof}
By Theorem~\ref{thm:Grothendieck-Schur}
there exist a Hilbert space $H$ and
bounded sequences $f_\nu, g_\mu \in H$
such that $b_{\mu\nu}=\INNER{g_\mu}{f_\nu}_H$
for all $\nu,\mu\in\N$.
\\
By the Banach-Alaoglu Theorem
bounded sequences in Hilbert space are weakly compact.
Hence there exist $f,g\in H$ and subsequences
$f_{\nu_j}$ and $g_{\mu_j}$ such that
$$
   \lim_{j\to\infty}\INNER{f_{\nu_j}}{h}=\INNER{f}{h}
   ,\qquad
   \lim_{j\to\infty}\INNER{g_{\mu_j}}{h}=\INNER{g}{h} ,
$$
for every $h\in H$.
Computing both limits
$$
   \lim_{\mu\to\infty}\lim_{\nu\to\infty} b_{\mu\nu}
   =\lim_{i\to\infty}\lim_{j\to\infty} b_{\mu_i\nu_j}
   =\lim_{i\to\infty}
   \underbrace{\lim_{j\to\infty} 
      \INNER{g_{\mu_i}}{f_{\nu_j}}}
   _{\stackrel{h=g_{\mu_i}}{=}\INNER{g_{\mu_i}}{f}}
   \stackrel{h=f}{=}\INNER{g}{f}
$$
and
$$
   \lim_{\nu\to\infty}\lim_{\mu\to\infty} b_{\mu\nu}
   =\lim_{j\to\infty}\lim_{i\to\infty} b_{\mu_i\nu_j}
   =\lim_{j\to\infty}
   \underbrace{\lim_{i\to\infty} 
      \INNER{g_{\mu_i}}{f_{\nu_j}}}
   _{\stackrel{h=f_{\nu_j}}{=}\INNER{g}{f_{\nu_j}}}
   \stackrel{h=g}{=}\INNER{g}{f}
$$
shows that they are equal.
\end{proof}

\begin{corollary}\label{cor:obstruction}
The matrix $\mat{b}=(\tfrac{\mu}{\mu+\nu})$ is not a Schur multiplier.
\end{corollary}

\begin{proof}
The limits
$$
   \lim_{\mu\to\infty} \lim_{\nu\to\infty} \tfrac{\mu}{\mu+\nu}
   =\lim_{\mu\to\infty} 0=0
   ,\qquad
   \lim_{\nu\to\infty}\lim_{\mu\to\infty} \tfrac{\mu}{\mu+\nu}
   =\lim_{\nu\to\infty} 1=1
$$
are different. So $\mat{b}$ it is not a Schur multiplier
by Corollary~\ref{cor:Schur-criterium}.
\end{proof}

For the significance of Corollary~\ref{cor:obstruction}
see Remark~\ref{rmk:obstruction}.

\section{Proof of Theorem~\ref{thm:B}}
\label{sec:thmB}

\subsection{Spectral projections}\label{sec:spectral-proj}

In the present Section~\ref{sec:thmB}
we make the following assumptions.
By $A\colon H_1\to H_0$ we denote an $H_0$-symmetric isomorphism whose
spectrum $\spec A$ consists of infinitely many positive and infinitely
many negative eigenvalues, of finite multiplicity each.
Let the eigenvalues $a_\ell$ of $A$, where $\ell\in\Z^*$, be enumerated as
in~(\ref{eq:A^0-spectrum}), so $a_1$ is the smallest positive
eigenvalue of $A$ and $a_{-1}$ is the largest negative one.
Let $\sigma=\sigma(A):=\min\{a_1,-a_{-1}\}>0$
be the \textbf{spectral gap} of $A$. 

Furthermore, in Section~\ref{sec:spectral-proj}
-- {\color{red}to avoid annoying} constants -- we
endow the space $H_1$ with the (since
$A$ is an isomorphism) equivalent inner product
\begin{equation}\label{eq:A-norm}
   \INNER{\cdot}{\cdot}_A:=\INNER{A\cdot}{A\cdot}_0\sim\INNER{\cdot}{\cdot}_1.
\end{equation}
%
%
Because the inclusion $H_1\INTO H_0$
is compact, the resolvent $R\colon H_0\to H_0$ of $A$ is compact and symmetric.
Hence there exists an orthonormal basis of $H_0$ composed
of eigenvectors of $A\colon H_1\to H_0$, in symbols
\begin{equation}\label{eq:ONB-Vv-H_0}
   \text{$\Vv=\Vv(A)=(v_\ell)_{\ell\in\Z^*}$ ONB of $H_0$}
   ,\quad
   v_\ell\in H_1
   ,\quad
   A v_\ell =a_\ell v_\ell .
\end{equation}
Since $\abs{v_\ell}_1^2=\INNER{v_\ell}{v_\ell}_1
\stackrel{\text{(\ref{eq:A-norm})}}{=}
\INNER{Av_\ell}{Av_\ell}_0=a_\ell^2 \INNER{v_\ell}{v_\ell}_0=a_\ell^2$
we obtain
\begin{equation*}
   \abs{v_\ell}_1
   =\abs{a_\ell}
   =\begin{cases}a_\ell&\text{, if $\ell>0$,}\\
   -a_\ell&\text{, if $\ell<0$.}\end{cases}
\end{equation*}

\begin{definition}[$\gamma_\pm$]
Let $\gamma_+\colon[0,1]\to\C\setminus(\spec A^{-1})$ be a continuous
loop such that the winding number of $\gamma_+$ with respect to
$1/a_\ell$ is $1$ and with respect to $1/a_{-\ell}$ is $0$ for all $\ell\in \N$.
Since the positive eigenvalues $1/a_\ell$, as well as the negative
eigenvalues $1/a_{-\ell}$, of $A^{-1}$ accumulate at the origin, the
path $\gamma_+$ has to run through the origin.
For definiteness and later use we choose the concatenation
$\gamma_+=\beta\#\alpha_+$ as indicated in
Figure~\ref{fig:fig-concatenation}
(which also shows an analogously defined path $\gamma_-=\bar \beta\#\alpha_-$
where $\bar\beta$ denotes $\beta$ traversed in the opposite
direction).
\end{definition}
\begin{figure}
  \centering
  \includegraphics
                             {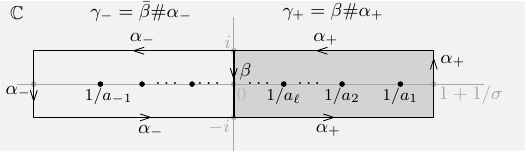}
  \caption{Loops $\gamma_\pm$ enclosing the positive/negative spectrum
                $(1/a_{\pm\nu})_{\nu\in\N}$ of $A^{-1}$}
   \label{fig:fig-concatenation}
\end{figure} 

The following map selects for each $A$
the projection $\Pi_\pm^A\in\Ll(H_0)$ to the positive/negative
eigenspaces of $A$, namely
\begin{equation}\label{eq:Pi+}
   \Pi_\pm^A:=
   \Pi_\pm(A)
   =\frac{i}{2\pi}\int_{\gamma_\pm}
   \underbrace{A\left(\Id-\zeta A\right)^{-1}}_{=:K_\zeta(A)} d\zeta
   =\Pi_\pm(A^{-1}) .
\end{equation}
This is
formula~\cite[(5.22) for the resolvent of $A^{-1}$ instead of $A$]{Kato:1995a}
where we use the identity in the footnote\footnote{
  $R^{A^{-1}}(\zeta)=\left(A^{-1}-\zeta\Id\right)^{-1}
  =\left((\Id-\zeta A) A^{-1}\right)^{-1}
  =A\left(\Id-\zeta A\right)^{-1}$
  }
and the fact that the positive eigenspaces of $A$ and $A^{-1}$ are equal.

\begin{remark}\label{rmk:444}
Let $\Ll_{sym_0}(H_1,H_0)\subset \Ll(H_1,H_0)$ be the subspace
of $H_0$-symmetric elements, and $\Ll_{sym_0}^*\subset \Ll_{sym_0}$
the open subset of isomorphisms.

a)~For $A\in\Ll_{sym_0}^*(H_1,H_0)$ the map $\Pi_\pm$
takes values in $\Ll(H_0)\cap\Ll(H_1)$.
Indeed $\Pi_\pm(A)$
is the projection in $H_0$ to the positive eigenspaces of $A$
and each of them lies in $H_1$, so $\Pi_\pm(A)$ is also a projection in $H_1$.

b)~The intersection $\Ll(H_0)\cap\Ll(H_1)$ is a subset of $\Ll(H_\frac12)$,
by Appendix~\ref{sec:interpolation}.

c)~The intersection $\Ll(H_0)\cap\Ll(H_1)$ is a subset of
$\Ll(H_1,H_0)$, by composing a map $H_1\to H_1$ with the inclusion
$\iota\colon H_1\INTO H_0$.
\\
So the projection selectors $\Pi_\pm$ are maps between the spaces
shown in Figure~\ref{fig:Pi_+}.
\end{remark}

\begin{figure}[h]
  \centering
  \begin{equation*}
  \begin{tikzcd} 
   &&&&
    \Ll(H_{\frac12})
    \arrow[dd, hookrightarrow, "\,\text{\rm$\Ii$~(\ref{eq:Ii})}"]
  \\
  \Pi_\pm\colon
  \Ll_{sym_0}^*(H_1,H_0)
  \arrow[urrrr, dashed]
  \arrow[rrr, "\text{(\ref{eq:Pi+})}", "\text{Rmk.\,\ref{rmk:444}\,a)}"']
  \arrow[drrrr]
    &&&
    \Ll(H_0)\cap\Ll(H_1)
    \arrow[ur, symbol=\subset, "\,\,\text{b)}"']
    \arrow[dr, symbol=\subset, "\,\,\text{c)}"]
  \\
    &&&&
    \Ll(H_1,H_0)
  \end{tikzcd}
  \end{equation*}
  \caption{The spectral projection selectors $\Pi_\pm\colon A\mapsto \Pi_\pm^A$}
   \label{fig:Pi_+}
\end{figure} 

There is a natural embedding $\Ii\colon \Ll(H_\frac12)\INTO\Ll(H_1,H_0)$
via the compact dense inclusions $H_1\subset H_\frac12\subset H_0$.
More precisely, there is the commuting diagram
\begin{equation}\label{eq:Ii}
\begin{tikzcd} [row sep=small] 
H_1
\arrow[rrr, dashed, "\Ii(B)"]
  \arrow[d, hookrightarrow, "\;\iota_1"']
  &&&
  H_0
\\
H_{\frac12}
\arrow[rrr, "B"]
  &&&
  H_{\frac12}
  \arrow[u, hookrightarrow, "\;\iota_0"']
\end{tikzcd} 
\end{equation}
where the operator defined by $\Ii(B):=\iota_0\circ B \circ\iota_1$ is compact.

\boldmath
\subsection{Differentiability of the projection selectors $\Pi_\pm$}
\label{sec:diff-Pi}
\unboldmath

\begin{proposition}\label{prop:dPi+}
The maps $\Pi_\pm$ defined by~(\ref{eq:Pi+}), viewed as maps
$$
   \Pi_\pm\colon\Ll_{sym_0}^*(H_1,H_0)
   \to \Ll(H_1,H_0) ,
$$
are differentiable and the derivative is given by
\begin{equation}\label{eq:dPi+}
\begin{split}
   d\Pi_\pm(A)\colon \Ll_{sym_0}(H_1,H_0)&\to \Ll(H_1,H_0)\\
   \Delta&\mapsto
   \frac{i}{2\pi}\int_{\gamma_\pm}\underbrace{\left(\Id-\zeta A\right)^{-1}\Delta
   \left(\Id-\zeta A\right)^{-1}}_{=: \Dd^A_{1,\zeta}(\Delta)} d\zeta .
\end{split}
\end{equation}
\end{proposition}

That~(\ref {eq:dPi+}) is the derivative of $\Pi_\pm$
can formally be seen as follows.
Write $B=B(A):=\Id-\zeta A$.
Then $dB(A)\Delta=-\zeta\Delta$,
so we get that $dB^{-1}(A)\Delta=-B^{-1}(dB(A)\Delta)B^{-1}
=B^{-1}\zeta\Delta B^{-1}$. 
Clearly $AB=BA$, conjugation with $B^{-1}$ shows that $AB^{-1}=B^{-1}A$.
Put these facts together to obtain
\begin{equation*}
\begin{split}
   d\left(AB^{-1}\right)\Delta
   &=\Delta B^{-1}+A\, dB^{-1}(A)\Delta\\
   &=B^{-1} B \Delta B^{-1}
   +A B^{-1}\zeta\Delta B^{-1}\\
   &=B^{-1} (\Id\underline{-\zeta A}) \Delta B^{-1}
   +B^{-1}\underline{A\zeta}\Delta B^{-1}\\
   &=B^{-1} \Delta B^{-1} .
\end{split}
\end{equation*}

In order to make this rigorous, for $A\in\Ll_{sym_0}^*(H_1,H_0)$ and
$\zeta\in\C$ we abbreviate
$$
   K_\zeta(A):=A(\Id-\zeta A)^{-1}.
$$
For $k\in\N$ and $\Delta\in\Ll_{sym_0}(H_1,H_0)$
we abbreviate further
\begin{equation}\label{eq:Dd_k_A}
   \Dd^A_{k,\zeta}(\Delta)
   :=   \zeta^{k-1} (\Id-\zeta A)^{-1}
   \left(\Delta\left(\Id-\zeta A\right)^{-1}\right)^k .
\end{equation}
The map $\Delta\mapsto \Dd^A_{k,\zeta}(\Delta)$ is $k$-homogenous;
$k=1$ is the integrand in~(\ref{eq:dPi+}), i.e.
\begin{equation}\label{eq:Dd_1_A}
   \Dd^A_{1,\zeta}(\Delta)
   =\left(\Id-\zeta A\right)^{-1}\Delta
   \left(\Id-\zeta A\right)^{-1} .
\end{equation}

\begin{lemma}\label{eq:K_zeta-diff}
Let $A\in\Ll_{sym_0}^*(H_1,H_0)$. If
$\Delta\in\Ll_{sym_0}(H_1,H_0)$ and $\zeta\in\C$ satisfy
\begin{equation}\label{eq:cond-GS}
   \Norm{\zeta\Delta(\Id-\zeta A)^{-1}}_{\Ll(H_0)}
   <1 ,
\end{equation}
then
$$
   K_\zeta(A+\Delta)-K_\zeta(A)
   =\sum_{k=1}^\infty \Dd^A_{k,\zeta}(\Delta) .
$$
\end{lemma}

\begin{proof}
We consider the operator difference
$$
   \Id-\zeta(A+\Delta)
   =\Id-\zeta A-\zeta\Delta
   =\left(\Id-\zeta\Delta(\Id-\zeta A)^{-1}\right)(\Id-\zeta A) .
$$
In view of~(\ref{eq:cond-GS}) this difference is invertible
and its inverse can be constructed with the help of the geometric
series as follows
\begin{equation*}
\begin{split}
   {\color{cyan}\left(\Id-\zeta(A+\Delta)\right)^{-1}}
   &=(\Id-\zeta A)^{-1}\left(\Id-\zeta\Delta(\Id-\zeta A)^{-1}\right)^{-1}
\\
   &=(\Id-\zeta A)^{-1}
   \sum_{k=0}^\infty \left(\zeta\Delta\left(\Id-\zeta A\right)^{-1}\right)^k .
\end{split}
\end{equation*}
In what follows we use the definition of $K_\zeta$ and {\color{gray} add zero}
to obtain equality {\color{gray} 1}, then we use the
{\color{cyan} above formula} to get equality {\color{gray} 2}
\begin{equation*}
\begin{split}
   &K_\zeta(A+\Delta)-K_\zeta(A)
\\
   &\stackrel{{\color{gray}1}}{=} 
   (A+\Delta) {\color{cyan}\, (\Id-\zeta (A+\Delta))^{-1}}
   {\color{gray}-(A+\Delta)(\Id-\zeta A)^{-1}+(A+\Delta)(\Id-\zeta A)^{-1}}
\\
   &\quad
   -A(\Id-\zeta A)^{-1}
\\
   &\stackrel{{\color{gray}2}}{=} 
   (A+\Delta) (\Id-\zeta A)^{-1}
   \sum_{k={\color{red}1}}^\infty \left(\zeta\Delta\left(\Id-\zeta A\right)^{-1}\right)^k
   +\Delta (\Id-\zeta A)^{-1}
\\
   &\stackrel{{\color{gray}3}}{=}
   \sum_{k=2}^\infty \zeta^k (A+\Delta) (\Id-\zeta A)^{-1}
   \left(\Delta\left(\Id-\zeta A\right)^{-1}\right)^k
   \\
   &\quad +(\underline{\underline{A}}+\underline{\Delta}) (\Id-\zeta A)^{-1}
   \zeta\Delta\left(\Id-\zeta A\right)^{-1}
   \\
   &\quad
   +{\color{gray}(\Id-\zeta A)^{-1}(\Id-\zeta A)\;}\Delta (\Id-\zeta A)^{-1}
\\
   &\stackrel{{\color{gray}4}}{=}
   \sum_{k=2}^\infty \zeta^{k-1} 
   (
   \underbrace{\zeta A{\color{gray}\,-\Id}}_{sum_1}
   \underbrace{{\color{gray}+\Id}}_{sum_2}
   +\underbrace{\zeta\Delta }_{sum_3}
   )
   (\Id-\zeta A)^{-1}
   \left(\Delta\left(\Id-\zeta A\right)^{-1}\right)^k
   \\
   &\quad +\zeta \underline{\Delta} (\Id-\zeta A)^{-1} \Delta (\Id-\zeta A)^{-1}
   +(\Id-\zeta A)^{-1}{\color{red}\zeta \underline{\underline{A}} \Delta\,} (\Id-\zeta A)^{-1}
   \\
   &\quad +(\Id-\zeta A)^{-1} ({\color{cyan}\,\Delta} {\color{red}\,-\zeta A \Delta}) (\Id-\zeta A)^{-1}
\\
\end{split}
\end{equation*}
\begin{equation*}
\begin{split}
   &\stackrel{{\color{gray}5}}{=}
   {\color{brown}
   \zeta \left(\Delta (\Id-\zeta A)^{-1}\right)^2
   }
   +(\Id-\zeta A)^{-1}{\color{cyan}\,\Delta} (\Id-\zeta A)^{-1}
   \\
   &\quad
   {\color{brown}
   -\sum_{k=2}^\infty \zeta^{k-1} \left(\Delta\left(\Id-\zeta A\right)^{-1}\right)^k
   }
   {\color{gray}
      \left(=-\sum_{k=1}^\infty \zeta^k 
      \left(\Delta\left(\Id-\zeta A\right)^{-1}\right)^{k+1}\right)
   }
   \\
   &\quad+\sum_{k=2}^\infty \zeta^{k-1} (\Id-\zeta A)^{-1}
   \left(\Delta\left(\Id-\zeta A\right)^{-1}\right)^k
   \\
   &\quad
   {\color{brown}
   +\sum_{k=2}^\infty \zeta^k
   \left(\Delta\left(\Id-\zeta A\right)^{-1}\right)^{k+1}
   }
\\
   &\stackrel{{\color{gray}6}}{=}
   {\color{brown}
   0+\,
   }
   \sum_{k={\color{cyan}1}}^\infty
   \underbrace{
   \zeta^{k-1} (\Id-\zeta A)^{-1}
   \left(\Delta\left(\Id-\zeta A\right)^{-1}\right)^k
   }_{=:\Dd^A_{k,\zeta}(\Delta)} .
\end{split}
\end{equation*}
In equality {\color{gray} 3} we write the summand for $k=1$
separately and multiply by the {\color{gray} identity}.
Equality {\color{gray} 4} uses that $A$ commutes with $(\Id-\zeta A)^{-1}$.
Equality {\color{gray} 5} uses that the two terms involving
${\color{red} \zeta A \Delta}$ cancel, moreover, we split the sum into three sums.
In equality {\color{gray} 6} {\color{brown} three terms} cancel.
\end{proof}

\begin{proof}[Proof of Proposition~\ref{prop:dPi+}]
The proof for $\Pi_+$ is in four steps.

\smallskip\noindent
\textsc{Step~1.}
Condition~(\ref{eq:cond-GS}) holds for $\zeta$ along the image of
the path $\gamma_+=\beta\#\alpha_+$;
see Figure~\ref{fig:fig-concatenation}.

\smallskip\noindent
\textbf{Along \boldmath$\alpha_+$.}
Since the path $\alpha_+$ is uniformly bounded away from
$\spec A^{-1}$, there is a constant
$c_1>0$ with $\norm{(\Id-\zeta A)^{-1}}_{\Ll(H_0,H_1)}\le c_1$
for any $\zeta$ on the image of the path $\alpha_+$.
Because the image of $\alpha_+$ is compact, there is a constant $c_2>0$
such that $\abs{\zeta}\le c_2$ for any $\zeta \in\Im\,\alpha_+$.
Therefore, if $\norm{\Delta}_{\Ll(H_1,H_0)}< 1/c_1c_2$, then
condition~(\ref{eq:cond-GS}) is satisfied for any $\zeta \in\Im\,\alpha_+$.

\smallskip\noindent
\textbf{Along \boldmath$\beta$.}
In this case $\zeta=it$ for $t\in[-1,1]$.
Then there exists a constant $c_3>0$ such that
$\norm{(\Id-it A)^{-1}}_{\Ll(H_0,H_1)}\le c_3/\abs{t}$.
If $\norm{\Delta}_{\Ll(H_1,H_0)}<1/c_3$, then
$$
   \norm{it\Delta(\Id-it A)^{-1}}_{\Ll(H_0)}
   \le\abs{t} \norm{\Delta}_{\Ll(H_1,H_0)}\norm{(\Id-itA)^{-1}}_{\Ll(H_0,H_1)}
   <1.
$$

So, from now on we suppose that
$\norm{\Delta}_{\Ll(H_1,H_0)}<\min\{\frac{1}{c_1c_2},\frac{1}{c_3}\}$.
By the above this guarantees that~(\ref{eq:cond-GS}) is satisfied for any $\zeta$ in the
image of the path~${\gamma_+}$.

By definition~(\ref{eq:Pi+}) of $\Pi_\pm(A+\Delta)$ and $\Pi_\pm(A)$,
together with formula~(\ref{eq:dPi+}) for the {\color{red}candidate} for the
derivative $d\Pi_\pm(A) \Delta$, we obtain the first equality in
\begin{equation}\label{eq:step1-f}
\begin{split}
   &\Pi_\pm(A+\Delta)-\Pi_\pm(A)
   -{\color{red}\tfrac{i}{2\pi}\int_{\gamma_\pm} \Dd^A_{1,\zeta}(\Delta) \, d\zeta}\\
   &=\frac{i}{2\pi}\int_{\gamma_\pm} K_\zeta(A+\Delta)-K_\zeta(A)
   -\Dd^A_{1,\zeta}(\Delta) \, d\zeta\\
   &\stackrel{{\color{gray}2}}{=} 
   \frac{i}{2\pi} \int_{\gamma_\pm} \sum_{k=2}^\infty
   \Dd^A_{k,\zeta}(\Delta) \, d\zeta
\end{split}
\end{equation}
where equality {\color{gray} 2} is by Lemma~\ref{eq:K_zeta-diff}.

\smallskip\noindent
\textsc{Step~2.}
We investigate the norm of the terms $\Dd^A_{k,\zeta}(\Delta)$
for $k\ge 2$ and for $\zeta$ along the path $\gamma_+=\beta\#\alpha_+$.

\smallskip\noindent
\textbf{Along \boldmath$\alpha_+$.}
By the continuous inclusion $H_1\INTO H_0$ and
by definition~(\ref{eq:Dd_k_A}) of $\Dd^A_{k,\zeta}$ as a composition
and $c_1,c_2$ in Step~1 we get, respectively, inequalities one and two
$$
   \norm{\Dd^A_{k,\zeta}(\Delta)}_{\Ll(H_1,H_0)}
   \le \norm{\Dd^A_{k,\zeta}(\Delta)}_{\Ll(H_0,H_1)}
   \le c_2^{k-1} c_1^{k+1}
   \norm{\Delta}_{\Ll(H_1,H_0)}^k.
$$
%
\textbf{Along \boldmath$\beta$.}
There is a constant $c_4>0$ such that
$$
   \norm{(\Id - it A)^{-1}}_{\Ll(H_0)}\le c_4
   ,\qquad
   \norm{(\Id - it A)^{-1}}_{\Ll(H_1)}\le c_4 ,
$$
whenever $t\in[-1,1]$. Therefore, for $\zeta=it$ and using the fact
that~(\ref{eq:Dd_k_A}) is a composition of operators
$(\Id-it A)^{-1}\Delta\colon H_1\to H_0\to H_1$, we get
$$
   \norm{\Dd^A_{k,\zeta}(\Delta)}_{\Ll(H_1,H_0)}
   \le \abs{t}^{k-1} c_4^2 \bigl(\tfrac{c_3}{\abs{t}}\bigr)^{k-1}
   \norm{\Delta}_{\Ll(H_1,H_0)}^k
   = c_4^2 c_3^{k-1} \norm{\Delta}_{\Ll(H_1,H_0)}^k .
$$

\smallskip\noindent
\textsc{Step~3.}
Fix $\eps\in(0,1)$.
Let $C:=\max\{1,c_1c_2,c_1^2c_2,c_3,c_3c_4\}$.
If $\norm{\Delta}_{\Ll(H_1,H_0)}=\frac{\eps}{2C}$, then
\begin{equation}\label{eq:step3-ff}
   \norm{\Dd^A_{k,\zeta}(\Delta)}_{\Ll(H_1,H_0)}
   \le \tfrac{\eps^k}{2^k}
\end{equation}
for any $\zeta$ along the path ${\gamma_+}$. This is true by Step~2.

\smallskip\noindent
\textsc{Step~4.} We prove the proposition.
By~(\ref{eq:step1-f}) we get equality {\color{gray} 1} in what follows
\begin{equation*}
\begin{split}
   &\Bigl\|\Pi_\pm(A+\Delta)-\Pi_\pm(A)
   -{\color{red}\tfrac{i}{2\pi}\int_{\gamma_\pm} \Dd^A_{1,\zeta}(\Delta) \, d\zeta}
   \Bigr\|_{\Ll(H_1,H_0)}\\
   &=\Bigl\|\tfrac{1}{2\pi} \sum_{k=2}^\infty \int_{\gamma_\pm}
   \Dd^A_{k,\zeta}(\Delta)\, d\zeta\Bigr\|_{\Ll(H_1,H_0)}\\
   &\stackrel{{\color{gray}1}}{\le} \tfrac{1}{2\pi} \sum_{k=2}^\infty \int_{\gamma_\pm}
   \norm{\Dd^A_{k,\zeta}(\Delta)}_{\Ll(H_1,H_0)}\, d\abs{\zeta}\\
   &\stackrel{{\color{gray}2}}{\le} \tfrac{1}{2\pi} \sum_{k=2}^\infty \int_{\gamma_\pm}
   \bigl(\tfrac{\eps}{2}\bigr)^k \, d\abs{\zeta}\\
   &= \tfrac{\mathrm{length}({\gamma_+})}{2\pi} \sum_{k=2}^\infty 
   \bigl(\tfrac{\eps}{2}\bigr)^k \\
   &= \tfrac{\mathrm{length}({\gamma_+})}{2\pi}
   \bigl(\tfrac{\eps}{2}\bigr)^2
   \frac{1}{1-\tfrac{\eps}{2}}\\
   &\le \tfrac{\mathrm{length}({\gamma_+})}{4\pi} \eps^2\\
   &=\tfrac{\mathrm{length}({\gamma_+})}{\pi} C^2
   \norm{\Delta}_{\Ll(H_1,H_0)}^2 .
\end{split}
\end{equation*}
Here inequality {\color{gray} 2} is by estimate (\ref{eq:step3-ff}).
The final equality holds by Step~3.

\smallskip
This proves Proposition~\ref{prop:dPi+} for $\Pi_+$.
The proof for $\Pi_-$ is analogous.
\end{proof}

\boldmath
\subsection{Directional derivative of $\Pi_\pm$ extends from $H_1$ to $H_{\frac12}$}
\label{sec:directional}
\unboldmath

By Proposition~\ref{prop:dPi+} the derivative of
$\Pi_\pm$ in~(\ref{eq:Pi+}) and Figure~\ref{fig:Pi_+}, viewed as map
$$
   \Pi_\pm \colon\Ll_{sym_0}^*(H_1,H_0)
   \to \Ll(H_1,H_0),
$$
exists at $A\in\Ll_{sym_0}^*$ in direction $\Delta\in \Ll_{sym_0}(H_1,H_0)$
and is a bounded linear map $S_\pm :=d\Pi_\pm (A)\Delta\colon H_1\to H_0$.
The theorem below shows, in particular, that $S_\pm $ extends to a bounded linear map
$\tilde S_\pm $ on $H_{\frac12}$, more precisely the diagram
\begin{equation*}
\begin{tikzcd} [row sep=small] 
H_{\frac12}
\arrow[rrr, dashed, "\tilde S_\pm "]
  &&&
  H_{\frac12}
  \arrow[d, hookrightarrow, "\;\iota_0"]
\\
H_1
  \arrow[u, hookrightarrow, "\;\iota_1"]
\arrow[rrr, "{S_\pm =d\Pi_\pm (A)\Delta}"]
  &&&
  H_0
\end{tikzcd} 
\end{equation*}
commutes, in symbols $S_\pm =\iota_0\circ \tilde S_\pm  \circ\iota_1$.

\begin{theorem}\label{thm:dPi+-ext}
Let $A\in \Ll_{sym_0}^*(H_1,H_0)$.
Then the derivative $d\Pi_\pm (A)$ 
from Proposition~\ref{prop:dPi+} factorizes through $\Ll(H_\frac12)$ 
via a map $T_A^\pm$, namely
\begin{equation*}
\begin{tikzcd} 
 &&&&
  \Ll(H_{\frac12})
  \arrow[d, hookrightarrow, "\text{\rm$\Ii$ cf.~(\ref{eq:Ii})}"]
\\
\Ll_{sym_0}(H_1,H_0)
\arrow[rrrr, "d\Pi_\pm (A)"]
\arrow[urrrr, dashed, "{T_A^\pm}"]
  &&&&
  \Ll(H_1,H_0)
\end{tikzcd}
\end{equation*}
still given by formula~(\ref{eq:dPi+}), namely
\begin{equation}\label{eq:dPi+-ext}
\begin{split}
   T_A^\pm\Delta
   =d\Pi_\pm(A)\Delta
   :=\frac{i}{2\pi}\int_{\gamma_\pm}\left(\Id-\zeta A\right)^{-1}\Delta
   \left(\Id-\zeta A\right)^{-1} d\zeta .
\end{split}
\end{equation}
Moreover, there is the bound
$$
   \bigl\|T_A^\pm\bigr\|_{\Ll(\Ll_{sym_0}(H_1,H_0), \Ll(H_{\frac12}))}
   \le \pi+\frac{4\sigma+2}{\pi\sigma^2} .
$$
\end{theorem}

As preparation for the proof of Theorem~\ref{thm:dPi+-ext}
we use additivity of the integral along the domain of integration
$\gamma_\pm=\beta\#\alpha_\pm$ (Figure~\ref{fig:fig-concatenation})
to decompose $d\Pi_\pm(A)$ into the sum of operators
\begin{equation}\label{eq:T_A^beta+T_A^alpha}
   d\Pi_+(A)=T_A^\beta+T_A^{\alpha_+}
   ,\qquad
   d\Pi_-(A)=T_A^{\bar \beta}+T_A^{\alpha_-}=-T_A^\beta+T_A^{\alpha_-} ,
\end{equation}
where the operators
$$
   T_A^*\colon \Ll_{sym_0}(H_1,H_0)\to\Ll(H_1,H_0)
$$
are for $*\in\{\beta,\alpha_+,\alpha_-\}$ defined by
\begin{equation}\label{eq:T_A^beta}
\begin{split}
   T_A^*\Delta
   :=\frac{i}{2\pi}\int_*\left(\Id-\zeta A\right)^{-1}\Delta
   \left(\Id-\zeta A\right)^{-1} d\zeta
   .
\end{split}
\end{equation}
We show that all these three operators factorize through $\Ll(H_\frac12)$
and by abuse of notation we keep on using the same notation.

\boldmath
\subsubsection{Block decomposing $d\Pi_\pm^A$ via positive/negative subspaces}
\unboldmath

Suppose that $A\in \Ll_{sym_0}^*(H_1,H_0)$ and pick
an ONB $\Vv=\Vv(A)$ of $H_0$ as in~(\ref{eq:ONB-Vv-H_0})
consisting of eigenvectors $v_\ell\in H_1$ of $A$ where $\ell\in\Z^*$ associated
to eigenvalues $a_\nu>0$ for and $a_{-\nu}<0$ for $\nu\in\N$.
For $\ell\in\Z^*$ let
$$
   E_\ell=\ker\left(a_\ell\,\Id-A\right)\subset H_1\subset H_0
$$
be the eigenspace for the eigenvalue $a_\ell$. For $r\in\R$ the
\textbf{positive/negative subspaces} of $H_r$, see~(\ref{eq:ell^2_f^r}),
are defined by taking the closure in the $H_r$ norm as follows
$$
   H_r^+:=\overline{\bigoplus_{\nu\in\N} E_\nu}^{H_r}
   ,\qquad
   H_r^-:=\overline {\bigoplus_{\nu\in\N} E_{-\nu}}^{H_r} .
$$
Thus there are the $H_r$ orthogonal decompositions
$$
   H_r=H_r^+\oplus H_r^-
   ,\quad r\in\R.
$$
The maps defined by~(\ref{eq:Pi+}) are projections
$\Pi_\pm^A=(\Pi_\pm^A)^2$ in $H_i$ for $i=0,1$, same notation.
These maps satisfy $\Pi_+^A+\Pi_-^A=\Id$ in $H_i$ and
$$
   H_i^+=\ker \Pi_-^A{\color{gray}\,=\Im\,\Pi_+^A=\Fix\,\Pi_+^A,}
   \qquad
   H_i^-=\ker \Pi_+^A{\color{gray}\,=\Im\,\Pi_-^A=\Fix\,\Pi_-^A.}
$$
For $\flat,\sharp\in\{-,+\}$ we abbreviate
\begin{equation*}
\begin{split}
   (d\Pi_\pm^A)^{\flat\sharp}\colon \Ll_{sym_0}(H_1,H_0)
   &\to\Ll_{sym_0}(H_1^\flat,H_0^\sharp)
\\
   \Delta&\mapsto \Pi_\sharp^A\circ d\Pi_\pm(A)\Delta \circ\Pi_\flat^A .
\end{split}
\end{equation*}

\begin{lemma}\label{le:p-dp-p}
For $A\in \Ll_{sym_0}^*(H_1,H_0)$ and $\Delta\in
\Ll_{sym_0}(H_1,H_0)$ the composition
\begin{equation*}
\begin{tikzcd} 
0\colon
H_1
\arrow[r, "\Pi_\pm^A"]
  &
  H_1
  \arrow[rr, "d\Pi_\pm(A)\Delta"]
    &&
    H_0
    \arrow[r, "\Pi_\pm^A"]
      &H_0
\end{tikzcd}
\end{equation*}
vanishes, in symbols $\Pi_\pm^A\circ d\Pi_\pm(A)\Delta\circ\Pi_\pm^A=0$.
\end{lemma}

\begin{proof}
Since $p=p(A):=\Pi_\pm^A\in \Ll(H_0)\cap\Ll(H_1)$ is a projection
the identity $p^2=p$ holds.
Differentiating this identity we get
$
   p\circ dp+dp\circ p=dp
$. Multiplying this identity from the right by $p$ yields
$p\circ dp\circ p+dp\circ p^2=dp\circ p$ which, since $p^2=p$, is equivalent to
$p\circ dp\circ p=0$.
\end{proof}

\begin{corollary}\label{cor:bjkhbjg6544}
For any $A\in \Ll_{sym_0}^*(H_1,H_0)$ the diagonal blocks vanish
$$
   (d\Pi_\pm^A)^{++}=0
   ,\qquad
   (d\Pi_\pm^A)^{--}=0 .
$$
\end{corollary}

\begin{proof}
Since $\Pi_+^A+\Pi_-^A=\Id$
it follows that $d\Pi_+(A)=-d\Pi_-(A)$.
Hence the corollary follows from Lemma~\ref{le:p-dp-p}.
\end{proof}

For $A\in \Ll_{sym_0}^*(H_1,H_0)$
and subpaths $*\in\{\beta,\alpha_+,\alpha_-\}$
and block indices $\flat,\sharp\in\{-,+\}$ we abbreviate
\begin{equation}\label{eq:4block-2}
\begin{split}
   (T_A^*)^{\flat\sharp}\colon \Ll_{sym_0}(H_1,H_0)
   &\to\Ll_{sym_0}(H_1^\flat,H_0^\sharp)
\\
   \Delta&\mapsto \Pi_\sharp^A\circ T_A^*\Delta \circ\Pi_\flat^A .
\end{split}
\end{equation}

\begin{corollary}
For any $A\in \Ll_{sym_0}^*(H_1,H_0)$ there are the block identities
$$
   (T_A^\beta)^{++}
   =-(T_A^{\alpha_+})^{++}
   ,\qquad
   (T_A^\beta)^{--}
   =-(T_A^{\alpha_+})^{--} .
$$
\end{corollary}

\begin{proof}
Corollary~\ref{cor:bjkhbjg6544} for $d\Pi_+^A$ and
the first identity in~(\ref{eq:T_A^beta+T_A^alpha}).
\end{proof}

\begin{lemma}\label{le:4block}
Let $H$ be any Hilbert space with an orthogonal decomposition
$H=H^+\oplus H^-$. Any operator $B\in\Ll(H)$ block decomposes into
four operators
\begin{equation}\label{eq:4block}
   B
   =
   \begin{pmatrix}
   B^{++}&B^{-+}\\
   B^{+-}&B^{--}
   \end{pmatrix}
   \colon H^+\oplus H^- \to H^+\oplus H^-
\end{equation}
and the operator norms are related by
$$
   \norm{B}
   \le
   \norm{B^{++}}+\norm{B^{-+}}+\norm{B^{+-}}+\norm{B^{--}}
   .
$$
\end{lemma}

\begin{proof}
Due to the decomposition $H=H^+\oplus H^-$
any $v\in H$ is a sum $v=v_++v_-$
for unique elements $v_\pm\in H^\pm$.
By orthogonality identity one (Pythagoras) holds
\begin{equation*}
\begin{split}
   \abs{Bv}^2
   &=\abs{(Bv)_+}^2+\abs{(Bv)_-}^2
\\
   &\le \bigl(\underbrace{\abs{B^{++} v_+}+\abs{B^{-+} v_-}}_{a}\bigr)^2
   +\bigl(\underbrace{\abs{B^{+-} v_+}+\abs{B^{--} v_-}}_{b}\bigr)^2
   {\color{gray}\small \le (a+b)^2}
\\
   &\le \bigl(
   \abs{B^{++} v_+}+\abs{B^{-+} v_-}
   +\abs{B^{+-} v_+}+\abs{B^{--} v_-}
   \bigr)^2.
\end{split}
\end{equation*}
Here 
inequality one holds since $(Bv)_+=B^{++} v_++B^{-+} v_-$,
by~(\ref{eq:4block}), together with the triangle inequality.
By definition of the operator norm and the previous estimate
we get
\begin{equation*}
\begin{split}
   \norm{B}:
   &=\sup_{\abs{v}=1}\abs{Bv}
\\
   &\le\sup_{1=\abs{v}\atop =\abs{v_+}+\abs{v_-}}
   \bigl(
   \abs{B^{++} v_+}
   +\abs{B^{-+} v_-}
   +\abs{B^{+-} v_+}
   +\abs{B^{--} v_-}
   \bigr)
\\
   &\le \norm{B^{++}}+\norm{B^{-+}}+\norm{B^{+-}}+\norm{B^{--}} .
\end{split}
\end{equation*}
where we used that $\abs{B^{++} v_+}\le \norm{B^{++}}\cdot 1$ and so on.
This proves Lemma~\ref{le:4block}.
\end{proof}

\boldmath
\subsection*{The operator $T_A^{\alpha_{\pm}}$}
\unboldmath

\boldmath
\subsubsection{Estimating the operator $T_A^{\alpha_{\pm}}$}
\unboldmath

\begin{proposition}\label{prop:T_A^alpha}
$
   \norm{T_A^{\alpha_{\pm}}}_{\Ll(\Ll_{sym_0}(H_1,H_0), \Ll(H_{\frac12}))}
   \le \frac{2\sigma+1}{\pi\sigma^2}
$.
\end{proposition}

\begin{lemma}
For any $\zeta\in \im \alpha_\pm$ 
it holds $\norm{(\Id-\zeta A)^{-1}}_{\Ll(H_0,H_1)} \le 1$.
\end{lemma}

\begin{proof}
The proof is in two steps.

\smallskip\noindent
\textsc{Step~1.}
It holds $\norm{(\Id-\zeta A)^{-1}}_{\Ll(H_0,H_1)}
=\max_{\ell\in\Z^*}\frac{\abs{a_\ell}}{\abs{1-\zeta a_\ell}}$.

\smallskip\noindent
\textsc{Proof.}
As $\abs{(\Id-\zeta A)^{-1}x}_{H_1}=\abs{A(\Id-\zeta A)^{-1}x}_{H_0}$,
by~(\ref{eq:A-norm}), we get equality {\color{gray}2}
\begin{equation*}
\begin{split}
   \norm{(\Id-\zeta A)^{-1}}_{\Ll(H_0,H_1)}
   &=\sup_{\norm{x}_0=1} \abs{(\Id-\zeta A)^{-1}x}_{H_1}\\
   &\stackrel{{\color{gray}2}}{=}
   \sup_{\norm{x}_0=1} \abs{A(\Id-\zeta A)^{-1}x}_{H_0}\\
   &=\norm{A(\Id-\zeta A)^{-1}}_{\Ll(H_0)}\\
   &=\max_{\ell\in\Z^*}\frac{\abs{a_\ell}}{\abs{1-\zeta a_\ell}}
\end{split}
\end{equation*}
where the last equality uses that the orthonormal basis $\Vv$ of $H_0$,
see~(\ref{eq:ONB-Vv-H_0}), consists of eigenvectors $v_\ell$ such that
operator $A$ with eigenvalues $a_\ell$ and therefore
$A(\Id-\zeta A)^{-1}v_\ell=a_\ell(1-\zeta a_\ell)^{-1} v_\ell$.

\smallskip\noindent
\textsc{Step~2.}
For all $\ell\in\Z^*$ and $\zeta\in\im \alpha_\pm$
it holds $\abs{a_\ell/(1-\zeta a_\ell)}\le 1$.

\smallskip\noindent
\textsc{Proof.}
We distinguish four cases (the four edges of the rectangle in
Figure~\ref{fig:fig-concatenation}).

\smallskip\noindent
\textbf{Case 1.} Let $\zeta=i+t$ where $t\in[-1-1/\sigma,1+1/\sigma]$.
We estimate for the denominator
$
   \abs{1-\zeta a_\ell}
   =\abs{1-t a_\ell-i a_\ell}\ge \abs{a_\ell}
$. Hence $\abs{a_\ell/(1-\zeta a_\ell)}\le 1$.

\smallskip\noindent
\textbf{Case 2.} Let $\zeta=1+1/\sigma+it$ where $t\in[-1-1]$.
We estimate for the denominator
$
   \abs{1-\zeta a_\ell}
   =\abs{1-a_\ell-a_\ell/\sigma-it a_\ell}
   \ge\abs{1-a_\ell-a_\ell/\sigma}
$.
\\
Case $\ell>0$. Then $a_\ell\ge a_1>0$, hence $a_\ell/\sigma\ge a_\ell/a_1\ge 1$.
Thus $\abs{1-a_\ell-a_\ell/\sigma}=\abs{a_\ell+a_\ell/\sigma-1} \ge \abs{a_\ell}$.
Hence $\abs{a_\ell/(1-\zeta a_\ell)}\le 1$.
\\
Case $\ell<0$. Then $1-a_\ell-a_\ell/\sigma
=1+\abs{a_\ell}+\abs{a_\ell}/\sigma \ge \abs{a_\ell}$.
Hence $\abs{a_\ell/(1-\zeta a_\ell)}\le 1$.

\smallskip\noindent
\textbf{Case 3.} Let $\zeta=-1-1/\sigma+it$ where $t\in[-1-1]$.
Analogue to Case~2.

\smallskip\noindent
\textbf{Case 4.} Let $\zeta=-i+t$ where $t\in[-1-1/\sigma,1+1/\sigma]$.
Analogue to Case~1.
\end{proof}

\begin{corollary}\label{cor:hgh5465vgv}
For all $\Delta\in \Ll_{sym_0}(H_1,H_0)$ and $\zeta\in\im \alpha_\pm$ it holds
\begin{equation}\label{eq:hgh5465vgv}
   \norm{(\Id-\zeta A)^{-1}\Delta (\Id-\zeta A)^{-1}}_{\Ll(H_i)}
   \le \tfrac{1}{\sigma}\norm{\Delta}_{\Ll(H_1,H_0)}
\end{equation}
whenever $i\in\{0,\frac12,1\}$.
\end{corollary}

\begin{proof}
The operator norm of the inclusion $\iota\colon H_1\INTO H_0$ is
bounded by
\begin{equation*}
\begin{split}
   \norm{\iota}_{\Ll(H_1,H_0)}
   =\sup_{\Abs{x}_1=1} \Abs{x}_0
   =\tfrac{1}{\sigma}\sup_{\Abs{x}_1=1} \sigma\Abs{x}_0
   \le\tfrac{1}{\sigma}
\end{split}
\end{equation*}
where we used that $\sigma \Abs{x}_0\le \Abs{Ax}_0=\Abs{x}_1=1$.
In case $i=0$ we get
\begin{equation*}
\begin{split}
   &\norm{\iota(\Id-\zeta A)^{-1}\Delta (\Id-\zeta A)^{-1}}_{\Ll(H_0)}\\
   &\le \norm{\iota}_{\Ll(H_1,H_0)}
   \norm{(\Id-\zeta A)^{-1}}_{\Ll(H_0,H_1)}
   \norm{\Delta}_{\Ll(H_1,H_0)}
   \norm{(\Id-\zeta A)^{-1}}_{\Ll(H_0,H_1)}\\
   &\le\tfrac{1}{\sigma} \norm{\Delta}_{\Ll(H_1,H_0)}
\end{split}
\end{equation*}
and in case $i=1$ we get
\begin{equation*}
\begin{split}
   &\norm{(\Id-\zeta A)^{-1}\Delta (\Id-\zeta A)^{-1}\iota}_{\Ll(H_1)}\\
   &\le
   \norm{(\Id-\zeta A)^{-1}}_{\Ll(H_0,H_1)}
   \norm{\Delta}_{\Ll(H_1,H_0)}
   \norm{(\Id-\zeta A)^{-1}}_{\Ll(H_0,H_1)}
   \norm{\iota}_{\Ll(H_1,H_0)}\\
   &\le\tfrac{1}{\sigma} \norm{\Delta}_{\Ll(H_1,H_0)} .
\end{split}
\end{equation*}
The case $i=\frac12$ follows from $i=0$ and $i=1$ by
interpolation (Prop.~\ref{prop:Stein}).
\end{proof}

\begin{proof}[Proof of Proposition~\ref{prop:T_A^alpha}]
Pick $\Delta\in \Ll_{sym_0}(H_1,H_0)$ and estimate
\begin{equation*}
\begin{split}
   \norm{T_A^{\alpha_{\pm}}\Delta}_{\Ll(H_{\frac12})}
   &\stackrel{\text{(\ref{eq:T_A^beta})}}{\le}
   \frac{1}{2\pi} \int_{\alpha_\pm}
   \norm{(\Id-\zeta A)^{-1}\Delta (\Id-\zeta A)^{-1}}_{\Ll(H_{\frac12})} d\abs{\zeta}\\
   &\stackrel{\text{(\ref{eq:hgh5465vgv})}}{\le}
   \frac{1}{2\pi\sigma}\cdot\norm{\Delta}_{\Ll(H_1,H_0)} 
   \cdot\underbrace{\mathrm{length}(\alpha_\pm)}_{2+2(1+1/\sigma)}
\end{split}
\end{equation*}
where the length calculation is illustrated by Figure~\ref{fig:fig-concatenation}.
\end{proof}

\boldmath
\subsection*{The operator $T_A^\beta$}
\unboldmath

\boldmath
\subsubsection{Matrix representation and Hadamard product}
\unboldmath

If $\Delta\in \Ll_{sym_0}(H_1,H_0)$ we identify
$\Delta$ with its matrix representation
$$
   \mat{d}:=[\Delta]_\Vv\in\R^{\Z^*\times\Z^*}=\Map(\Z^*\times\Z^*,\R)
$$
whose entries are given by
\begin{equation}\label{eq:d_ell_m}
   d_{n \ell}:=\INNER{v_n}{\Delta v_\ell}_0
\end{equation}
where $\Vv=(v_\ell)_{\ell\in\Z^*}$ is the ONB of $H_0$ composed of
eigenvectors of $A$; see~(\ref{eq:ONB-Vv-H_0}).
Hence
\begin{equation}\label{eq:Delta-v_m}
   \Delta v_\ell=\sum_{n\in\Z^*} d_{n\ell} v_n .
\end{equation}

\begin{definition}
Let $\mat{b}$ be the (infinite) matrix whose entries are defined by
\begin{equation}\label{eq:b-1}
   b_{\ell m}
   :=\begin{cases}
   \frac{\arctan a_m-\arctan a_\ell}{a_m-a_\ell}&\text{, if $\ell\not=m$,}\\
   \frac{1}{a_\ell^2+1}&\text{, if $\ell=m$ .}
   \end{cases}
\end{equation}
Then the \textbf{Hadamard}, i.e. entry-wise, \textbf{product}
of the matrizes $\mat{b}$ and $\mat{d}$ is the matrix $\mat{c}$
defined by
$$
   \mat{c}=\mat{b} \odot \mat{d}
   ,\qquad
   \text{$c_{\ell m}:=b_{\ell m} d_{\ell m}$, for $\ell, m\in\Z^*$.}
$$
\end{definition}

\boldmath
\subsubsection{Matrix block decomposition}
\unboldmath

Since the spectrum of $A$ has a positive and a negative part,
and since we use the convention $\ell^2_{\Z^*}=\ell^2_\N\oplus\ell^2_{-\N}$,
our matrices decompose into four blocks
\begin{equation*}
\begin{split}
   \mat{d}=
   (d_{mn})_{m,n\in\Z^*}
   &=
   \begin{pmatrix}
   (d_{\mu\nu})&(d_{\mu(-\nu)})\\
   (d_{(-\mu)\nu})&(d_{(-\mu)(-\nu)})
   \end{pmatrix}
\\
   &=
   \begin{pmatrix}
   \;\;(d_{\mu\nu}^{++})&\;\quad(d_{\mu\nu}^{-+})\quad\mbox{ }\\
   \;\;(d_{\mu\nu}^{+-})&\;\quad(d_{\mu\nu}^{--})\quad\mbox{ }
   \end{pmatrix}
   \in\Ll(\ell^2_{+\N}\oplus\ell^2_{-\N})
\end{split}
\end{equation*}
where $d^{++}=(d_{\mu\nu}^{++})$ abbreviates
$(d_{\mu\nu}^{++})_{\mu,\nu\in\N}$ and so on.
Analogously for
$$
   \mat{b}=
   (b_{mn})_{m,n\in\Z^*}
   =
   \begin{pmatrix}
   \mat{b}^{++}&\mat{b}^{-+}\\
   \mat{b}^{+-}&\mat{b}^{--}
   \end{pmatrix}
   \colon \ell^2_{+\N}\oplus\ell^2_{-\N} \to
   \ell^2_{+\N}\oplus\ell^2_{-\N} .
$$
Being entry-wise, the \textbf{Hadamard multiplication operator} defined by
\begin{equation*}
\begin{split}
   H_{\mat{b}}\colon \Ll_{sym_0}(H_1,H_0)\to \Ll_{sym_0}(H_0),\quad
   \mat{d}\mapsto \mat{b}\odot\mat{d} =: \mat{c}
\end{split}
\end{equation*}
respects the four block decomposition $++$, $+-$, $-+$, $--$.

\boldmath
\subsubsection{Representing $T_A^\beta$ as Hadamard multiplication
  operator $H_{\mat{b}}$}
\unboldmath

\begin{proposition}
Applying the operator $T_A^\beta$ to $\Delta\in \Ll_{sym_0}(H_1,H_0)$,
in symbols $\Delta\mapsto T_A^\beta\Delta$,
corresponds to Hadamard multiplying $\mat{d}$ by $\mat{b}$,
in symbols $\mat{d}\mapsto \mat{b}\odot\mat{d}$.
\end{proposition}

\begin{proof}
Pick elements $v_\ell ,v_m$ of the orthonormal basis 
$\Vv=\{v_\ell\}_{\ell\in\Z^*}\subset H_1$ of $H_0$.
Recall that $Av_\ell=a_\ell v_\ell$ for $\ell\in\Z^*$ where $a_\ell\in\R\setminus\{0\}$.
%
%
For $\ell\in\Z^*$ we have
\begin{equation*}
   (\Id-\zeta A)^{-1}v_\ell=\frac{1}{1-\zeta a_\ell} v_\ell .
\end{equation*}
As in~(\ref{eq:d_ell_m}) we denote the matrix elements of $\Delta$ by 
$d_{n\ell}:=\INNER{v_n}{\Delta v_\ell}_0$ for $\ell,m\in\Z^*$.
Hence $\Delta v_\ell=\sum_{n\in\Z^*} d_{n\ell} v_n$.
By definition of $\Dd^a_{1,\zeta}(\Delta)$, by symmetry of $A$, and
since the inner product now being complex (say complex anti-linear in
the second variable), we get equality one
\begin{equation*}
\begin{split}
   c_{\ell m}
   :&=\INNER{
   \left(\frac{i}{2\pi}\int_{\gamma_1}\Dd^A_{1,\zeta}(\Delta)\, d\zeta\right) v_\ell
   }{v_m}_0
\\
   &\stackrel{\text{{\color{gray} 1}}}{=}
   \frac{i}{2\pi}\int_{\gamma_1}
   \INNER{\Delta\left(\Id-\zeta A\right)^{-1} v_\ell}
   {\left(\Id-\bar \zeta A\right)^{-1} v_m}_0 d\zeta
\\
   &=\frac{i \inner{\Delta v_\ell}{v_m}_0}{2\pi}\int_{\gamma_1}
   \frac{1}{(1-\zeta a_\ell)(1-\zeta a_m)}
   d\zeta
\\
   &\stackrel{\text{{\color{gray} 3}}}{=}
   \frac{i d_{\ell m}}{2\pi}\int_{-1}^1
   \frac{1}{(1+it a_\ell)(1+it a_m)}
   (-i\, dt)
\\
   &\stackrel{\text{{\color{gray} 4}}}{=}
   \frac{d_{\ell m}}{2\pi}\int_{-1}^1
   \frac{(1-ita_\ell)(1-ita_m)}{(1+t^2a_\ell^2)(1+t^2a_m^2)}
   dt
\\
   &=\frac{d_{\ell m}}{2\pi}\int_{-1}^1
   \frac{1{\color{red}\,-it(a_\ell+a_m)}-t^2a_\ell a_m}{(1+t^2a_\ell^2)(1+t^2a_m^2)}
   dt
\\
   &\stackrel{\text{{\color{gray} 6}}}{=}
   \frac{d_{\ell m}}{2\pi}\int_{-1}^1
   \frac{1-t^2a_\ell a_m}{(1+t^2a_\ell^2)(1+t^2a_m^2)}
   dt
\\
   &=\frac{d_{\ell m}}{\pi}
   \begin{cases}
   \frac{\arctan(a_m)-\arctan(a_\ell)}{a_m-a_\ell}&\text{, if $\ell\not=m$,}\\
   \arctan^\prime(a_\ell)=\frac{1}{a_\ell^2+1}&\text{, if $\ell=m$ .}
   \end{cases}
\end{split}
\end{equation*}
In equality {\color{gray}3} we parametrize the path $\gamma_1$ by $[-1,1]\ni
t\mapsto -it \in [i,-i]\subset\C$.
In equality {\color{gray}4} we expand the fraction by the complex conjugate of
the denominator.
In equality {\color{gray}6} the {\color{red}imaginary part} vanishes during
integration by symmetry.
Finally integration was carried out by
\href{https://www.integralrechner.de}{integralrechner.de}.
\end{proof}

\boldmath
\subsubsection{Estimating the off-diagonal block $(T_A^\beta)^{-+}$ alias $H_{\mat{b}^{-+}}$}
\unboldmath

The main result of this subsection is

\begin{proposition}\label{prop:Had-+-block}
Hadamard multiplication
\begin{equation*}
\begin{split}
   H_{\mat{b}^{-+}}\colon \Ll(H_1^-,H_0^+)
   \to \Ll(H_{{\color{red}\frac12}}^-, H_{{\color{red}\frac12}}^+),\quad
   \mat{d}^{-+}\mapsto \mat{b}^{-+}\odot\mat{d}^{-+}
\end{split}
\end{equation*}
is a continuous linear map with bound
$\norm{H_{\mat{b}^{-+}}}_{\Ll(\Ll(H_1^-,H_0^+),
\Ll(H_{{\color{red}\frac12}}^-), H_{{\color{red}\frac12}}^+))}\le \pi/2$.
The same estimate holds true for the operator $H_{\mat{b}^{+-}}$.
\end{proposition}

Note that the constant in Proposition~\ref{prop:Had-+-block}
is universal and does not depend on the growth types
$\gf{f}_A^+$ and $\gf{f}_A^-$
of the positive and negative eigenvalues of $A$.

\medskip
We prove Proposition~\ref{prop:Had-+-block}
for $H_{\mat{b}^{-+}}$, the case $H_{\mat{b}^{+-}}$ being analogous
by interchanging the roles of $a_\mu$ and $b_\nu$:
Recall that for $\mu\in\N$ we have $a_\mu>0$.
We define $b_\nu:=-a_{-\nu}>0$ for $\nu\in\N$.
We further define
\begin{equation}\label{eq:b-1-+}
   b_{\mu\nu}^{-+}
   :=b_{\mu(-\nu)}
   =
   \frac{\arctan a_\mu+\arctan b_\nu}{a_\mu+b_\nu}
\end{equation}
for all $\mu,\nu\in\N$.

To prove the proposition we need the following lemma and isometries.

\begin{lemma}\label{le:tilde-b}
The matrix $\mat{\tilde b}$ whose entries for $\mu,\nu\in\N$ are defined by
\begin{equation}\label{eq:b-tilde}
   \tilde b_{\mu\nu}
   :=\sqrt{a_\mu b_\nu} \cdot b_{\mu\nu}^{-+}
   =\frac{\sqrt{a_\mu b_\nu}}{a_\mu+b_\nu}\left(\arctan a_\mu+\arctan b_\nu\right)
\end{equation}
is a Schur multiplier 
whose Schur operator satisfies
$\norm{S_{\mat{\tilde b}}}_{\Ll(\Ll(\ell^2))}\le\frac{\pi}{2}$.
\end{lemma}

\begin{proof}
Example~\ref{ex:Schur} part~III.
\end{proof}

\begin{definition}[Isometries]\label{def:isometries}
Given the standard orthonormal bases
$(e_\nu)_{\nu\in\N}$ of $\ell^2=\ell^2(\N)$
and $\Vv=(v_\ell)_{\ell\in\Z^*}$ of $H_0$, see~(\ref{eq:ONB-Vv-H_0}),
we introduce for $r\in\R$ the isometries
$$
   \Psi^+_r\colon\ell^2\to H^+_r
   ,\quad
   e_\mu\mapsto \tfrac{1}{{a_\mu}^r}\, v_\mu
$$
and
$$
   \Psi^-_r\colon\ell^2\to H^-_r
   ,\quad
   e_\nu\mapsto \tfrac{1}{{b_\mu}^r}\, v_{-\nu} .
$$
\end{definition}

\begin{proof}[Proof of Proposition~\ref{prop:Had-+-block}]
Given $\mat{d}^{-+}\in \Ll(H_1^-,H_0^+)$, we use the isometries $\Psi$
to define a map
$$
   \mat{D}
   :=(\Psi_0^+)^{-1}\mat{d}^{-+}\Psi_1^-
   \colon\ell^2
   \stackrel{\Psi_1^-}{\longrightarrow} H_1^-
   \stackrel{\mat{d}^{-+}}{\longrightarrow} H_0^+
   \stackrel{(\Psi_0^+)^{-1}}{\longrightarrow} \ell^2 .
$$
By the isometry property of the $\Psi$'s the operator norms are equal
$$
   \norm{\mat{D}}_{\Ll(\ell^2)}
   =\norm{\mat{d}^{-+}}_{\Ll(H_1^-,H_0^+)} .
$$
The matrix elements of $\mat{D}$ and the entries $d_{\mu\nu}^{-+}$ defined by
$$
   D_{\mu\nu}:=\INNER{e_\mu}{D e_\nu}_{\ell^2_\N}
   ,\qquad
   d_{\mu\nu}^{-+}:=\INNER{v_\mu}{\mat{d}^{-+} v_{-\nu}}_{\ell^2_{\Z^*}} ,
$$
for $\mu,\nu\in\N$ are related by
\begin{equation}\label{eq:D_munu}
\begin{split}
   D_{\mu\nu}
   &=\INNER{e_\mu}{(\Psi_0^+)^{-1}\mat{d}^{-+}\Psi_1^- e_\nu}_{\ell^2_\N}\\
   &\stackrel{{\color{gray}2}}{=}
   \INNER{\Psi_0^+ e_\mu}{\mat{d}^{-+}\Psi_1^- e_\nu}_{\ell^2_{\Z^*}}\\
   &\stackrel{{\color{gray}3}}{=}
   \INNER{\tfrac{1}{{a_\mu}^0}\, v_\mu}
   {\mat{d}^{-+}\tfrac{1}{{b_\nu}^1}\, v_\nu}_{\ell^2_{\Z^*}}\\
   &=\tfrac{1}{b_\nu} d_{\mu\nu}^{-+}.
\end{split}
\end{equation}
Here identity {\color{gray}2} is by the isometry property of $\Psi_0^+$
and identity {\color{gray}3} is by Definition~\ref{def:isometries}
of $\Psi_0^+$ and $\Psi_1^-$.

Since $\mat{D}$ is a bounded linear map on $\ell^2$
and $\mat{\tilde b}$ is a Schur multiplier by Lemma~\ref{le:tilde-b},
we obtain a bounded linear map $\mat{C}$ on $\ell^2$
characterized by the condition that
$\INNER{e_\mu}{\mat{C} e_\nu}_{\ell^2}=C_{\mu\nu}:=D_{\mu\nu}\tilde b_{\mu\nu}$
for all $\mu,\nu\in\N$. Moreover, as by Lemma~\ref{le:tilde-b} the
norm of $\mat{\tilde b}$ as a Schur multiplier is bounded above by
$\tfrac{\pi}{2} $, we get
$$
   \norm{\mat{C}}_{\Ll(\ell^2)}
   =\norm{S_{\mat{\tilde b}}(\mat{D})}_{\Ll(\ell^2)}
   \le\tfrac{\pi}{2} \norm{\mat{D}}_{\Ll(\ell^2)}
   =\tfrac{\pi}{2} \norm{\mat{d}^{-+}}_{\Ll(H_1^-,H_0^+)} .
$$
Observe that
\begin{equation}\label{eq:C_munu}
\begin{split}
   C_{\mu\nu}:
   =
   D_{\mu\nu}\tilde b_{\mu\nu}
   \stackrel{\text{(\ref{eq:D_munu})}}{=}
   \tfrac{1}{b_\nu} d_{\mu\nu}^{-+} \tilde b_{\mu\nu}
   \stackrel{\text{(\ref{eq:b-tilde})}}{=}   
   =\sqrt{\tfrac{a_\mu}{b_\nu}} d_{\mu\nu}^{-+} b_{\mu\nu}^{-+} .
\end{split}
\end{equation}
Here identity two is by the calculation above and identity three by
definition of $\tilde b_{\mu\nu}$ in Lemma~\ref{le:tilde-b}.
We use again the isometries $\Psi$ to define a linear map
$$
   \mat{C}^{-+}:=\Psi_{\frac12}^+ \mat{C}(\Psi_{\frac12}^-)^{-1}  
   \colon H^-_{\frac12}
   \stackrel{(\Psi_{\frac12}^-)^{-1}}{\longrightarrow} \ell^2
   \stackrel{\mat{C}}{\longrightarrow} \ell^2
   \stackrel{\Psi_{\frac12}^+}{\longrightarrow} H^+_{\frac12} .
$$
By the isometry property of the $\Psi$'s the operator norms are equal, hence
\begin{equation}\label{eq:C-+}
   \norm{\mat{C}^{-+}}_{\Ll(H_{\frac12}^-,H_{\frac12}^+)}
   =\norm{\mat{C}}_{\Ll(\ell^2)}
   \le \tfrac{\pi}{2} \norm{\mat{d}^{-+}}_{\Ll(H_1^-,H_0^+)} .
\end{equation}
The matrix entries of $\mat{C}^{-+}$ are for $\mu,\nu\in\N$ defined and given by
\begin{equation*}
\begin{split}
   C_{\mu\nu}^{-+}:
   &=\INNER{v_\mu}{\mat{C}^{-+} v_{-\nu}}_{\ell^2_{\Z^*}}\\
   &\stackrel{{\color{gray}2}}{=}\big\langle v_\mu , \Psi_{1/2}^+
   \underbrace{\mat{C}\,(\Psi_{1/2}^-)^{-1} v_{-\nu}}
      _{\mat{C}\sqrt{b_\nu} e_\nu=\sqrt{b_\nu}\sum_\lambda C_{\lambda\nu} e_\lambda}
   \big\rangle_{\ell^2_{\Z^*}}\\
   &\stackrel{{\color{gray}3}}{=}
   \Big\langle
   v_\mu ,
   \sqrt{b_\nu}\sum_{\lambda\in\N} C_{\lambda\nu}
   \underbrace{\Psi_{1/2}^+ e_\lambda}_{\tfrac{1}{\sqrt{a_\lambda}} v_\lambda}
   \Big\rangle_{\ell^2_{\Z^*}}\\
   &=\sqrt{b_\nu}\sum_{\lambda\in\N}
   \frac{C_{\lambda\nu}}{\sqrt{a_\lambda}}
   \underbrace{\INNER{v_\mu}{v_{\lambda}}_{\ell^2_{\Z^*}}}_{\delta_{\mu\lambda}}
   \\
   &=\sqrt{b_\nu}\frac{C_{\mu\nu}}{\sqrt{a_\mu}}\\
   &\stackrel{{\color{gray}6}}{=} 
   d_{\mu\nu}^{-+} b_{\mu\nu}^{-+} .
\end{split}
\end{equation*}
Here identity {\color{gray} 2} is by definition of $\mat{C}^{-+}$
and identity {\color{gray} 3} is by Definition~\ref{def:isometries}
of $\Psi_{1/2}^-$ and $\Psi_{1/2}^+$ and by expanding
$\mat{C} e_\nu=\sum_{\lambda\in\N} C_{\lambda\nu} e_\lambda$.
Identity {\color{gray} 6} holds by~(\ref{eq:C_munu}).

\smallskip
Thus $\mat{C}^{-+}$ is the Hadamard product
$\mat{C}^{-+}=\mat{b}^{-+} \odot \mat{d}^{-+} =:H_{\mat{b}^{-+}}(\mat{d}^{-+})$.
Hence
$$
   \norm{H_{\mat{b}^{-+}}(\mat{d}^{-+})}_{\Ll(H_\frac12^-,H_\frac12^+)}
   \stackrel{\text{(\ref{eq:C-+})}}{\le}
   \tfrac{\pi}{2} \norm{\mat{d}^{-+}}_{\Ll(H_1^-,H_0^+)} .
$$
by~(\ref{eq:C-+}).
Therefore Hadamard multiplication by $\mat{b}^{-+}$ is a bounded linear map
$$
   H_{\mat{b}^{-+}}\colon \Ll(H_1^-,H_0^+)
   \to \Ll(H_{{\color{red}\frac12}}^-, H_{{\color{red}\frac12}}^+)
$$
of norm
$$
   \norm{H_{\mat{b}^{-+}}}_{\Ll(\Ll(H_1^-,H_0^+),\Ll(H_{{\color{red}\frac12}}^-, H_{{\color{red}\frac12}}^+))}
   \le \tfrac{\pi}{2}.
$$
This concludes the proof of Proposition~\ref{prop:Had-+-block}.
\end{proof}

\begin{remark}\label{rmk:obstruction}
Because of Corollary~\ref{cor:obstruction}
the proof of Proposition~\ref{prop:Had-+-block}
cannot be strengthened to give an extension to $\Ll(H_0,H_0)$.
\end{remark}

\boldmath
\subsection*{The operator sum $T_A^\beta+T_A^{\alpha_+}=d\Pi_+(A)$}
\unboldmath

\boldmath
\subsubsection{Proof of Theorem~\ref{thm:dPi+-ext}}
\unboldmath

Consider the operator $d\Pi_+(A)=T_A^\beta+T_A^{\alpha_+}$
given by~(\ref{eq:T_A^beta+T_A^alpha})
and suppose that $\Delta\in\Ll_{sym_0}(H_1,H_0)$.
By Corollary~\ref{cor:bjkhbjg6544}
the diagonal blocks vanish
$$
   d\Pi_+(A)\Delta
   \stackrel{\text{(\ref{eq:4block})} \atop \text{(\ref{eq:4block-2})}}{=}
   \begin{pmatrix}
   0&(T_A^\beta)^{-+}\Delta+(T_A^{\alpha_+})^{-+}\Delta\\
   (T_A^\beta)^{+-}\Delta+(T_A^{\alpha_+})^{+-}\Delta&0
   \end{pmatrix} .
$$
We use Lemma~\ref{le:4block} to estimate
the operator norm by the norms of the four blocks,
then we use the triangle inequality to get inequality two
\begin{equation*}
\begin{split}
   &\norm{d\Pi_+(A)\Delta}_{\Ll(H_{\frac12})}
\\
   &\le
   \norm{(T_A^\beta)^{-+}\Delta+(T_A^{\alpha_+})^{-+}\Delta}
      _{\Ll(H_{\frac12}^-, H_{\frac12}^+)}\\
   &\quad
   +\norm{(T_A^\beta)^{+-}\Delta+(T_A^{\alpha_+})^{+-}\Delta}
      _{\Ll(H_{\frac12}^+, H_{\frac12}^-)}
\\
   &\le
   \underbrace{\norm{(T_A^\beta)^{-+}\Delta}_{\Ll(H_{\frac12}^-,H_{\frac12}^+)}}
      _{\text{$\le\frac{\pi}{2}$ by Prop.~\ref{prop:Had-+-block}}}
   +\underbrace{\norm{(T_A^{\alpha_+})^{-+}\Delta}_{\Ll(H_{\frac12}^-,H_{\frac12}^+)}}
      _{\text{$\le \frac{2\sigma+1}{\pi\sigma^2}$ by Prop.~\ref{prop:T_A^alpha}}}
\\
   &\quad
   +\underbrace{\norm{(T_A^\beta)^{+-}\Delta}_{\Ll(H_{\frac12}^+,H_{\frac12}^-)}}
      _{\text{$\le\frac{\pi}{2}$ by Prop.~\ref{prop:Had-+-block}}}
   +\underbrace{\norm{ (T_A^{\alpha_+})^{+-}\Delta}_{\Ll(H_{\frac12}^+,H_{\frac12}^-)}}
      _{\text{$\le \frac{2\sigma+1}{\pi\sigma^2}$ by Prop.~\ref{prop:T_A^alpha}}}
\\
   &\le\pi+\tfrac{4\sigma+2}{\pi\sigma^2}.
\end{split}
\end{equation*}
Here we also used that the norm of each block is bounded from above by
the norm of the total operator, in symbols
$\norm{(T_A^{\alpha_+})^{-+}\Delta}_{\Ll(H_{\frac12}^-,H_{\frac12}^+)}
\le \norm{T_A^{\alpha_+}\Delta}_{\Ll(H_{\frac12})}$.
This proves Theorem~\ref{thm:dPi+-ext}.

\boldmath
\subsection{We prove Theorem~\ref{thm:B}}
\label{sec:thmB-2}
\unboldmath

Let $A$ be in the open subset $\Ff_\gf{h}^*$ of $\Ff_\gf{h}$
consisting of invertible weak Hessians. Invertibility is an open condition.
Since eigenvalues depend continuously on $A$ we can
choose an open convex neighborhood 
$\Uu$ of $A$ in $\Ff_\gf{h}^*$ with the property that there exists a
\textbf{uniform spectral gap} $\sigma_0>0$ in the sense that the
spectral gap of every $B\in\Uu$ is at least $\sigma_0$.
Let $\Delta \in\Ff_\gf{h}$ be such that $A+\Delta\in \Uu$.
Since $\Uu$ is convex it follows that $A+t\Delta\in \Uu$ for every $t\in[0,1]$
and therefore the spectral gap of $A+t\Delta$ is at least $\sigma_0$.

Recall from~(\ref{eq:dPi+-ext}) that $T_A^+$ denotes the extension of
$d\Pi_+(A)$ from $\Ll(H_1,H_0)$ to $\Ll(H_{\frac12})$.
The fundamental theorem of calculus is used in the first inequality
\begin{equation*}
\begin{split}
   &\norm{T_{A+\Delta}^+-T_A^+}_{\Ll(H_{\frac12})}
\\
   &\le\int_0^1
   \norm{T_{A+t\Delta}^+}_{\Ll(\Ll_{sym_0}(H_1,H_0), \Ll(H_{\frac12}))}
   \norm{\Delta}_{\Ll(H_1,H_0)}\, dt
\\
   &\le
   \left(\pi+\frac{4\sigma_0+2}{\pi\sigma_0^2}\right) 
   \norm{\Delta}_{\Ll(H_1,H_0)}
\end{split}
\end{equation*}
and inequality two holds by Theorem~\ref{thm:dPi+-ext}.
This proves continuity in $\Ll(H_{\frac12})$.
By scale shift invariance, as explained in Example~\ref{ex:h-pair},
continuity in $\Ll(H_{\frac{3}{2}})$ follows as well.
This completes the proof of Theorem~\ref{thm:B}.

\appendix
\section{Hilbert space pairs}\label{sec:HS-pairs}

\begin{definition}
A \textbf{Hilbert space pair} $H=(H_0,H_1)$ consists of
\emph{infinite} dimensional\footnote{
  In finite dimension the dense inclusion would imply that $H_1=H_0$.
  }
Hilbert spaces such that $H_1$ is a dense subset of $H_0$
and the inclusion map is compact.
An \textbf{isomorphism of Hilbert space pairs} $H=(H_0,H_1)$  and
$W=(W_0,W_1)$ is a Hilbert space isomorphism  $T\colon H_0\to W_0$
whose restriction to $H_1$, notation $T_1$, takes values in $W_1$
and such that $T_1\colon H_1\to W_1$ is a Hilbert space isomorphism.
Two Hilbert space pairs are called \textbf{isomorphic}
if there exists an isomorphism of Hilbert space pairs.
\end{definition}

\begin{definition}
An \textbf{isometry of Hilbert space pairs}
is an isomorphism $\Phi \colon P\to Q$ of Hilbert space pairs
with the additional property that $\Phi \colon H_0\to W_0$
as well as its restriction $\Phi _1\colon H_1\to W_1$ are Hilbert space
isometries.
Two Hilbert space pairs are called \textbf{isometric}
if there exists an isometry of Hilbert space pairs.
\end{definition}

\begin{example}[Scale shift invariance]\label{ex:h-pair}
Let $\gf{h}$ be a growth function. Consider the Hilbert space pair
$H=(H_0,H_1)=(\ell^2,\ell^2_{\gf{h}})$ and,
for any real $r\in\R$, consider the Hilbert space pair
$W=(W_0,W_1)=(\ell^2_{\gf{h}^{r}},\ell^2_{\gf{h}^{r+1}})$.
These two pairs are isometric Hilbert pairs.
An isometry is given by the linear map determined on basis vectors by
$$
   \Phi\colon H_0\to W_0
   ,\quad
   e_\nu\mapsto \gf{h}(\nu)^{-r/2} e_\nu .
$$
\end{example}

\begin{theorem}\label{thm:HS_pair-growth_fct}
Given a Hilbert space pair $H=(H_0,H_1)$, there exists a growth function
$\gf{h}\colon\N\to(0,\infty)$
such that the pair $(H_0,H_1)$ is isometric to $(\ell^2,\ell^2_{\gf{h}})$.
\end{theorem}

\begin{corollary}\label{cor:HS_pair-growth_fct}
In a Hilbert space pair $(H_0,H_1)$ both Hilbert spaces $H_0$ and
$H_1$ are separable.
\end{corollary}

\begin{remark}[Uniqueness of pair growth]\label{rmk:pair-gf-unique}
The pair growth function $\gf{h}=\gf{h}(H_0,H_1)$ is unique up to pair
isometry.
The growth type $[\gf{h}]$ is unique up to pair isomorphism
as we show in Theorem~\ref{thm:f-g-equiv} below.
\end{remark}

\begin{proof}[Proof of Theorem~\ref{thm:HS_pair-growth_fct}]
On $H_1$ we have two inner products, namely $\INNER{\cdot}{\cdot}_1$
and the restriction to $H_1$ of the inner product of $H_0$,
still denoted by $\INNER{\cdot}{\cdot}_0$.
By the Theorem of Riesz there exists a bounded linear map
$T\colon H_1\to H_1$ such that
\begin{equation}\label{eq:Riesz}
   \INNER{\xi}{\eta}_0=\INNER{\xi}{T\eta}_1
\end{equation}
for all $\xi,\eta\in H_1$.
We need the following lemma.

\begin{lemma}\label{le:compact}
The operator $T\in\Ll(H_1)$ is compact, symmetric, positive definite.
\end{lemma}

\begin{proof}[Proof of Lemma~\ref{le:compact}]
Symmetry, resp. positive definiteness, of the inner product
$\INNER{\cdot}{\cdot}_1$ imply symmetry, resp. positive
definiteness, of the operator $T$.

To prove compactness of $T$ we
pick a unit ball sequence $\xi_\nu\in H_1$, that is $\norm{\xi_\nu}_1\le 1$.
By compactness of the inclusion $\iota\colon H_1\INTO H_0$ there
is a subsequence $\xi_{\nu_j}$ which converges in $H_0$.
\\
We show that $T\xi_{\nu_j}$ converges in $H_1$:
Since the sequence $\xi_{\nu_j}$ converges in $H_0$,
the sequence $\xi_{\nu_j}$ is Cauchy in $H_0$.
Given $\eps>0$, there exists $j_0=j_0(\eps)\in\N$
such that for all $j,j^\prime\ge j_0$ we have
$\norm{\xi_{\nu_j}-\xi_{\nu_{j^\prime}}}_0\le\eps/\sqrt{\norm{T}}$
where $\norm{T}=\norm{T}_{\Ll(H_1)}$.

In the following we show that $\norm{T(\xi_{\nu_j}-\xi_{\nu_{j^\prime}})}_1\le \eps$.
For this purpose we abbreviate $v:=\xi_{\nu_j}-\xi_{\nu_{j^\prime}}$.
Hence $\norm{v}_0^2\le \eps^2/\norm{T}$.
We estimate
\begin{equation*}
\begin{split}
   0
   &\le\Norm{v-\norm{T}^{-1} Tv}_0^2
\\
   &=\INNER{v-\norm{T}^{-1} Tv}{v-\norm{T}^{-1} Tv}_0
\\
   &=\Norm{v}_0^2-\tfrac{2}{\norm{T}} \INNER{v}{Tv}_0
   +\tfrac{1}{\norm{T}^2}\INNER{Tv}{Tv}_0
\\
   &\stackrel{{\color{gray} 4}}{=}
   \Norm{v}_0^2
   -\tfrac{2}{\norm{T}} \INNER{Tv}{Tv}_1
   +\tfrac{1}{\norm{T}^2}\INNER{Tv}{TTv}_1
\\
   &\stackrel{{\color{gray} 5}}{\le} \tfrac{\eps^2}{\norm{T}}
   -\tfrac{2}{\norm{T}} \Norm{Tv}_1^2
   +\tfrac{1}{\norm{T}^2} \Norm{Tv}_1\underbrace{\Norm{TTv}_1}_{\le\norm{T}\norm{Tv}_1}
\\
   &\le  \tfrac{\eps^2}{\norm{T}}
   -\tfrac{1}{\norm{T}} \Norm{Tv}_1^2
\end{split}
\end{equation*}
and this proves that $\norm{Tv}_1\le\eps$.
In step {\color{gray} 4} we used~(\ref{eq:Riesz}) and symmetry of $T$.
In step {\color{gray} 5} we used $\norm{v}_0^2\le \eps^2/\norm{T}$
and Cauchy-Schwarz. Thus the sequence $T\xi_{\nu_j}$ is Cauchy in
$H_1$, but $H_1$ is complete. This proves Lemma~\ref{le:compact}.
\end{proof}

By Lemma~\ref{le:compact}
the spectral theorem applies and yields
an orthonormal basis 
\begin{equation}\label{eq:E-T}
   \Ee:=\{E_\nu\}_{\nu\in\N}\subset H_1
   ,\qquad
   TE_\nu=\kappa_\nu E_\nu,\quad \kappa_\nu\searrow 0,
\end{equation}
of $H_1$ consisting of eigenvectors of $T\colon H_1\to H_1$
whose eigenvalues are positive reals and
form a monotone decreasing sequence converging to zero.

\medskip
\noindent
\textbf{Claim.}
The appropriately rescaled $E_\nu$'s, namely
\begin{equation}\label{eq:e-E-T}
   e_\nu:=\tfrac{1}{\sqrt{\kappa_\nu}}\, E_\nu
   {\color{gray}\;\in H_1\subset H_0 ,}
\end{equation}
form an orthonormal basis of $H_0$.

\begin{proof}[Proof of the claim]
Note that $T e_\nu=\kappa_\nu e_\nu$.
Use~(\ref{eq:Riesz}) in the first step to get
\begin{equation*}
\begin{split}
   \INNER{e_\nu}{e_\mu}_0
   &=\INNER{e_\nu}{Te_\mu}_1
\\
   &=\INNER{e_\nu}{\kappa_\mu e_\mu}_1
\\
   &=\INNER{\tfrac{1}{\sqrt{\kappa_\nu}}\, E_\nu}
   {\kappa_\mu\tfrac{1}{\sqrt{\kappa_\mu}}\, E_\mu}_1
\\
   &=\sqrt{\tfrac{\kappa_\mu}{\kappa_\nu}}\INNER{E_\nu}{E_\mu}_1
\\
   &=\sqrt{\tfrac{\kappa_\mu}{\kappa_\nu}} \delta_{\nu\mu}
\\
   &=\delta_{\nu\mu}.
\end{split}
\end{equation*}
This proves that $e_\nu$ and $e_\mu$ are orthonormal in $H_0$.
Since the inclusion of $H_1$ in $H_0$ is dense and the $E_\nu$ form a basis
of $H_1$ it follows that the $e_\nu$ form an orthonormal basis of $H_0$.
This proves the claim.
\end{proof}

By the claim the following linear map is an isometry
\begin{equation}\label{eq:Phi-H}
   \Phi\colon \ell^2\to H_0,\quad
   (x_\nu)\mapsto \sum_{\nu=1}^\infty x_\nu e_\nu .
\end{equation}
The \textbf{growth function associated to the Hilbert pair \boldmath$H$}
is defined by
\begin{equation}\label{eq:f_H}
\boxed{
   \gf{h}(\nu):=\tfrac{1}{\kappa_\nu} .
}
\end{equation}
We calculate
\begin{equation*}
\begin{split}
   \INNER{e_\nu}{e_\mu}_1
   &=\INNER{\tfrac{1}{\sqrt{\kappa_\nu}}  E_\nu}
   {\tfrac{1}{\sqrt{\kappa_\mu}} E_\mu}_1
\\
   &=\tfrac{1}{\sqrt{\kappa_\nu\kappa_\mu}}\INNER{E_\nu}{E_\mu}_1
\\
   &=\tfrac{1}{\sqrt{\kappa_\nu\kappa_\mu}}\,\delta_{\nu\mu}
\\
   &=\tfrac{1}{\kappa_\nu}\,\delta_{\nu\mu}
\\
   &=\gf{h}(\nu) \delta_{\nu\mu}.
\end{split}
\end{equation*}
Hence $\Phi$ restricts to an isometry $\Phi|\colon\ell^2_{\gf{h}}\to H_1$.
This finishes the proof of Theorem~\ref{thm:HS_pair-growth_fct}.
\end{proof}

\begin{theorem}\label{thm:f-g-equiv}
Two Hilbert space pairs $(\ell^2,\ell^2_\gf{f})$ and
$(\ell^2,\ell^2_\gf{g})$ are isomorphic
iff the growth functions $\gf{f}$ and $\gf{g}$ are equivalent.
\end{theorem}

\begin{proof}
'$\Leftarrow$' Suppose $\gf{f}\sim \gf{g}$.
Then the inner products~(\ref{eq:ell^2_f}) on $\ell^2_\gf{f}$ and
$\ell^2_\gf{g}$ are equivalent.
Then the identity map $\Id\colon \ell^2_\gf{f}\to\ell^2_\gf{g}$
is an isomorphism of Hilbert spaces.
Hence $\Id\colon (\ell^2,\ell^2_\gf{f}) \to
(\ell^2,\ell^2_\gf{g})$, $x\mapsto x$,
is a Hilbert space pair isomorphism.

\smallskip
\noindent
'$\Rightarrow$'
Suppose there is a Hilbert space pair isomorphism
$$
   T\colon (\ell^2,\ell^2_\gf{f})\to
   (\ell^2,\ell^2_\gf{g}).
$$
Since $T\colon \ell^2\to\ell^2$ is an isomorphism,
there is a constant $c_0$ such that
\begin{equation}\label{eq:T}
   \frac{1}{c_0}\norm{x}_{\ell^2}\le \norm{Tx}_{\ell^2}
   \le c_0 \norm{x}_{\ell^2},\quad \forall x\in\ell^2.
\end{equation}
Since the restriction $T_1=T|\colon
\ell^2_\gf{f}\to\ell^2_\gf{g}$ is an isomorphism
and choosing $c_0$ larger, if necessary, it holds that
\begin{equation}\label{eq:T-2}
   \frac{1}{c_0}\norm{x}_{\ell^2_\gf{f}}\le \norm{Tx}_{\ell^2_\gf{g}}
   \le c_0 \norm{x}_{\ell^2_\gf{f}},\quad \forall x\in\ell^2_\gf{f}.
\end{equation}
The map $\Pi_n\colon \ell^2\to\ell^2$ projects
a sequence $x=(x_i)$ to its first $n$ members and sets all others
equal zero, in symbols $x\mapsto (x_1,\dots,x_n,0,\dots)$.
For natural numbers $n,m\in\N$ define the linear map
$$
   A^n_m:=\Pi_nT^{-1}\Pi_{m-1}T\colon \Pi_n\ell^2\to\Pi_n \ell^2 .
$$

\medskip
\noindent
\textbf{Claim.}
If $\gf{f}(n)<\frac{\gf{g}(m)}{{c_0}^4}$
then $A^n_m$ is bijective, hence a vector space isomorphism.

\begin{proof}[Proof of the claim].
Since $\Pi_n\ell^2$ is a finite dimensional vector space
it suffices to show that $A^n_m$ is injective.
To show this we pick $\xi\in \ker A^n_m$.
Therefore $\xi\in\ell^2$ satisfies the following equations
$A^n_m\xi=0$ and since $\xi$ lies in the image of the projection $\Pi_n=\Pi_n^2$
we have in addition that $\xi=\Pi_n\xi$.
Hence we obtain
\begin{equation*}
\begin{split}
   \xi
   &=\Pi_n\xi\\
   &=\Pi_nT^{-1}T\xi\\
   &=\Pi_nT^{-1}T\xi-A^n_m\xi\\
   &=\Pi_nT^{-1}T\xi-\Pi_nT^{-1}\Pi_{m-1}T\xi\\
   &=\Pi_nT^{-1}\left(\Id-\Pi_{m-1}\right) T\xi .
\end{split}
\end{equation*}
Now we estimate
\begin{equation*}
\begin{split}
   \Norm{\xi}_{\ell^2}
   &=\Norm{\Pi_nT^{-1}\left(\Id-\Pi_{m-1}\right) T\xi}_{\ell^2}\\
   &\le \Norm{\Pi_n}_{\Ll(\ell^2,\ell^2)}
   \Norm{T^{-1}}_{\Ll(\ell^2,\ell^2)}
   \Norm{\left(\Id-\Pi_{m-1}\right) T\xi}_{\ell^2} .
\end{split}
\end{equation*}
Observe that $\norm{\Pi_n}_{\Ll(\ell^2,\ell^2)}=1$, by orthogonality
of the projection $\Pi_n$, and $\norm{T^{-1}}_{\Ll(\ell^2,\ell^2)}\le c_0$
by the first inequality in~(\ref{eq:T}).
Therefore
\begin{equation}\label{eq:id-Pi-ell2}
   \Norm{\xi}_{\ell^2}
   \le c_0 \Norm{\left(\Id-\Pi_{m-1}\right) T\xi}_{\ell^2} .
\end{equation}
We abbreviate
\begin{equation}\label{eq:x-T}
   x:=\left(\Id-\Pi_{m-1}\right) T\xi .
\end{equation}
Since $\Id-\Pi_{m-1}$ is the orthogonal projection
which makes the first $m-1$ entries zero,
the element $x$ is of the form
\[
   x=\left(0,\dots,0,x_m,x_{m+1},\dots\right) .
\]
By definition of the $\ell^2_\gf{g}$ norm,
compare~(\ref{eq:ell^2_f}), and using that $\gf{g}$ is monotone
increasing, we have
\begin{equation*}
\begin{split}
   \Norm{x}_{\ell^2_\gf{g}}^2
   &=\sum_{i=m}^\infty \gf{g}(i) x_i^2 \\
   &\ge \gf{g}(m) \sum_{i=m}^\infty x_i^2 \\
   &=\gf{g}(m) \Norm{x}_{\ell^2}^2 .
\end{split}
\end{equation*}
Therefore any $x$ of the form
$\left(0,\dots,0,x_m,x_{m+1},\dots\right)$ satisfies the inequality
$$
   \Norm{x}_{\ell^2}
   \le\frac{1}{\sqrt{\gf{g}(m)}} \Norm{x}_{\ell^2_\gf{g}} .
$$
Using the particular form~(\ref{eq:x-T}) of $x$ we have
\begin{equation*}
\begin{split}
   \Norm{\left(\Id-\Pi_{m-1}\right) T\xi}_{\ell^2}
   &\le \frac{1}{\sqrt{\gf{g}(m)}}
   \Norm{\left(\Id-\Pi_{m-1}\right) T\xi}_{\ell^2_\gf{g}} \\
   &\le \frac{1}{\sqrt{\gf{g}(m)}}
   \Norm{\Id-\Pi_{m-1}}_{\Ll(\ell^2_\gf{g}, \ell^2_\gf{g})}
   \Norm{T}_{\Ll(\ell^2_\gf{f}, \ell^2_\gf{g})}
   \Norm{\xi}_{\ell^2_\gf{f}}
\end{split}
\end{equation*}
Observe that $\norm{T}_{\Ll(\ell^2_\gf{f}, \ell^2_\gf{g})}\le c_0$,
by the second inequality in~(\ref{eq:T-2}), and that
$\norm{\Id-\Pi_{m-1}}_{\Ll(\ell^2_\gf{g}, \ell^2_\gf{g})}\le 1$
since the projection $\Id-\Pi_{m-1}$ is also orthogonal in $\ell^2_\gf{g}$.
Thus
\begin{equation}\label{eq:id-Pi}
   \Norm{\left(\Id-\Pi_{m-1}\right) T\xi}_{\ell^2}
   \le
   \frac{c_0}{\sqrt{\gf{g}(m)}}\Norm{\xi}_{\ell^2_\gf{f}} .
\end{equation}
Using that $\xi$ lies in the image of $\Pi_n$
it is of the form $\xi=\left(\xi_1,\dots,\xi_n,0,\dots\right)$.
\begin{equation*}
\begin{split}
   \Norm{\xi}_{\ell^2_\gf{f}}^2
   &=\sum_{i=1}^n \gf{f}(i) \xi_i^2 \\
   &\le \gf{f}(n) \sum_{i=1}^n \xi_i^2 \\
   &=\gf{f}(n)\Norm{\xi}_{\ell^2}^2 .
\end{split}
\end{equation*}
Therefore
\begin{equation}\label{eq:xi-ell2}
   \Norm{\xi}_{\ell^2_\gf{f}}
   \le \sqrt{\gf{f}(n)}\Norm{\xi}_{\ell^2} .
\end{equation}
Therefore we are now in position to estimate
\begin{equation*}
\begin{split}
   \Norm{\xi}_{\ell^2}
   &\stackrel{(\ref{eq:id-Pi-ell2})}{\le}
   c_0 \Norm{\left(\Id-\Pi_{m-1}\right) T\xi}_{\ell^2}
\\
   &\stackrel{(\ref{eq:id-Pi})}{\le}
   \frac{c_0^2}{\sqrt{\gf{g}(m)}}\Norm{\xi}_{\ell^2_\gf{f}}
\\
   &\stackrel{(\ref{eq:xi-ell2})}{\le}
   a_{nm}\Norm{\xi}_{\ell^2}
   ,\quad
   a_{nm}:=c_0^2\sqrt{\frac{\gf{f}(n)}{\gf{g}(m)}} .
\end{split}
\end{equation*}
By assumption of the claim we have $a_{nm}<1$.
Therefore $\norm{\xi}_{\ell^2}=0$, hence $\xi=0$
This means that $A^n_m$ is injective and therefore
a vector space isomorphism of the finite ($n$) dimensional vector space
$\Pi_n\ell^2$ in itself. This proves the claim.
\end{proof}

Since, under the assumption of the claim, the composition
$A^n_m=(\Pi_nT^{-1})\circ(\Pi_{m-1}T)\colon \Pi_n\ell^2\to\Pi_n\ell^2$
is bijective, it follows that the second part
$B:=\Pi_nT^{-1}|_{\Pi_{m-1}\ell^2}\colon \Pi_{m-1}\ell^2\to \Pi_{n}\ell^2$
is surjective.
By the dimension theorem for finite dimensional vector spaces
$$
   m-1=\dim\Pi_{m-1}\ell^2=\dim\ker B+\dim \Pi_{n}\ell^2\ge \dim
   \Pi_{n}\ell^2=n .
$$
Thus $n<m$.

So far we have shown the following implication
\[
   \gf{f}(n)<\frac{\gf{g}(m)}{{c_0}^4}
   \quad\Rightarrow\quad
   n<m.
\]
Therefore at the same point $n=m$ both functions are related by
\[
   \gf{f}(n)\ge\frac{\gf{g}(n)}{{c_0}^4}.
\]
Note that $T^{-1}\colon\ell^2\to\ell^2$ is as well an isomorphism
satisfying the inequalities~(\ref{eq:T}) and~(\ref{eq:T-2}) with
$\gf{f}$ and $\gf{g}$ interchanged.
Hence we obtain as well the inequality
\[
   \gf{g}(n)\ge\frac{\gf{f}(n)}{{c_0}^4}.
\]
Now set $c:=c_0^4$ to obtain
\[
   \frac{1}{c} \gf{f}(n)\le \gf{g}(n)
   \le c \gf{f}(n).
\]
This concludes the proof of Theorem~\ref{thm:f-g-equiv}.
\end{proof}

\section{Interpolation}\label{sec:interpolation}

For convenience of the reader we include
the proof of the Stein interpolation theorem~\cite[Theorem~2]{Stein:1946a}
in our situation of the Hilbert spaces $\ell^2_{\gf{f}}$
with inner product~(\ref{eq:ell^2_f}).
The proof is based on the three line theorem
whose usage in interpolation theory goes back to Thorin;
we recommend the article~\cite{Bondesson:2008a}
on Thorin's life and work.

\begin{proposition}\label{prop:Stein}
Assume that $T$ is a linear map with bounds
\begin{equation*}
\begin{split}
   T\colon \ell^2
   &\to \ell^2,
   \qquad
   \Norm{T}_{\Ll(\ell^2,\ell^2)}\le M_0 ,
\\
   T\colon \ell^2_{\gf{f}}
   &\to \ell^2_{\gf{g}},
   \qquad
   \Norm{T}_{\Ll(\ell^2_{\gf{f}},\ell^2_{\gf{g}})}\le M_1 .
\end{split}
\end{equation*}
Then $T$ is also bounded as a linear map
\begin{equation*}
\begin{split}
   T\colon \ell^2_{\sqrt{\gf{f}}}
   &\to \ell^2_{\sqrt{\gf{g}}} ,
   \qquad
   \Norm{T}_{\Ll(\ell^2_{\sqrt{\gf{f}}},\ell^2_{\sqrt{\gf{g}}})}\le
   \sqrt{M_0 M_1}.
\end{split}
\end{equation*}
\end{proposition}

\begin{proof}
Pick two unit vectors $\xi\in\ell^2_{\sqrt{\gf{f}}}$
and $\eta\in\ell^2_{\sqrt{\gf{g}}}$.
For $k\in\N$ let $e_k$ be the sequence of reals whose
members are all zero except member $k$ which is $1$,
in symbols $e_k=(0,\dots,0,1,0,\dots)$.
Write $\xi=\sum_k \xi_ke_k$ and $\xi=\sum_\ell \xi_\ell e_\ell$
and set $T_{k,l}:=\INNER{Te_k}{e_\ell}_0$.
On the strip $S=\{\theta+it\mid\text{$\theta\in[0,1]$, $t\in\R$}\}$
consider the function $F\colon S\to\C$ defined by
\[
   F(\theta+it)
   :=\sum_{k,\ell} T_{k,\ell} \xi_k
   \gf{f}(k)^{\quarter-\frac12\theta-\frac12 it}
   \eta_\ell \gf{g}(\ell)^{\quarter+\frac12\theta+\frac12 it} .
\]
\medskip
\noindent
\textbf{Claim.}
(i)~$F(\frac12)=\INNER{T\xi}{\eta}_{\sqrt{\gf{g}}}$.
(ii)~$\abs{F(it)}\le M_0$.
(iii)~$\abs{F(1+it)}\le M_1$.

\begin{proof}[Proof of the claim]
(i)~We have $F(\frac12)=\sum_{k,\ell} T_{k,\ell} \xi_k\eta_\ell \sqrt{\gf{g}(\ell)}
=\INNER{T\xi}{\eta}_{\sqrt{\gf{g}}}$.
\\
(ii)~We have $F(it)=\sum_{k,\ell} T_{k,\ell} \xi_k
   \gf{f}(k)^{\quarter-\frac12 it}
   \eta_\ell \gf{g}(\ell)^{\quarter+\frac12 it}$.
The element
$$
   x:=\sum_{k} \underbrace{\xi_k \gf{f}(k)^{\quarter-\frac12
       it}}_{=: x_k} e_k\in\ell^2\otimes \C
$$
has unit norm square
$
   \norm{x}_{\ell^2}^2
   =\sum_{k} \abs{\xi_k \gf{f}(k)^{\quarter-\frac12 it}}^2
   = \sum_{k} \xi_k^2 \sqrt{\gf{f}(k)}
   =\norm{\xi}_{\sqrt{\gf{f}}}^2
   =1
$.
The element
$$
   y:=\sum_\ell \underbrace{\eta_\ell \gf{g}(\ell)^{\quarter+\frac12 it}}_{=: y_\ell} e_\ell\in\ell^2\otimes \C
$$
has unit norm square
$
   \norm{y}_{\ell^2}^2
   =\sum_\ell \abs{\eta_\ell \gf{g}(\ell)^{\quarter+\frac12 it}}^2
   = \sum_\ell \eta_\ell^2 \sqrt{\gf{g}(\ell)}
   =\norm{\eta}_{\sqrt{\gf{g}}}^2
   =1
$.
Now
$$
   \Abs{F(it)}
   =\Bigl|\sum_{k,\ell} T_{k,\ell} x_k y_\ell\Bigr|
   =\Abs{\INNER{Tx}{y}_{\ell^2}}
   \le\Norm{T}_{\Ll(\ell^2,\ell^2)}\Norm{x}_{\ell^2}\Norm{y}_{\ell^2}
   \le M_0 .
$$
(iii)~We have $F(1+it)
   :=\sum_{k,\ell} T_{k,\ell} \xi_k
   \gf{f}(k)^{-\quarter-\frac12 it}
   \eta_\ell \gf{g}(\ell)^{\frac{3}{4}+\frac12 it} .$
The element
$$
   x:=\sum_{k} \underbrace{\xi_k \gf{f}(k)^{-\quarter-\frac12
       it}}_{=: x_k} e_k\in\ell^2_{\gf{f}}\otimes \C
$$
has $\ell^2_{\gf{f}}$-norm square
$
   \norm{x}_{\ell^2_{\gf{f}}}^2
   =\sum_{k} \abs{\xi_k \gf{f}(k)^{-\quarter-\frac12 it}}^2 \gf{f}(k)
   = \sum_{k} \xi_k^2 \sqrt{\gf{f}(k)}
   =\norm{\xi}_{\sqrt{\gf{f}}}^2
   =1
$.
The element
$$
   y:=\sum_\ell \underbrace{\eta_\ell \gf{g}(\ell)^{-\quarter+\frac12 it}}_{=: y_\ell} e_\ell\in\ell^2_{\gf{g}}\otimes \C
$$
has $\ell^2_{\gf{g}}$-norm square
$
{   \norm{y}_{\ell^2_{\gf{g}}}^2
   =\sum_\ell \abs{\eta_\ell \gf{g}(\ell)^{-\quarter+\frac12 it}}^2\gf{g}(\ell)
   = \sum_\ell \eta_\ell^2 \sqrt{\gf{g}(\ell)}
   =\norm{\eta}_{\sqrt{\gf{g}}}^2
   =1
}
$.
Now
$$
   \Abs{F(1+it)}
   =\Bigl|\sum_{k,\ell} T_{k,\ell} x_k \gf{g}(\ell) y_\ell\Bigr|
   =\Abs{\INNER{Tx}{y}_{\ell^2_{\gf{g}}}}
   \le\Norm{T}_{\Ll(\ell^2_{\gf{f}},\ell^2_{\gf{g}})}\Norm{x}_{\ell^2_{\gf{f}}}\Norm{y}_{\ell^2_{\gf{g}}}
   \le M_1 .
$$
This proves the claim.
\end{proof}

\begin{lemma}[The three line theorem]
Let $F(z)$ be analytic on the open strip
$0<\Re z<1$ and bounded and continuous on the closed strip $0\le \Re
z\le 1$. If
$$
   \abs{F(it)}\le M_0,\qquad
   \abs{F(1+it)}\le M_1,\qquad -\infty<t<\infty ,
$$
then we have
$$
   \abs{F(\tfrac12+it)}\le \sqrt{M_0 M_1},\qquad -\infty<t<\infty .
$$
\end{lemma}

\begin{proof}
See e.g.~\cite[Le.\,1.1.2]{Bergh:1976a}.
\end{proof}

By parts~(ii) and~(iii) of the claim
the three line theorem applies and gives
$$
   \abs{F(\tfrac12)}\le \sqrt{M_0 M_1} .
$$
By part~(i) of the claim 
$$
   \abs{\INNER{T\xi}{\eta}_{\sqrt{\gf{g}}}}\le \sqrt{M_0 M_1} .
$$
Since $\xi$ and $\eta$ were arbitrary unit vectors in
$\ell^2_{\sqrt{\gf{f}}}$, respectively $\ell^2_{\sqrt{\gf{g}}}$,
this implies
$$
   \Norm{T}_{\Ll(\ell^2_{\sqrt{\gf{f}}},
     \ell^2_{\sqrt{\gf{g}}})}^2
   :=\sup_{\norm{\xi}_{\ell^2_{\sqrt{\gf{f}}}}=1}\Norm{T\xi}_{\ell^2_{\sqrt{\gf{g}}}}
   =\sup_{\norm{\xi}_{\ell^2_{\sqrt{\gf{f}}}}=1\atop\norm{\eta}_{\ell^2_{\sqrt{\gf{g}}}}=1}
   \INNER{T\xi}{\eta}_{\sqrt{\gf{g}}}
   \le \sqrt{M_0 M_1} .
$$
Here the equality holds since ``$\ge$'' is true by Cauchy-Schwarz,
but ``$=$'' holds by choosing
$\eta=\norm{T\xi}_{\ell^2_{\sqrt{\gf{g}}}}^{-1} T\xi$.
This concludes the proof of
Proposition~\ref{prop:Stein}.
\end{proof}

\bibliographystyle{alpha}
\addcontentsline{toc}{section}{References}
\bibliography{$HOME/Dropbox/0-Libraries+app-data/Bibdesk-BibFiles/library_math,$HOME/Dropbox/0-Libraries+app-data/Bibdesk-BibFiles/library_math_2020,$HOME/Dropbox/0-Libraries+app-data/Bibdesk-BibFiles/library_physics}{}

%


\end{document}